\documentclass[12pt]{article}
\title{Center stable manifolds around line solitary waves of the Zakharov--Kuznetsov equation
\amssubj{35B35, 35Q53.}
}
\author{YOHEI YAMAZAKI\footnote{{\it E-mail addresses:} yohei-yamazaki@hiroshima-u.ac.jp} \\ {\footnotesize Department of Mathematics, Kyoto University,} \\ {\footnotesize Kitashirakawa-Oiwakecho, Sakyo, Kyoto 606-8502, Japan,}\\ {\footnotesize Osaka City University Advanced Mathematical Institute, }\\ {\footnotesize 3-3-138 Sugimoto, Sumiyoshi-ku Osaka 558-8585, Japan, }  \\ {\footnotesize University of Cergy-Pontoise,
 } \\  {\footnotesize  UMR CNRS 8088, Cergy-Pontoise, F-95000, France} \\
{\footnotesize and } \\ {\footnotesize Hirhoshima University,
 } \\  {\footnotesize  1-3-2 Kagamiyama, Higashi-Hiroshima City, Hiroshima, 739-8511 Japan} }
\usepackage{amssymb}
\usepackage{mathrsfs}
\usepackage{amsthm}
\usepackage{amsmath}
\usepackage{titlefoot}
\usepackage{bm}
\def\pdfliteral #1 {}

\numberwithin{equation}{section}

\newtheorem{theorem}{Theorem}[section]
\newtheorem{corollary}[theorem]{Corollary}

\newtheorem{lemma}[theorem]{Lemma}
\newtheorem{proposition}[theorem]{Proposition}

\theoremstyle{definition}
\newtheorem{definition}[theorem]{Definition}
\newtheorem{remark}[theorem]{Remark}

\renewcommand{\eqref}[1]{(\ref{#1})}

\renewcommand{\bigskip}{\vspace{0.3cm}}

\newcommand{\N}{{\mathbb N}}
\newcommand{\R}{{\mathbb R}}

\newcommand{\Z}{{\mathbb Z}}
\newcommand{\T}{{\mathbb T}}

\newcommand{\norm}[1]{{\left \lVert #1 \right \rVert}}

\newcommand{\rbr}[1]{\left( #1 \right)}
\newcommand{\tbr}[1]{\langle #1 \rangle}
\newcommand{\gbr}[1]{\lceil #1 \rfloor}
\newcommand{\Tbr}[1]{\left\langle #1 \right\rangle}

\newcommand{\RTL}{\mathbb {R} \times \mathbb {T}_L}

\date{}

\begin{document}

\maketitle


\begin{abstract}
In this paper, we construct  center stable manifolds of  unstable line solitary waves for the Zakharov--Kuznetsov equation on $\RTL$ and show the orbital stability of the unstable line solitary waves on the center stable manifolds, which yields the asymptotic stability of unstable solitary waves on the center stable manifolds near by stable line solitary waves.
The construction is based on the graph transform approach by  Nakanishi--Schlag \cite{N S 2}.
Applying the bilinear estimate on Fourier restriction spaces by Molinet--Pilod \cite{M P} and modifying the mobile distance in \cite{N S 2}, we construct a contraction map on the graph space.
\end{abstract}


\section{Introduction}

We consider the two dimensional Zakharov--Kuznetsov equation
\begin{equation}\label{ZKeq}
 u_t + \partial_x(\Delta u + u^2)=0, \quad (t,x,y) \in \R\times \RTL,
 \end{equation}
where $\Delta=\partial_x^2 + \partial_y^2$, $u=u(t,x,y)$ is an unknown real-valued function, $\T_L=\R/2\pi L \Z$ and $L>0$.
The equation \eqref{ZKeq} preserves the mass and the energy:
\begin{align*}
M(u) = \int_{\RTL} |u|^2 dxdy,
\end{align*}
\begin{align*}
E(u)= \int_{\RTL}\Bigl( \frac{1}{2}|\nabla u|^2 - \frac{1}{3}u^3 \Bigr) dxdy,
\end{align*}
where $u \in H^1(\RTL)$.

The Zakharov--Kuznetsov equation was derived by Zakharov and Kuznetsov \cite{Z K} to describe the propagation of ionic-acoustic waves in uniformly magnetized plasma.
The rigorous derivation of the Zakharov--Kuznetsov equation from the Euler--Poisson system was proved by Lannes, Linares and Saut \cite{L L S}.
The Cauchy problem of the Zakharov--Kuznetsov equation has been studied in many paper \cite{AG1,AVF,G H,L P 1,L P 2,L P S,L S,M P,R V}.
The global well-posedness of the Zakharov--Kuznetsov equation in $H^s(\RTL)$ for $s>\frac{3}{2}$ has been proved by Linares, Pastor and Saut \cite{L P S} to study of the transverse instability of the $N$-soliton of the Korteweg-de Vries equation.
By proving a bilinear estimate in the context of Bourgain's spaces $X^{s,b}$, Molinet and Pliod \cite{M P} improved the result of the well-poseness on $\RTL$ in \cite{L P S} and showed the global well-posedness in $H^1(\RTL)$.


The Zakharov--Kuznetsov equation is one of multi-dimensional model of the Korteweg--de Vries equation.
The Korteweg--de Vries equation has the one soliton 
\[ Q_c(x-ct)=\frac{3c}{2}\cosh^{-2}\Bigl( \frac{\sqrt{c}(x-ct)}{2} \Bigr),\]
where $c>0$.
The orbital stability of the one soliton was showed by Benjamin \cite{TBB}.
Pego and Weinstein \cite{P W 2} proved the asymptotic stability of the one soliton on the exponentially weighted space.
To treat solutions including a small soliton, Mizumachi \cite{TM1} showed the the asymptotic stability of the one soliton on a polynomial weighed spaces.
In \cite{M M 1,M M 2,M M 3}, Martel and Merle proved the asymptotic stability on the energy space by using the Liouville type theorem and the monotonicity property.

We regard the one soliton $Q_c(x-ct)$ of the Korteweg--de Vries equation as a line solitary wave of \eqref{ZKeq}.
In \cite{K P}, to study the stability of line solitary wave with weak transverse perturbation of the KdV flow, Kadomtsev and Petviashvili derived the two-dimensional models of the KdV equation which is the KP equation.
In \cite{VEZ}, Zakharov obtained the proof of the instability of line solitary waves on the KP-I flow which was based on the existence
of a Lax pair for the KP-I equation.
The spectral stability of line solitary waves as the KP equation was obtained by Alexander, Pego and Sachs \cite{A P S}.
In \cite{R T 0,R T 1,R T 3}, by using the argument which is applicable to show  the transverse instability of the dispersive equations without integrable structure, Rousset and Tzvetkov proved the stability and instability of line solitary waves of the KP-I equation on $\R^2$ and $\RTL$.
Applying the inverse scattering method, Villarroel and Ablowitz showed the stability of line solitary waves of the KP-II equation for the decaying perturbations in \cite{V A}. 
The orbital stability and the asymptotic stability of line solitary waves for the KP-II equation on $\RTL$ was showed by Mizumachi and Tzvetkov \cite{M T}.
Using the local modulations of the amplitude and the phase shift of line solitary waves which behaves like a self-similar solution of the Burgers equation, Mizumachi proved the asymptotic stability of line solitary waves of the KP-II equation on $\R^2$ in \cite{TM2,TM3}.  
The instability of the line solitary waves of the Zakharov--Kuznetsov equation on $\R^2$ was proved by Rousset and Tzvetkov \cite{R T 0}.
On $\T_{L_1}\times \T_{L_2}$ with sufficiently large $L_2$, the linear instability of line periodic solitary waves of the Zakharov--Kuznetsov equation have been showed by Johnson \cite{MJ} by using Evan's function method.
The instability of the line solitary waves of the Zakharov--Kuznetsov equation on $\RTL$ with large traveling speed was showed by Bridges \cite{TJB}.
In \cite{YY4}, the author proved that the line solitary waves $Q_c(x-ct)$ of the Zakharov--Kuznetsov equation on $\RTL$ is orbitally stable and the asymptotically stable for $0<c\leq \frac{4}{5L^2}$ and is unstable for $c>\frac{4}{5L^2}$.
The proof of the asymptotic stability in \cite{YY4} is based on a Liouville type theorem and virial type estimates in \cite{M M 3}.
By the Lyapunov--Schmidt reduction method, the existence and the stability of the transversely modulated solitary waves was showed in \cite{YY4}.
Using the normal form which describes the motion of the amplitude of the transversely modulated solitary waves, Pelinovsky proved the asymptotic behavior of solutions near by the transversely modulated solitary waves for the Zakharov--Kuznetsov equation in \cite{DP}.
Moreover, Pelinovsky showed the asymptotic stability of the transversely modulated solitary waves in the sense by Pego and Winstein \cite{P W 2}.

In this paper, we construct a center stable manifold in energy space to study the behavior of solutions near by unstable line solitary waves.
There have been many papers \cite{B J,MB0,MB,J L Z,K S,K N S,M M N R,N S,N S 2,WS} for constructing the center stable manifold for various equations.
In \cite{B J}, developing the Hadamard method, Bates and Jones and constructed invariant manifolds in abstract setting  for nonlinear partial differential equations.
Moreover, applying the construction in abstract setting, Bates and Jones proved the existence of a Lipschitz center stable manifold of the nonlinear Klein--Gordon equation under the radial symmetry restriction with the power nonlinearity which satisfies the nonlinearity is local Lipschitz $H^1 \to L^2$.
To treat a derivative loss term due to the translation, Nakanishi and Schlag \cite{N S 2} proved the existence of a center stable manifold for the nonlinear Klein--Gordon equation around solitary waves which generated from radial stationary solutions by the action of Lorentz transforms and spatial translations  by introducing the mobile distance.
Using a framework based on vector bundle coordinates, Jin, Lin and Zeng \cite{J L Z} constructed the center stable manifold for the 3D Cross-Pitaecskii equation around solitary waves.
By using the Strichartz estimate of the linear evolution around ground states, Schlag constructed a center stable manifold for the 3D cubic nonlinear Schr\"odinger equation around ground states in $W^{1,1}\cup W^{1,2}$ and proved the scattering on the center stable manifold in \cite{WS}.
Improving the result \cite{WS}, Beceanu \cite{MB} constructed a center stable manifold for the 3D cubic nonlinear Schr\"odinger equation around ground states in the critical space $\dot{H}^{\frac{1}{2}}$.
Applying the argument \cite{WS}, Krieger and Schlag \cite{K S} constructed a center stable manifold for 1D nonlinear Schr\"odinger equation with $L^2$-critical nonlinearity around ground states.
By investigating the ejection of solutions near the ground state and using trichotomy results in \cite{K N S 0}, Krieger, Nakanishi and Schlag construct a center stable manifold for the energy critical nonlinear wave equation.
Martel, Merle, Nakanishi and Rapha\"el \cite{M M N R} constructed a center stable manifold for the $L^2$-critical generalized KdV equation around one solition on weighted space by applying trichotomy result which classifies initial datum near one soliton by the asymptotic behavior of solution.

To state the main result, we define some notations.
The solitary wave manifold of $Q_c$ is defined as
\[ S(c)=\{\tau_{q}Q_c: q \in \R\}\]
and the neighborhood of the solitary wave manifold $S(c)$ is defined as
\[N_{c}(\delta)=\{ u \in H^1(\RTL): \inf_{ q \in \R}\norm{u-\tau_qQ_c}_{H^1}<\delta\},\]
where $(\tau_qu)(x,y)=u(x-q,y)$.
Let $\mathbb{L}_c=-\Delta +c -2 Q_c$.
Then, the linearized operator of \eqref{ZKeq} around the line solitary wave $Q_c(x-ct)$ as a relative equilibrium point is $\partial_x \mathbb{L}_c$.
By the global well-posedness result in \cite{M P}, we define $U(t)$ as the flow map of \eqref{ZKeq} at time $t$.
The following theorem is the main theorem.
\begin{theorem}\label{thm-main}
Let $c_0 \in \{ c >4/5L^2: c \neq 4n^2/5L^2 \mbox{ for } n \in \Z\}$.
Then there exists $C^1$ manifold $\mathcal{M}_{cs}(c_0)$ in $H^1(\R)$ containing the solitary wave manifold $S(c_0)$ with the following properties:
\begin{enumerate}
  \setlength{\parskip}{0.1cm} 
  \setlength{\itemsep}{0.05cm} 
\renewcommand{\labelenumi}{\rm (\roman{enumi})}
\item The codimension of $\mathcal{M}_{cs}(c_0)$ in $H^1(\RTL)$ equals $2$ times the integer part of $\frac{\sqrt{5c_0}L}{2}$ which is the total dimension of the eigenspaces of the linearized operator $\partial_x\mathbb{L}_{c_0}$ corresponding to eigenvalues with positive real part.
\item $\tau_q U(t) \mathcal{M}_{cs}(c_0) \subset \mathcal{M}_{cs}(c_0)$ for $q \in \R $ and $t \geq 0$.
\item $\mathcal{M}_{cs}(c_0)$ is normal at $Q_{c_0}$ to the eigenspaces corresponding to eigenvalues of $\partial_x\mathbb{L}_{c_0}$ with positive real part.
\item For any $\varepsilon >0$, there exists $\delta>0$ such that $U(t)(\mathcal{M}_{cs}(c_0)\cap N_{c_0}(\delta)) \subset \mathcal{M}_{cs}(c_0) \cap N_{c_0}(\varepsilon)$ for $t \geq 0$.
\item There is $\varepsilon_0>0$ such that for $u_0 \in N_{c_0}(\varepsilon_0) \setminus \mathcal{M}_{cs}$ there exists $t_0>0$ satisfying $ U(t_0)u_0 \notin N_{c_0}(\varepsilon_0)$.
\end{enumerate}
\end{theorem}
\begin{remark}
The construction of the center stable manifold $\mathcal{M}_{cs}(c)$ is based on the global well-posedness of \eqref{ZKeq} in \cite{M P}.
Therefore, in this paper, we only consider the solutions of \eqref{ZKeq} in \cite{M P} which are in the Bourgain spaces in local time.
If the unconditional uniqueness of the solutions of \eqref{ZKeq} in $C(\R,H^1(\RTL))$ is proved, we can show the center stable manifold without the restriction of the class of solutions.
\end{remark}
\begin{remark}
For any positive integer $n$, $\partial_x\mathbb{L}_{4n^2/5L^2}$ has extra eigenfunctions corresponding to $0$ eigenvalue which follows a bifurcation of the branch $\{Q_c:c>0\}$ at $c=4n^2/5L^2$.
Therefore, it is difficult to construct a manifold around $S(4n^2/5L^2)$ satisfying (i)--(v) in Theorem \ref{thm-main}.    
\end{remark}
Applying the asymptotic stability result of the line solitary wave with the critical speed $4/5L^2$, we obtain the asymptotic behavior of solutions on the center stable manifold  $\mathcal{M}_{cs}(c)$ near by $Q_{4/5L^2}$.
\begin{corollary}\label{cor-main}
For any $\beta>0$ there exists $c_{\beta}>4/5L^2$ such that for $4/5L^2<c<c_{\beta}$ there exists $\varepsilon_{\beta,c}>0$ satisfying the following.
For any  solution $u_0 \in \mathcal{M}_{sc}(c)\cap N_{c}(\varepsilon_{\beta,c})$ there exist $\rho_1 \in C^1([0,\infty))$ and $c_+>0$ such that 
\[\lim_{ t \to \infty} \int_{\{x>\beta t\}}\Bigl(U(t)u_0-\tau_{\rho_1(t)}Q_{c_+}\Bigr) \, dxdy=0, \quad \lim_{t \to \infty} \dot{\rho}_1(t)=c_+,\]
\[|c-c_+|\lesssim \norm{u_0-Q_c}_{H^1(\RTL)}.\]
\end{corollary}

The proof of the existence of the center stable manifold $\mathcal{M}_{cs}(c)$ is based on the argument by \cite{N S 2}.
Since the translation of functions $\tau_q u$ in the energy space are not Lipschitz continuous generally by the energy norm with respect to the translation parameter $q$, it is difficult to show the smallness of the difference between functions and the translated function uniformly.
In \cite{N S 2}, Nakanishi and Schlag introduced the mobile distance on the energy space which measures the translation of functions as a Lipschitz continuous term.

There are two difficulties to apply the argument by \cite{N S 2}.
In the case of \eqref{ZKeq}, the generalized eigenfunction of the adjoint operator of the linearized operator of \eqref{ZKeq} around the line solitary wave is not in $L^2(\RTL)$.
Therefore, we can not use the suitable symplectic spectral decomposition of functions on the energy space.
By applying anther decomposition, a modulation term has same order of the difference between a line solitary wave and the solution of \eqref{ZKeq}.
In \cite{N S 2}, Nakanishi and Schlag used that the order of modulation term is higher than the order of the difference to show the estimate for the contraction map on the set of graphs which is in Lemma 3.2 of \cite{N S 2}.
Therefore, we can not show the estimate for the contraction map by the mobile distance in \cite{N S 2}.
This difficulty appeared to construct a center stable manifold of the energy critical wave equation and was overcame in \cite{K N S} by the ignition lemma and that the codimension of the center stable manifold is one.
To construct the center stable manifold with high codimension and to modify the argument in \cite{N S 2}, we adjust the scaling of the correction term in the mobile distance and show the estimate for the contraction in Lemma \ref{lem-est-md}.
The equation \eqref{ZKeq} has the nonlinear term $\partial_x u^2$ which has a derivative loss.
To control the nonlinear team, in \cite{M P}, Molinet and Pliod proved a bilinear estimate on Fourier restriction spaces introduced by Bourgain \cite{JB}.
To treat the nonlinear term $\partial_x u^2$ and to construct the center stable manifold, we apply the argument in \cite{N S 2} by using space-time estimates with space-time derivatives.

Our plan of the present paper is as follows.
In Section 2, we define a spectral decomposition with respect to the linearized operator of \eqref{ZKeq}.
In Section 3, we show the estimate of the difference between solutions of a linearized equation of \eqref{ZKeq} and solutions of a localized equation of \eqref{ZKeq} by the mobile distance.
In Section 4, we construct the center stable manifold by applying the argument in \cite{N S 2}.
In Section 5, we prove the $C^1$ regularity of the center stable manifold which follows the argument in \cite{K N S}.

\section{Preliminaries}

Let $c^*>0$.
In this section, we assume that there exists a positive integer $n_0$ such that $\frac{2n_0}{\sqrt{5c^*}}<L<\frac{2(n_0+1)}{\sqrt{5c^*}}$.
We denote the inner product in $L^2(X)$ by 
\[(u,v)_{L^2(X)}=\int_X uv\,dx, \quad u,v \in L^2(X),\]
and the coupling between $H^{1}(X)$ and $H^{-1}(X)$ by
\[\tbr{u,v}=\tbr{u,v}_{H^1(X), H^{-1}(X)}, \quad u \in H^1(X), v \in H^{-1}(X),\]
where $X=\R$ or $\RTL$.
In particular, we denote $(\cdot,\cdot)_{L^2}=(\cdot,\cdot)_{L^2(\RTL)}$, $\norm{\cdot}_{L^2}=\norm{\cdot}_{L^2(\RTL)}$ and $\norm{\cdot}_{H^s}=\norm{\cdot}_{H^s(\RTL)}$.

We define the linearized operator $\mathbb{L}_c$ of the stationary equation of \eqref{ZKeq} around $Q_c$ as
\begin{align*}
\mathbb{L}_c=-\Delta +c -2Q_c
\end{align*}
and the linearized operator $\mathcal{L}_c$ of the stationary equation of the Korteweg--de Vries equation around $Q_c$ as
\begin{align*}
\mathcal{L}_c=-\partial_x^2 +c -2 Q_c.
\end{align*}
Then, the linearized operator of \eqref{ZKeq} around $Q_c(x-ct)$ is $\partial_x\mathbb{L}_c$.
By the Fourier expansion, we have for $u \in H^1(\RTL)$ 
\[(\mathbb{L}_cu)(x,y)=\sum_{n=-\infty}^{\infty}\Bigl(\mathcal{L}_c+\frac{n^2}{L^2} \Bigr)u_n(x) e^{\frac{iny}{L}},\]
where
\[u(x,y)=\sum_{n=-\infty}^{\infty} u_n(x)e^{\frac{iny}{L}}.\]

The following proposition follows in Proposition 3.1 in \cite{YY3}.
\begin{proposition}\label{prop-2-1}
Let $c>0$.
The following holds.
\begin{enumerate}
  \setlength{\parskip}{0.1cm} 
  \setlength{\itemsep}{0.05cm} 
\renewcommand{\labelenumi}{\rm (\roman{enumi})}
\item $\mathcal{L}_c$ has the only one negative eigenvalue $-\frac{5c}{4}$. Moreover, the set of eigenfunctions of $\mathcal{L}_c$ corresponding to $-\frac{5c}{4}$ is spanned by $Q_c^{\frac{3}{2}}$.
\item If $0<a< \frac{5c}{4}$, then $\partial_x(\mathcal{L}_c+a)$ has the only one positive eigenvalue $\lambda(a)$ $($resp. the only one negative eigenvalue $-\lambda(a)$ $)$ which is simple. Moreover, eigenfunctions of $\partial_x(\mathcal{L}_c+a)$ corresponding to $\lambda(a)$ $($resp. $-\lambda(a)$ $)$ is in $H^{\infty}(\R)$.
\item If $a>\frac{5c}{4}$, then $\partial_x(\mathcal{L}_c+a)$ has no positive eigenvalues.
\item If $a\neq \frac{5c}{4}$ and $a>0$, then the kernel of $\partial_x(\mathcal{L}_c+a)$ is trivial.
\item If $a=\frac{5c}{4}$, then the kernel of $\partial_x(\mathcal{L}_c+a)$ is spanned by $Q_c^{\frac{3}{2}}$.
\end{enumerate}
\end{proposition}
Let $\lambda_k=\lambda(\frac{k^2}{L^2})$ 
for integer $k$ with $0<k\leq n_0$.
We define 
\begin{align}\label{min-max-eigen}
k^*=\max_{0<k\leq n_0}\lambda_k , \quad k_*=\min_{0<k\leq n_0} \lambda_{k}.
\end{align}
For $0<k\leq n_0$, we define the eigenfunction of $\partial_x(\mathcal{L}_{c^*}+\frac{k^2}{L^2})$ corresponding to $ \lambda_k$ (resp. $-\lambda_k$) as $f_{k}^+$ (resp. $f_{k}^-$).
Then, we show the following property of $f_k^+$.
\begin{proposition}
Let $g_k^+(x)=f_k^+(-x)$.
The following holds.
\[\partial_x\Bigl(\mathcal{L}_{c^*}+\frac{k^2}{L^2}\Bigr)g_k^+=-\lambda_kg_k^+, \quad \Bigl(f_k^+,\Bigl(\mathcal{L}_{c^*}+\frac{k^2}{L^2}\Bigr)g_k^+\Bigr)_{L^2(\R)}<0.\]
\end{proposition}
\proof
$\partial_x(\mathcal{L}_{c^*}+\frac{k^2}{L^2})g_k^+=-\lambda_kg_k^+$ follows the definition of $g_k^+$.
By the spectral decomposition corresponding to $\mathcal{L}_{c^*}$, we denote $f_k^+$ by
\[f_k^+=cQ_{c^*}^{\frac{3}{2}}+\mu \partial_xQ_{c^*}+\gamma^+.\]
Then, the above equation yields
\[g_k^+=cQ_{c^*}^{\frac{3}{2}}-\mu \partial_xQ_{c^*}+\gamma^-,\]
where $\gamma^-(x)=\gamma^+(-x)$.
Let $\mathcal{L}_{c^*}(k)=\mathcal{L}_{c^*}+\frac{k^2}{L^2}$.
Since
\begin{align*}
(f_k^+,\mathcal{L}_{c^*}(k)f_k^+)_{L^2(\R)}=& \frac{1}{\lambda_k}(\partial_x\mathcal{L}_{c^*}(k)f_k^+,\mathcal{L}_{c^*}(k)f_k^+)_{L^2(\R)},
\end{align*}
we have
\[(f_k^+,\mathcal{L}_{c^*}(k)f_k^+)_{L^2(\R)}=(g_k^+,\mathcal{L}_{c^*}(k)g_k^+)_{L^2(\R)}=0.\]
On the other hand, from the equations
\begin{align*}
(f_k^+,\mathcal{L}_{c^*}(k)f_k^+)_{L^2(\R)}=\Bigl(\frac{k^2}{L^2}-\frac{5c^*}{4}  \Bigr) c^2\norm{Q_{c^*}^{\frac{3}{2}}}_{L^2(\R)}^2+\frac{k^2\mu^2}{L^2}\norm{\partial_xQ_{c^*}}_{L^2(\R)}^2+(\gamma^+,\mathcal{L}_{c^*}(k)\gamma^+)_{L^2(\R)}
\end{align*}
and
\begin{align*}
(g_k^+,\mathcal{L}_{c^*}(k)g_k^+)_{L^2(\R)}=\Bigl(\frac{k^2}{L^2}-\frac{5c^*}{4}  \Bigr) c^2\norm{Q_{c^*}^{\frac{3}{2}}}_{L^2(\R)}^2+\frac{k^2\mu^2}{L^2}\norm{\partial_xQ_{c^*}}_{L^2(\R)}^2+(\gamma^-,\mathcal{L}_{c^*}(k)\gamma^-)_{L^2(\R)},
\end{align*}
we obtain
\begin{align*}
(g_k^+,\mathcal{L}_{c^*}(k)f_k^+)_{L^2(\R)}=&\Bigl(\frac{k^2}{L^2}-\frac{5c^*}{4}  \Bigr) c^2\norm{Q_{c^*}^{\frac{3}{2}}}_{L^2(\R)}^2-\frac{k^2\mu^2}{L^2}\norm{\partial_xQ_{c^*}}_{L^2(\R)}^2+(\gamma^-,\mathcal{L}_{c^*}(k)\gamma^+)_{L^2(\R)}\\
\leq & -2\frac{k^2}{L^2}\mu^2\norm{\partial_xQ_{c^*}}_{L^2(\R)}^2.
\end{align*}
We assume $\mu=0$ and 
\[(\gamma^-,\mathcal{L}_{c^*}(k)\gamma^+)_{L^2(\R)}=(\gamma^+,\mathcal{L}_{c^*}(k)\gamma^+)_{L^2(\R)},\]
then $\gamma_+=\gamma_-$ and $g_k^+=f_k^+$.
This contradicts $\lambda_k>0$.
Thus, we obtain the conclusion.
\qed

After an appropriate normalization of $f_{k}^+$ and $f_{k}^-$, we have
\[ \Bigl( f_{k}^+, \Bigl(\mathcal{L}_{c^*}+\frac{k^2}{L^2}\Bigr) f_{k}^-\Bigr)_{L^2(\R)}=\frac{1}{\pi L}, \quad \norm{f_{k}^+}_{L^2(\R)}=\norm{f_{k}^-}_{L^2(\R)}.\]
Then, from the uniqueness of $f_{k}^+$ and $f_{k}^-$, we obtain that 
\[f_{k}^+(x)=-f_{k}^-(-x), \quad x\in \R.\]
The unstable and stable eigenfunctions of $\partial_x\mathbb{L}_{c^*}$ are denoted by
\[F_k^{\pm,0}(x,y)=f_k^{\pm}(x) \cos \frac{ky}{L}, \quad  F_k^{\pm,1}(x,y)=f_k^{\pm}(x) \sin \frac{ky}{L}.\]
The functions $F_k^{\pm,0}$ and $F_k^{\pm,1}$ satisfy
\[(F_k^{+,j},\mathbb{L}_{c^*}F_k^{-,j})_{L^2}=1, \quad \partial_x\mathbb{L}_{c^*}F_k^{\pm,j}=\pm \lambda_k F_k^{\pm,j}, \quad j=0,1.\] 
We consider the decomposition in $L^2(\RTL)$ such that
\begin{align}\label{2-decom}
u=\sum_{\substack{ j=0,1 \\ k=1,2, \dots, n_0}} (\Lambda_k^{+,j} F_{k}^{+,j} + \Lambda_k^{-,j} F_{k}^{-,j}) + \mu_1 \partial_x Q_{c^*} + \mu_2 \partial_c Q_{c^*}+\gamma,
\end{align}
where 
\[ \Lambda_k^{\pm ,j}= \Lambda_k^{\pm ,j}(u)=(u, \mathbb{L}_{c^*} F_k^{\mp,j})_{L^2}, \quad \mu_1=\mu_1(u)=\frac{(u,\partial_x Q_{c^*})_{L^2}}{\norm{\partial_x Q_{c^*}}_{L^2}^2}, \]
\[\mu_2=\mu_2(u)=\frac{(u,Q_{c^*})_{L^2}}{(\partial_cQ_{c^*},Q_{c^*})_{L^2}},\]
\[\gamma=\gamma(u)=u- \sum_{\substack{ j=0,1 \\ k=1,2, \dots, n_0}} (\Lambda_k^{+,j} F_{k}^{+,j} + \Lambda_k^{-,j} F_{k}^{-,j}) - \mu_1 \partial_x Q_{c^*} - \mu_2 \partial_c Q_{c^*}.\]
We define the projections corresponding to \eqref{2-decom} as
\[P_{\pm}u=\sum_{\substack{ j=0,1 \\ k=1,2, \dots, n_0}}\Lambda_k^{\pm,j}(u) F_k^{\pm,j}, \quad P_0 u=\mu_1(u)\partial_x Q_{c^*}+\mu_2(u)\partial_cQ_{c^*}, \quad P_1u=\mu_1(u)\partial_xQ_{c^*}, \]
\[ P_2u=\mu_2(u)\partial_cQ_{c^*}, \quad P_{\gamma}u=\gamma(u), \quad P_d=Id-P_{\gamma}.\]
The orthogonality of the projections yields the properties:
\begin{align}
P_{\pm}\partial_x\mathbb{L}_{c^*}=\partial_x\mathbb{L}_{c^*}P_{\pm},\label{orth-pro-1}
\\
\mathbb{L}_{c^*}P_1=0, \quad P_2\partial_x\mathbb{L}_{c^*}=0 \mbox{ and }\mathbb{L}_{c^*}\partial_x\mathbb{L}_{c^*}P_2=0. \label{orth-pro-2}
\end{align}
The following proposition follows the proof of Lemma 2.2 in \cite{N S}.
\begin{proposition}\label{prop-energy-norm}
There exists $C>0$ such that for $u \in H^1(\RTL)$ we have
\[\tbr{P_{\gamma}u,\mathbb{L}_{c^*}P_{\gamma} u}_{H^1,H^{-1}} \geq C \norm{P_{\gamma}u}_{H^1}^2,\]
 where $\tbr{\cdot,\cdot}_{H^1,H^{-1}}$ is the coupling between $H^1(\RTL)$ and $H^{-1}(\RTL)$.
\end{proposition}
\proof
We assume that there exists $P_{\gamma}u_0 \neq 0$ satisfying $\tbr{P_{\gamma}u_0,\mathbb{L}_{c^*}P_{\gamma}u_0}_{H^1,H^{-1}}\leq 0$.
Then, for $0<k\leq n_0$ and $j=0,1$ we have 
\[\tbr{\mathbb{L}_{c^*}P_{\gamma}u_0,\partial_xQ_{c^*}}_{H^{-1},H^1}=\tbr{\mathbb{L}_{c^*}P_{\gamma}u_0, F_k^{+,j}}_{H^{-1},H^1}=0.\]
Since for $0<k_1,k_2\leq n_0$ and $ j_1,j_2\in \{0,1\}$ 
\[(\mathbb{L}_{c^*}\partial_xQ_{c^*},\partial_xQ_{c^*})_{L^2}=(\mathbb{L}_{c^*}F_{k_1}^{+,j_1},F_{k_2}^{+,j_2})_{L^2}=0,\]
 we obtain
\[\Tbr{\mathbb{L}_{c^*}\Bigl(aP_{\gamma}u_0+a_0 \partial_x Q_{c^*}+\sum_{j,k} a_{k,j}F_k^{+,j}\Bigr), aP_{\gamma}u_0+a_0\partial_xQ_{c^*}+\sum_{j,k}a_{k,j}F_k^{+,j}}_{H^1,H^{-1}}\leq 0\]
for $a, a_0, a_{k,j} \in \R$.
Therefore, the dimension of the non-positive eigenspace of $\mathbb{L}_{c^*} $ is more than $2n_0+1$,  which contradicts the non-positive eigenspace of $\mathbb{L}_{c^*} $ is spanned by 
\[\Bigl\{Q_{c^*}^{\frac{3}{2}} \cos \frac{ky}{L}, Q_{c^*}^{\frac{3}{2}} \sin \frac{ky}{L},\partial_xQ_{c^*}; 0<k\leq n_0, k \in \Z \Bigr\}.\]
Therefore, for any $u \in H^1(\RTL)$ we have $\tbr{P_{\gamma}u_0,\mathbb{L}_{c^*}P_{\gamma}u_0}_{H^1,H^{-1}}\geq 0$.
By Weyl's theorem on essential spectrum, the essential spectrum of $P_{\gamma}^*\mathbb{L}_{c^*}P_{\gamma}$ is $[c^*,\infty)$, where $P_{\gamma}^*$ is the adjoint operator of $P_{\gamma}$.
Thus, we obtain the conclusion.
\qed

The linearized energy norm is defined on $H^1(\RTL)$ by
\begin{align*}
\norm{u}_{E}^2=\sum_{\substack{ j=0,1 \\ k=1,2, \dots, n_0}} \bigl((\Lambda_k^{+,j})^2+(\Lambda_k^{-,j})^2\bigr)+\mu_1^2+\mu_2^2+\tbr{\gamma,\mathbb{L}_{c^*}\gamma}_{H^1,H^{-1}}.
\end{align*}
Then, the linearized energy norm $\norm{\cdot}_{E}$ is equivalent to the energy norm $\norm{\cdot}_{H^1}$.

Let $u$ be a solution to \eqref{ZKeq} and $\tau_{\rho}u(t,x,y)=u(t,x-\rho,y)$ for $\rho \in \R$.
Then, $v(t)=\tau_{-\rho(t)} u(t)-Q_{c(t)}$ solves 
\begin{align}\label{vZKeq}
v_t=\partial_x\mathbb{L}_{c^*}v+(\dot{\rho}-c)\partial_xQ_{c^*}-\dot{c}\partial_cQ_{c^*} + N(v,c,\rho),
\end{align}
where
\[N(v,c,\rho)=\partial_x[-v^2+(\dot{\rho}-c^*)v+2(Q_{c^*}-Q_c)v+(\dot{\rho}-c)(Q_{c}-Q_{c^*})] -\dot{c}\partial_c(Q_{c}-Q_{c^*}).\]
Then, using the following lemma, we choose $c$ and $\rho$ which satisfy the orthogonality condition
\begin{align}\label{cond-orth}
(v,\partial_xQ_{c^*})_{L^2}=(v,Q_{c^*})_{L^2}=0
\end{align}
for $u \in \mathcal{N}_{\delta,c^*}$, where $v=\tau_{-\rho}u-Q_{c}$ and 
\[\mathcal{N}_{\delta,c^*} =\{u \in H^1(\RTL); \inf_{\rho \in \R} \norm{u-\tau_{\rho}Q_{c^*}}_{H^1}<\delta\}. \]
\begin{lemma}\label{lem-orth}
There exist $\delta_0, C_{\delta_0}>0$ and smooth maps $\rho : \mathcal{N}_{\delta_0, c^*} \to \R$ and $c:\mathcal{N}_{\delta_0,c^*} \to (0,\infty)$  such that for $u \in \mathcal{N}_{\delta_0,c^*}$, $v=\tau_{-\rho(u)}u-Q_{c(u)}$ satisfies the orthogonality condition $\eqref{cond-orth}$ and
\begin{align}\label{eq-lem-orth-1}
\norm{v}_{H^1}+|c(u)-c^*|< C_{\delta_0} \inf_{q\in \R} \norm{u-\tau_qQ_{c^*}}_{H^1}.
\end{align}
\end{lemma}
\proof
We define $G$ by 
\[G(u,c,\rho)=
\begin{pmatrix}
(\tau_{-\rho}u-Q_{c}, Q_{c^*})_{L^2} \\
(\tau_{-\rho}u-Q_{c}, \partial_xQ_{c^*})_{L^2}
\end{pmatrix}.\]
Then, $G(Q_{c^*},c^*,0)= {}^t(0,0)$ and 
\[\frac{\partial G}{\partial c \partial \rho}(Q_{c^*},c^*,0)=\mbox{diag}(-(\partial_{c}Q_{c^*},Q_{c^*})_{L^2}, -\norm{\partial_x Q_{c^*}}_{L^2}^2).\]
By $(\partial_{c}Q_{c^*},Q_{c^*})_{L^2}>0$ and the implicit function theorem, we obtain that there exists $\delta_0, C_{\delta_0}>0$ such that  there is a unique smooth map $(c(u), \rho(u))$ satisfying \eqref{cond-orth} and
\[\norm{\tau_{-\rho(u)}u-Q_{c(u)}}_{H^1}+|c(u)-c^*|<C_{\delta_0}\norm{u-Q_{c^*}}_{H^1}\]
for $u \in \{u \in H^1(\RTL); \norm{u-Q_{c^*}}_{H^1}<\delta_0 \}$.
By the uniqueness of the map $(c(u), \rho(u))$ and the invariance with respect to $\tau_\rho$, expanding the map $(c(u), \rho(u))$, we obtain the unique smooth map $(c(u), \rho(u))$ satisfying the orthogonality condition $\eqref{cond-orth}$ and \eqref{eq-lem-orth-1}.
\qed

Let $c(t)=c(u(t))$ and $\rho(t)=\rho(u(t))$, where $c(u(t))$ and $\rho(u(t))$ are defined in Lemma \ref{lem-orth}.
Then, from the orthogonality condition, we have 
\begin{align*}
0=\frac{d}{dt}(v,\partial_xQ_{c^*})_{L^2}
=& (\dot{\rho}-c) \bigl( \norm{\partial_x Q_{c^*}}_{L^2}^2 +(\partial_x v,\partial_xQ_{c^*})_{L^2}+(\partial_x(Q_c-Q_{c^*}),\partial_xQ_{c^*})_{L^2}\bigr )\\
&+ (c-c^*)(\partial_x v,\partial_xQ_{c^*})_{L^2} - (v,\mathbb{L}_{c^*}\partial_x^2Q_{c^*})_{L^2}  - (\partial_x v^2, \partial_x Q_{c^*})_{L^2}\\
&-(2(Q_{c^*}-Q_c)v,\partial_x^2Q_{c^*})_{L^2},
\end{align*}
\begin{align*}
0=\frac{d}{dt}(v,Q_{c^*})_{L^2}
=&-\dot{c} (\partial_cQ_{c}, Q_{c^*}) _{L^2}- (\partial_x v^2, Q_{c^*})_{L^2}-(2(Q_{c^*}-Q_c)v,\partial_x Q_{c^*})_{L^2}.
\end{align*}
Therefore, $(c(t),\rho(t))$ satisfies 
\begin{align}\label{eq-orth-1}
\begin{pmatrix}
\dot{\rho}-c\\
\dot{c}
\end{pmatrix}
=
\begin{pmatrix}
\norm{\partial_x Q_{c^*}}_{L^2}^{-2}(v, \mathbb{L}_{c^*}\partial_x^2Q_{c^*})_{L^2}\\
0
\end{pmatrix}
+ \bm{N}(v,c),
\end{align}
where
\begin{align*}
&\bm{N}(v,c)\\
=&
\begin{pmatrix}
\norm{\partial_xQ_{c^*}}_{L^2}^2 +(\partial_x v, \partial_xQ_{c^*})_{L^2} + (\partial_x(Q_c-Q_{c^*}),\partial_xQ_{c^*})_{L^2} & 0\\
0 & -(\partial_c Q_c, Q_{c^*})_{L^2}
\end{pmatrix}^{-1}\\
& \times
\begin{pmatrix}
((c-c^*)\partial_x v-\partial_x v^2+2\partial_x((Q_{c^*}-Q_c)v),\partial_xQ_{c^*})_{L^2}- (v,\mathbb{L}_{c^*}\partial_x^2Q_{c^*})_{L^2}  \\
- (\partial_x v^2, Q_{c^*})_{L^2}-(2(Q_{c^*}-Q_c)v,\partial_x Q_{c^*})_{L^2}
\end{pmatrix}\\
&-
\begin{pmatrix}
\norm{\partial_x Q_{c^*}}_{L^2}^{-2}(v, \mathbb{L}_{c^*}\partial_x^2Q_{c^*})_{L^2}\\
0
\end{pmatrix}\\
=& O(\norm{v}_{L^2}^2+\norm{v}_{L^2}|c-c^*|) \mbox{ as } \norm{v}_{L^2}+|c-c^*| \to 0.
\end{align*}
On the tubular neighborhood $\mathcal{N}_{\delta_0,c^*}$, $u=\tau_{\rho}(v+Q_{c})$ solves \eqref{ZKeq} with $(c(t),\rho(t)) $ satisfying the orthogonality condition \eqref{cond-orth} if and only if $v=\tau_{-\rho}u-Q_c$ solves \eqref{vZKeq} with $(c(t),\rho(t))$ satisfying \eqref{eq-orth-1} and $(v(0),\partial_xQ_{c^*})_{L^2}=(v(0),Q_{c^*})_{L^2}=0$.

\section{Localized equation}
Let $c^*>0$.
In this section, we assume that there exists a positive integer $n_0$ such that $\frac{2n_0}{\sqrt{5c^*}}<L<\frac{2(n_0+1)}{\sqrt{5c^*}}$.
Let $\chi \in C_0^{\infty}(\R)$ be a smooth function with
\[\chi(r)=\begin{cases} 1 & (|r|\leq 1) \\ 0 & (|r| \geq 2) \end{cases}, \quad 0 \leq \chi \leq 1.\]
Let 
\[\chi_{\delta}=\chi_{\delta}(v,c-c^*)=\chi\Bigl( \frac{\norm{v}_{H^1}^2+|c-c^*|^2}{\delta^2} \Bigr).\]
We define the localized system of \eqref{vZKeq} as
\begin{align}
v_t=&\partial_x \mathbb{L}_{c^*}v + (\dot{\rho}-c) \partial_x Q_{c^*}-\dot{c}\partial_cQ_{c^*} +\chi_{\delta}(v,c-c^*)N(v,c,\rho),\label{LZKeq-1}\\
\begin{pmatrix}
\dot{\rho}-c\\
\dot{c}
\end{pmatrix}
=&
\begin{pmatrix}
\norm{\partial_x Q_{c^*}}_{L^2}^{-2}(v, \mathbb{L}_{c^*}\partial_x^2Q_{c^*})_{L^2}\\
0
\end{pmatrix}
+ \bm{N}_{\delta}(v,c), \label{LZKeq-2}
\end{align}
where
\begin{align*}
&\bm{N}_{\delta}(v,c)\\
=&
\begin{pmatrix}
\norm{\partial_xQ_{c^*}}_{L^2}^2 +\chi_{\delta}(\partial_x v+\partial_x(Q_c-Q_{c^*}), \partial_xQ_{c^*})_{L^2}  & 0\\
0 & -(\partial_c Q_{c^*}+\chi_{\delta}\partial_c(Q_{c}-Q_{c^*}),Q_{c^*})_{L^2}
\end{pmatrix}^{-1}\\
& \times
\begin{pmatrix}
\chi_{\delta}((c-c^*)\partial_x v-\partial_x v^2+2\partial_x((Q_{c^*}-Q_c)v),\partial_xQ_{c^*})_{L^2}- (v,\mathbb{L}_{c^*}\partial_x^2Q_{c^*})_{L^2}  \\
-\chi_{\delta} (\partial_x v^2-2\partial_x((Q_{c^*}-Q_c)v), Q_{c^*})_{L^2}
\end{pmatrix}\\
&-
\begin{pmatrix}
\norm{\partial_x Q_{c^*}}_{L^2}^{-2}(v, \mathbb{L}_{c^*}\partial_x^2Q_{c^*})_{L^2}\\
0
\end{pmatrix}\\
=& \chi_{\delta}O(\norm{v}_{L^2}^2+\norm{v}_{L^2}|c-c^*|) \mbox{ as } \norm{v}_{L^2}+|c-c^*| \to 0.
\end{align*}
Then,  for a solution $(v,c,\rho) $ to the system \eqref{LZKeq-1}--\eqref{LZKeq-2} and $t \in \R$
\[P_1v(0)=P_1v(t), \quad P_2v(0)=P_2v(t).\]
Especially, for initial data $(v(0),c(0),\rho(0))$ satisfying the orthogonality condition \eqref{cond-orth}, the solution $(v,c,\rho)$ of the system \eqref{LZKeq-1}--\eqref{LZKeq-2} also satisfies the orthogonality condition \eqref{cond-orth}.

To solve the system \eqref{LZKeq-1}--\eqref{LZKeq-2}, we define the Bourgain space $X^{s,b}$ related to the linear part of \eqref{ZKeq} as the completion of the Schwartz space under the norm
\begin{align*}
\norm{u}_{X^{s,b}}=\Bigl( \int_{\R \times \R} \sum_{L \eta \in \Z} \tbr{\tau-\xi(\xi^2+\eta^2)}^{2b}\tbr{\sqrt{3\xi^2+\eta^2}}^{2s}|\tilde{u}(\tau,\xi,\eta)|^2\, d\tau\, d\xi \Bigr)^{\frac{1}{2}},
\end{align*}
where $\tbr{x}=1+|x|$ and $\tilde{u}$ is the space-time Fourier transform of $u$.
For $T>0$, we define the localized space $X^{s,b}_T$ of $X^{s,b}$ by the norm
\[\norm{u}_{X^{s,b}_T}=\inf\{\norm{v}_{X^{s,b}};v \in X^{s,b}, v(t)=u(t) \mbox{ for } t \in [-T,T]\}.\]
Let $\theta$ be a smooth function with $0 \leq \eta \leq 1$ and $\theta(t)=1$ for $|t|\leq 1$, $\theta(t)=0$ for $|t|\geq 2$, let $\theta_T(t)=\theta(t/T)$ for $T>0$.
Then, we have the following linear estimates in \cite{G T V,M P,M S T}.
\begin{proposition}\label{prop-linear-est}
Let $s\geq 0$, $T>0$, $b>\frac{1}{2}$ and $ b_1\leq 0 \leq  b_2 \leq b_1+1$.
Then
\begin{align}
&\norm{\theta(t) e^{-t\partial_x \Delta} u_0}_{X^{s,b}} \lesssim_{s,b} \norm{u_0}_{H^s}, \label{eq-linear-est-1}\\
&\norm{u}_{L^{\infty}H^s(\RTL)} \lesssim_{s,b} \norm{u}_{X^{s,b}}, \label{eq-linear-est-2}\\
&\norm{\theta_T(t) \int_0^t e^{-(t-t')\partial_x \Delta}g(t') \,dt'}_{X^{s, b_2}} \lesssim_{s,b_1,b_2} T^{1-b_2+b_1} \norm{g}_{X^{s, b_1}}, \label{eq-linear-est-3} \\
&\norm{\partial_x( Q u)}_{X^{s,0}} \lesssim_{s,b} ( \norm{\partial_x Q}_{L^{\infty}_tW^{s,\infty}_{x,y}}+ \sum_{|\alpha|\leq s} \norm{\partial^{\alpha} Q}_{L_x^2L_{ty}^{\infty}})\norm{u}_{X^{s,b}} \label{eq-linear-est-4},
\end{align}
for $u_0 \in H^s(\RTL)$, $u \in X^{s,b}$, $Q \in L_t^{\infty}((-\infty,\infty), W_{xy}^{s,\infty}\cap H_{xy}^s)$ and $g \in X^{s,b_1}$, where $\partial^{\alpha}=\partial_x^{\alpha_1}\partial_y^{\alpha_2}$ and $\alpha=(\alpha_1,\alpha_2)$.
\end{proposition}
To estimate the nonlinear term of the system \eqref{LZKeq-1}--\eqref{LZKeq-2}, we use the following bilinear estimate by Molinet and Pilod \cite{M P}.
\begin{proposition}\label{prop-bilinear-est}
Let $s \geq 1$. 
Then, there exists $b_* >0$ such that
\begin{align}\label{eq-bilinear-est}
\norm{\partial_x(uv)}_{X^{s,-\frac{1}{2}+2b}}\lesssim_{s, b} \norm{u}_{X^{s,\frac{1}{2}+b}}\norm{v}_{X^{s,\frac{1}{2}+b}},
\end{align}
for $0<b <b_*$, $u,v \in X^{s,\frac{1}{2}+b}$.
\end{proposition}

To prove the estimate for solutions to the system \eqref{LZKeq-1}--\eqref{LZKeq-2}, we show the conservation of the linearized energy norm for the linearized system of the system \eqref{LZKeq-1}--\eqref{LZKeq-2}.
\begin{lemma}\label{lem-linear-property}
Let $c_0>0$. The solution $(v,\rho)$ to the system
\begin{align}
v_t=\partial_x\mathbb{L}_{c^*}v+(\dot{\rho}-c_0)\partial_xQ_{c^*}\label{Lsys-LZK-1}\\
\dot{\rho}-c_0=\frac{(v,\mathbb{L}_{c^*}\partial_x^2Q_{c^*})_{L^2}}{\norm{\partial_xQ_{c^*}}_{L^2}^{2}}\label{Lsys-LZK-2}
\end{align} 
with an initial data $(v(0),\rho(0)) \in H^{1}(\RTL) \times \R$ satisfies
\begin{align}\label{Lproperty-gamma}
\begin{array}{l}
\norm{P_{\gamma}v(t)}_{E}=\norm{P_{\gamma}v(0)}_{E}, \quad (v(t),\partial_xQ_{c^*})_{L^2}=(v(0),\partial_xQ_{c^*})_{L^2},  \\
(v(t),Q_{c^*})_{L^2}=(v(0),Q_{c^*})_{L^2}
\end{array}
\end{align}
and
\begin{align}\label{linear-eq-1}
v(t)=e^{t\partial_x\mathbb{L}_{c^*}}v(0)+\int_0^t \frac{(e^{s\partial_x\mathbb{L}_{c^*}}v(0), \mathbb{L}_{c^*}\partial_x^2Q_{c^*})_{L^2}}{\norm{\partial_xQ_{c^*}}_{L^2}^2} \partial_x Q_{c^*} \, ds
\end{align}
\end{lemma}

\proof
Let $(v,\rho)$ be the solution to the system \eqref{Lsys-LZK-1}--\eqref{Lsys-LZK-2} with a smooth initial data $(v(0),\rho(0))$.
Then, we have
\begin{align*}
\partial_t \norm{P_{\gamma}v}_{E}^2
=& 2\rbr{\partial_x \mathbb{L}_{c^*} P_{\gamma} v - P_1(\partial_x \mathbb{L}_{c^*} P_{\gamma}v), \mathbb{L}_{c^*} P_{\gamma}v}_{L^2}=0.
\end{align*}
Since $\mathbb{L}_{c^*}\partial_x Q_{c^*}=0$, 
\begin{align*}
\partial_t (v(t),Q_{c^*})_{L^2}=(\partial_x\mathbb{L}_{c^*}v(t), Q_{c^*})_{L^2}=0.
\end{align*}
By the system \eqref{Lsys-LZK-1} and \eqref{Lsys-LZK-2}, we have
\begin{align*}
\partial_t (v(t),\partial_xQ_{c^*})_{L^2}=\Bigl(\partial_x\mathbb{L}_{c^*}v(t)+\frac{(v(t),\mathbb{L}_{c^*}\partial_x^2Q_{c^*})_{L^2}}{\norm{\partial_xQ_{c^*}}_{L^2}^{2}}, \partial_xQ_{c^*}\Bigr)_{L^2}=0.
\end{align*}
By the density argument we obtain \eqref{Lproperty-gamma}.
The orthogonality $(\partial_x Q_{c^*},\mathbb{L}_{c^*}\partial_x^2 Q_{c^*})_{L^2}=0$ yields the formula \eqref{linear-eq-1}.
\qed

Let $\mathcal{A}$ as 
\[\mathcal{A}v=\partial_x\mathbb{L}_{c^*}v+\frac{(v,\mathbb{L}_{c^*}\partial_x^2Q_{c^*})_{L^2}}{\norm{\partial_xQ_{c^*}}_{L^2}^2}\partial_xQ_{c^*}.\]
We define the system
\begin{align}
w_t =& -\partial_x \Delta w-2\partial_x((\tau_{\rho_*}Q_{c^*})w)+(\dot{\rho}-c)\tau_{\rho_*}\partial_xQ_{c^*} -\dot{c}\tau_{\rho_*}\partial_cQ_{c^*}\notag\\
&+\chi_\delta(w, c-c^*)\tilde{N}(w,\rho,c) \label{Leq-w-1}\\
\begin{pmatrix}
\dot{\rho}-c\\
\dot{c}
\end{pmatrix}
=&
\begin{pmatrix}
\norm{\partial_xQ_{c^*}}_{L^2}^{-2}(w,\tau_{\rho_*}(\mathbb{L}_{c^*}\partial_x^2Q_{c^*}))_{L^2}\\
0
\end{pmatrix}
+ \tilde{\bm{N}}_{\delta}(w,c,\rho),\label{Leq-w-2}
\end{align}
where 
\[ \tilde{N}(w,c,\rho)=\partial_x[ - w^2+2w\tau_{\rho_*}(Q_{c^*}-Q_c)+ (\dot{\rho}-c)\tau_{\rho_*}(Q_{c}-Q_{c^*})]-\dot{c}\tau_{\rho_*}\partial_c(Q_c-Q_{c^*}),\]
$\tilde{\bm{N}}_{\delta}(w,c,\rho)=\bm{N}_{\delta}(\tau_{-\rho_*}w,c)$ and 
\[\rho_*(w,c,\rho,t)=\rho_*(t)=c^*t+\int_0^{t} \chi_{\delta}(w(s),c(s)-c^*)(\dot{\rho}(s)-c^*)\,ds.\]
We solve the system \eqref{LZKeq-1}--\eqref{LZKeq-2} by the following theorem.
\begin{theorem}\label{thm-gwp-LZKeq}
The system $\eqref{LZKeq-1}$--$\eqref{LZKeq-2}$ is globally well-posed in $H^1(\RTL) \times (0,\infty) \times \R$.
More precisely there exists $b > \frac{1}{2}$ for every $(v_0,c_0,\rho_0) \in H^1(\RTL) \times (0,\infty) \times \R$ and $T>0$ there exists a unique solution $(w, c,\rho)$ of the system $\eqref{Leq-w-1}$--$\eqref{Leq-w-2}$ such that $(w(0),c(0),\rho(0))=(v_0,c_0,\rho_0)$,
\[(w,\dot{c},\dot{\rho}-c) \in X^{1,b}_T\times L^2(-T,T) \times L^2(-T,T),\]
and $(\tau_{-\rho_*}w,c,\rho)$ is a solution to the system $\eqref{LZKeq-1}$--$\eqref{LZKeq-2}$ with initial data $(v(0),c(0),\rho(0))$.
Moreover, the flow map of the system $\eqref{Leq-w-1}$--$\eqref{Leq-w-2}$ is Lipschitz continuous on bounded sets of $H^1(\RTL) \times (0,\infty) \times \R$ and  there exists $T^*>0$ such that for any $0< \delta<1$ and initial data $(v(0), c(0), \rho(0)) \in H^1(\RTL) \times (0,\infty) \times \R $ the solution $(v,c,\rho)$ to the system $\eqref{LZKeq-1}$--$\eqref{LZKeq-2}$ with the initial data $(v(0), c(0), \rho(0))$ satisfies
\begin{align}
\sup_{|t|\leq T^*} \norm{v(t)}_{E}+\norm{\dot{c}}_{L^2(-T^*,T^*)}+\norm{\dot{\rho}-c}_{L^2(-T^*,T^*)} \lesssim \norm{v(0)}_{E},\label{eq-gwp-est-1}\\
\sup_{|t|\leq T^*}\norm{P_d(v(t)-e^{t\mathcal{A}}v(0))}_{E} \lesssim \min \{\norm{v(0)}_{E}, \delta\}^2,\label{eq-gwp-est-2} \\
\sup_{|t|\leq T^*} \bigl|\norm{P_{\gamma}v(t)}_{E}^2-\norm{P_{\gamma}v(0)}_{E}^2\bigr| \lesssim \min \{\norm{v(0)}_{E}, \delta\}^3,\label{eq-gwp-est-3}
\end{align}
where the implicit constants and $T^*$ do not depend on $\delta$ and $(v(0),c(0),\rho(0))$.
\end{theorem}
\begin{remark}
The above theorem does not imply the Lipschitz continuity of the flow map of the system $\eqref{LZKeq-1}$--$\eqref{LZKeq-2}$.
\end{remark}
\proof
First, we consider the case $\norm{v(0)}_{H^1}^2+|c(0)-c^*|^2>16\delta^2$.
Since $\chi_{\delta}(v(0),c(0)-c^*)=0$, there exists $T_1>0$ such that the solution $(v_1(t), c_1(t), \rho_1(t))$ to the system \eqref{Lsys-LZK-1}--\eqref{Lsys-LZK-2} with the initial data $(v_1(0),c_1(0), \rho_1(0))=(v(0),c(0),\rho(0))$ is a solution to the system\eqref{LZKeq-1}--\eqref{LZKeq-2} on the time interval $ (-T_1,T_1)$ with the initial data $(v(0),\rho(0),c(0))$ and
\[\sup_{t \in (-T_1,T_1)} \norm{v_1(t)}_{H^1}^2+|c_1(t)-c^*|^2>8\delta^2.\]
If $|c(0)-c^*|>2\delta$, then $c(t)=c_1(t)=c(0)$ for $t >0$.
Thus, the solution $(v,c,\rho)=(v_1,c_1,\rho_1)$ is global in time.
Next we consider  the case $\norm{v(0)}_{H^1}>2\sqrt{3}\delta$.
Then, we have 
\[\norm{v_1(t)}_{H^1}\leq \norm{e^{t\mathcal{A}}}_{H^1 \to H^1} \norm{v(0)}_{H^1}.\]
Therefore, $T_1$ is independent of initial data $(v(0),c(0),\rho(0))$ and $\delta>0$ if  $\norm{v(0)}_{H^1}^2+|c(0)-c^*|^2>16\delta^2$.
By the continuity of solutions to the system \eqref{LZKeq-1}--\eqref{LZKeq-2} in time, solutions to the system \eqref{LZKeq-1}--\eqref{LZKeq-2} with the initial data $(v(0),c(0),\rho(0))$ is also a solution to the system \eqref{Lsys-LZK-1}--\eqref{Lsys-LZK-2} on the time interval $ (-T_1,T_1)$.
Therefore, we obtain the existence and the uniqueness of solution to the system \eqref{LZKeq-1}--\eqref{LZKeq-2} on $C((-T_1,T_1),H^1(\RTL))$ for initial data $(v(0),c(0),\rho(0))$ satisfying $\norm{v(0)}_{H^1}^2+|c(0)-c^*|^2>16\delta^2$.
From Lemma \ref{lem-linear-property} the solution $(v,c,\rho)$ to \eqref{LZKeq-1}--\eqref{LZKeq-2} with $\norm{v(0)}_{H^1}^2+|c(0)-c^*|^2>16\delta^2$ satisfies the estimates \eqref{eq-gwp-est-1}--\eqref{eq-gwp-est-3} on $[-T_1/2,T_1/2]$, where $T_1/2$ does not depend on $\delta$ and initial data $(v(0),c(0),\rho(0))$.



Second, we treat the case $|c(0)-c^*|\leq 4\delta$ and $\norm{v(0)}_{H^1}\leq 4\delta$. 
We consider solutions $(w,c,\rho)$ to the system \eqref{Leq-w-1}--\eqref{Leq-w-2}.
To show the global well-posedness of the system \eqref{Leq-w-1}--\eqref{Leq-w-2}, we apply the contraction mapping theorem to
\begin{align*}
\Phi_T(w,c,\rho)(t)=&e^{-t\partial_x\Delta }v(0)+\int_0^te^{-(t-s)\partial_x\Delta}[-2\partial_x((\tau_{\rho_*}Q_{c^*})w)\\
&+(\dot{\rho}-c)\tau_{\rho_*}\partial_xQ_{c^*}-\dot{c}\tau_{\rho_*}\partial_cQ_{c^*}+\chi_{\delta}(w,(c-c^*))\tilde{N}(w,c,\rho)]ds,\\
\Psi_T(w,c,\rho)=&
\begin{pmatrix}
\Psi_{T,1}(w,c,\rho)\\
\Psi_{T,2}(w,c,\rho)
\end{pmatrix}=
\begin{pmatrix}
\norm{\partial_xQ_{c^*}}_{L^2}^{-2}(w,\tau_{\rho_*}(\mathbb{L}_{c^*}\partial_x^2 Q_{c^*}))_{L^2}\\
0
\end{pmatrix}
+ \tilde{\bm{N}}_{\delta}(w,c,\rho)
\end{align*}
on $X^{1,b}_T\times L^2_T \times L^2_T$, where $L^2_T=L^2(-T,T)$.
Applying Proposition \ref{prop-linear-est}, for $\frac{1}{2}<b<\min\{\frac{3}{4},\frac{1}{2}+b_*\}$ and $0<T<1$ we have
\begin{align}
&\norm{\Phi_T(w,c,\rho)}_{X^{1,b}_T} \notag \\
\lesssim_{b}& \norm{v(0)}_{H^1(\RTL)} + T^{1-b}\norm{\partial_x((\tau_{\rho_*}Q_{c^*})w)}_{X^{1,0}_T}+T^{1-b}\norm{(\dot{\rho}-c)\tau_{\rho_*}\partial_xQ_{c^*}}_{X^{1,0}_T} \notag\\
&+ T^{1-b}\norm{\dot{c}\tau_{\rho_*}\partial_cQ_{c^*}}_{X^{1,0}_T}+T^{b-\frac{1}{2}}\norm{\chi_{\delta}\partial_x(w^2)}_{X^{1,2b-\frac{3}{2}}_T}+T^{1-b}\norm{\chi_{\delta}\partial_x(w\tau_{\rho_*}(Q_{c^*}-Q_{c}))}_{X^{1,0}_T} \notag\\
&+T^{1-b}\norm{\chi_{\delta}(\dot{\rho}-c)\tau_{\rho_*}\partial_x(Q_c-Q_{c^*})}_{X^{1,0}_T}+T^{1-b}\norm{\chi_{\delta}\dot{c}\tau_{\rho_*}\partial_c(Q_c-Q_{c^*})}_{X^{1,0}_T} \notag\\
\lesssim_{b} & \norm{v(0)}_{H^1(\RTL)} + T^{1-b}(1+\norm{\tau_{\rho_*}Q_{c^*}}_{L^2_xL^{\infty}_{t,y}}+\norm{\tau_{\rho_*} (Q_c-Q_{c^*})}_{L^2_xL^{\infty}_{t,y}}\norm{w}_{X^{1,b}_T} \notag\\
&+T^{1-b}(1+|c(0)-c^*|+T^{\frac{1}{2}}\norm{\dot{c}}_{L^{2}_T})(\norm{\dot{\rho}-c}_{L^{2}_T}+\norm{\dot{c}}_{L^{2}_T})+T^{b-\frac{1}{2}}\norm{\chi_{\delta}\partial_x(w^2)}_{X^{1,2b-\frac{3}{2}}_T}.\label{gwp-eq-0}
\end{align}
Since
\[ \sup_{-T\leq t \leq T, y \in \T_L} |\tau_{\rho_*}Q_{c^*}|(t,x,y) \leq \sup_{|a|\leq (Tc^*+T|c(0)-c^*|+T^\frac{1}{2}\norm{\dot{\rho}-c}_{L^{2}_T}+T^{\frac{3}{2}}\norm{\dot{c}}_{L^{2}_T})} Q_{c^*}(x+a),\]
we have 
\begin{align}
&\norm{\tau_{\rho_*}Q_{c^*}}_{L^2_xL^{\infty}_{t,y}((-T,T)\times \T_L)} \notag \\
\lesssim & \norm{Q_{c^*}}_{L^2}+(Tc^*+T|c(0)-c^*|+T^\frac{1}{2}\norm{\dot{\rho}-c}_{L^{2}_T}+T^{\frac{3}{2}}\norm{\dot{c}}_{L^{2}_T})^{1/2}\norm{Q_{c^*}}_{L^{\infty}(\RTL)} \notag \\
\lesssim &1+T^{\frac{1}{2}}|c(0)-c^*|+T^{\frac{1}{2}} +T^{\frac{1}{4}}\norm{\dot{\rho}-c}_{L^{2}_T}+T^{\frac{3}{4}}\norm{\dot{c}}_{L^{2}_T}.\label{gwp-eq-1}
\end{align}
By the similar calculation, we obtain 
\begin{align}\label{gwp-eq-2}
&\norm{\theta_{2T}\tau_{\rho_*}(Q_c-Q_{c^*})}_{L^2_xL^{\infty}_{t,y}((-T,T)\times \T_L)}\notag \\
&\lesssim (|c(0)-c^*|+T^{\frac{1}{2}}\norm{\dot{c}}_{L^{2}_T})(1+T^{\frac{1}{2}}+T^{\frac{1}{2}}|c(0)-c^*| +T^{\frac{1}{4}}\norm{\dot{\rho}-c^*}_{L^{2}_T}+T^{\frac{3}{4}}\norm{\dot{c}}_{L^{2}_T}).
\end{align}
Since 
\[ \norm{\partial_t \norm{u}_{H^1_{x,y}}}_{L^2_t} \leq \norm{\norm{\partial_t e^{t\partial_x\Delta}u}_{H^1_{x,y}}}_{L^2_t} \lesssim \norm{u}_{X^{1,1}}\]
for $u \in X^{1,1}$, by the interpolation theorem we have
\begin{align}\label{gwp-eq-2-1}
\norm{\norm{w}_{H^1_{x,y}}}_{H^b_t} \lesssim \norm{w}_{X^{1,b}}.
\end{align}
Applying Proposition \ref{prop-bilinear-est} to $\norm{\chi_{\delta}\partial_x(w^2)}_{X^{1,2b-\frac{3}{2}}_T}$, we obtain that there exists 
\[C=C(\delta^{-2}(\norm{w}_{L^\infty((-T,T) H^1_{x,y})}^2+|c(0)-c^*|^2+T\norm{\dot{c}}_{L^2_T}^2))>0\]
such that
\begin{align}
&\norm{\chi_{\delta}\partial_x(w^2)}_{X^{1,2b-\frac{3}{2}}_T}\notag \\
\lesssim \ &  (\norm{\chi_{\delta}( w, (c-c^*))-1}_{H^b_T}+1)\norm{w}_{X^{1,b}_T}^2\notag\\
\lesssim & C  (1+\delta^{-2}\norm{ w}_{X^{1,b}_T}^2+\delta^{-2}(1+T^{\frac{1}{2}})|c(0)-c^*|^2+ \delta^{-2}(1+T^{\frac{3}{2}})\norm{ \dot{c}}_{L^{2}_T}^2)\norm{w}_{X^{1,b}_T}^2,\label{gwp-eq-5}
\end{align}
where 
\[\norm{u}_{H^b_T}=\inf\{ \norm{v}_{H^b(\R)};v \in H^b(\R), v(t)=u(t) \mbox{ for } t \in [-T,T]\}.\]
By the simple calculation we have
\begin{align}\label{gwp-eq-6}
&\norm{\Psi_{T,1}(w,c,\rho)}_{L^{2}_T}+\norm{\Psi_{T,2}(w,c,\rho)}_{L^{2}_T} \notag \\
\lesssim& T^{\frac{1}{2}}\norm{w}_{X^{1,b}_T}+T^{\frac{1}{2}}\norm{w}_{X^{1,b}_T}^3+T^{\frac{1}{2}}|c(0)-c^*|\norm{w}_{X^{1,b}_T}+T^{\frac{1}{2}}|c(0)-c^*|^2\norm{w}_{X^{1,b}_T}\notag \\
&+T^{\frac{3}{2}}\norm{\dot{c}}_{L^{2}_T}^2+T^{2}\norm{\dot{c}}_{L^{2}_T}^3
\end{align}
for $b>\frac{1}{2}$ and $0<T<1$.
From \eqref{eq-linear-est-4} and \eqref{gwp-eq-0}--\eqref{gwp-eq-5} we obtain that 
\begin{align}
&\norm{\Phi_T(w)}_{X^{1,b}_T}\notag \\
\lesssim_b&  \norm{v(0)}_{H^1} +T^{1-b}(\norm{w}_{X^{1,b}_T}+\norm{\dot{\rho}-c}_{L^{2}_T}+\norm{\dot{c}}_{L^{2}_T})(1+|c(0)-c^*|+\norm{\dot{\rho}-c}_{L^{2}_T}+\norm{\dot{c}}_{L^{2}_T})^2\notag\\
&+T^{b-\frac{1}{2}}C\norm{w}_{X^{1,b}_T}^2\label{gwp-eq-7}
\end{align}
for $\frac{1}{2}<b<\min\{\frac{3}{4},\frac{1}{2}+b_*\}$ and $0<T<1$, where $C$  depend on $\delta^{-2}(\norm{w}_{X^{1,b}_T}^2+|c(0)-c^*|^2+\norm{\dot{c}}_{L^2_T}^2)$.
From the same calculation as \eqref{gwp-eq-0}--\eqref{gwp-eq-7} we obtain the estimate of the difference
\begin{align}
&\norm{\Phi_T(w_1,c_1,\rho_1)-\Phi_T(w_2,c_2,\rho_2)}_{X^{1,b}_T}+\norm{\Psi_{T,1}(w_0,c_0,\rho_0)-\Psi_{T,1}(w_1,c_1,\rho_1)}_{X^{1,b}_T}\notag\\
&+\norm{\Psi_{T,2}(w_0,c_0,\rho_0)-\Psi_{T,2}(w_1,c_1,\rho_1)}_{X^{1,b}_T}\notag \\
 \lesssim_b &(T^{1-b}+T^{b-\frac{1}{2}})C(\norm{w_1-w_2}_{X^{1,b}_T}+\norm{\dot{\rho}_1-\dot{\rho}_2-c_1+c_2}_{L^{2}_T}+\norm{\dot{c}_1-\dot{c}_2}_{L^{2}_T})\notag \\
& \times [1+|c(0)-c^*|+\max_{j=1,2}(\norm{w_j}_{X^{1,b}_T}+\norm{\dot{\rho}_j-c_j}_{L^{2}_T}+\norm{\dot{c}_j}_{L^{2}_T})]^2\notag 
\end{align} 
for $\frac{1}{2}<b<\min\{\frac{3}{4},\frac{1}{2}+b_*\}$ and $0<T<1$, where $C$ depend on $\delta^{-2}\max_{j=0,1}(\norm{w_j}_{X^{1,b}_T}^2+|c(0)-c^*|^2+\norm{\dot{c}_j}_{L^2_T})$.
Thus, there exist $C_1,C_2 >0$ and $T_{\delta,v(0)}=T(\delta,\norm{v(0)}_{H^1})>0$ such that the mapping $(\Phi_{T_{\delta,v(0)}}, \Psi_{T_{\delta,v(0)}})$ is the contraction mapping on $\mathcal{B}(C_1\norm{v(0)}_{H^1},C_2\norm{v(0)}_{H^1})$, where $C_1$ and $C_2$ do not depend on $\delta$ and the initial data $(v(0), c(0),\rho(0))$,
\[\mathcal{B}(r_1,r_2)=\{(w,\dot{c},\dot{\rho}-c) \in X^{1,b}_T\times (L^2_T)^2: \norm{w}_{X^{1,b}_T}<r_1, \norm{\dot{\rho}-c}_{L^{2}_T}+\norm{\dot{c}}_{L^{2}_T}<r_2\}.\]
Therefore, the system \eqref{Leq-w-1}--\eqref{Leq-w-2} is locally well-posed in $H^1(\RTL) \times (0,\infty) \times \R $ and $(\tau_{-\rho_*}w,c,\rho)$ is the solution to the system \eqref{LZKeq-1}--\eqref{LZKeq-2} with initial data $(v(0),c(0),\rho(0))$.
Moreover, for $0<\delta<1$ and $\frac{1}{2}<b<\min\{\frac{3}{4},\frac{1}{2}+b_*\}$ there exist $C, T_{\delta}>0$ such that for any $v(0) \in H^{1}(\RTL)$ with $\norm{v(0)}_{H^1}\leq 4\delta$ and $ (w,c,\rho) \in \mathcal{B}(C\delta,C\delta)$ we have
\begin{align}
&\norm{\Phi_{T_{\delta}}(w,c,\rho)}_{X^{1,b}_{T_\delta}}\notag\\
\lesssim_b&  \norm{v(0)}_{H^1} +(T_{\delta}^{1-b}+T_{\delta}^{b-\frac{1}{2}})(\norm{w}_{X^{1,b}_{T_\delta}}+\norm{\dot{\rho}-c}_{L^{2}_t}+\norm{\dot{c}}_{L^{2}_{T_\delta}})\notag \\
&\qquad \qquad\qquad\qquad\times (1+\norm{w}_{X^{1,b}_{T_\delta}}+|c(0)-c^*|+\norm{\dot{\rho}-c}_{L^{2}_{T_\delta}}+\norm{\dot{c}}_{L^{2}_{T_\delta}})^2
< \frac{C\delta}{2}, \notag
\end{align}
where $C$ does not depend on $\delta$.
 By the continuity argument, we obtain
that there exists $T^*>0$ such that  the solution $(v,c,\rho)$ to the system \eqref{LZKeq-1}--\eqref{LZKeq-2} exists on $[-T^*,T^*]$ and
\[\sup_{|t|\leq T^*} \norm{v(t)}_{E} \lesssim \sup_{|t|\leq T^*} \norm{v(t)}_{H^1} \lesssim \norm{w_{T^*}}_{X^{1,b}} \lesssim \norm{v(0)}_{H^1} \lesssim \norm{v(0)}_{E} \lesssim \delta,\]
where $w_{T^*}$ is the fixed point  by $(\Phi_{T^*},\Psi_{T^*})$ on $X^{1,b}_{T^*}\times (L^2_{T^*})^2$ with the initial data $(v(0),c(0),\rho(0))$ and, the implicit constants and $T^*$ do not depend on $\delta$.
Thus, \eqref{eq-gwp-est-1} holds.
Since
\[\norm{\dot{c}}_{L^2_{T^*}}=O(\norm{v}_{H^1}^2)\]
for the solution $(v,\rho,c)$ to the system \eqref{LZKeq-1}--\eqref{LZKeq-2}, we obtain \eqref{eq-gwp-est-2}.

The identity $\mathbb{L}_{c^*}[P_{\gamma},\partial_x\mathbb{L}_{c^*}]=0$ yields
\begin{align}\label{eq-gwp-9}
(P_{\gamma}(\partial_x \mathbb{L}_{c^*} v + (\dot{\rho}-c) \partial_xQ_{c^*}-\dot{c}\partial_cQ_{c^*}), \mathbb{L}_{c^*} P_{\gamma}v)_{L^2}=0
\end{align}
 for $(v,c,\rho) \in H^3(\RTL) \times C^1(0,\infty) \times C^1(\R)$.
By the Plancherel theorem, we have
\begin{align*}
&\Bigl| \int_0^T \int_{\RTL} (-\Delta u )v\,dtdxdy  \Bigr| \notag \\
=& \Bigl | \int_{-\infty}^{\infty} \int_{\R}\sum_{\eta} \tbr{\tau-\xi(\xi^2+\eta^2)}^{\beta}  \widetilde{ 1_{[0,T]} \nabla \tau_\rho u } \times  \tbr{\tau -\xi(\xi^2+\eta^2)}^{-\beta}\widetilde{\nabla \tau_\rho v}\, d\tau d\xi \Bigr| \notag \\
\lesssim & \Bigl(\norm{\tbr{\tau}^{\beta}\widetilde{ 1_{[0,T]}}}_{L^2_{\tau}}\norm{\tau_\rho u}_{X^{1,b}}+\norm{\widetilde{ 1_{[0,T]}}}_{L^p_{\tau}}\norm{\tbr{\tau}^{\beta-b}}_{L^{\frac{p}{p-1}}_{\tau}}\norm{\tau_\rho u}_{X^{1,b}}\Bigr)\norm{\tau_\rho v}_{X^{1,-\beta}}.\notag \\
\end{align*}
Since
\[ \int_{-\infty}^{\infty}1_{[0,T]} e^{-it\tau}dt=\frac{i(1-e^{-iT\tau})}{\tau},\]
we have
\begin{align}\label{eq-gwp-a-3}
\Bigl| \int_0^T \int_{\RTL} (-\Delta u )v\,dtdxdy  \Bigr| \lesssim \norm{\tau_\rho u}_{X^{1,b}}\norm{\tau_\rho v}_{X^{1,-\beta}}
\end{align}
for $\tau_\rho u \in X^{1,b}, \tau_\rho v \in X^{1,-\beta}$, $\rho \in L^{\infty}_t$, $p>1$ and $0\leq \beta <\frac{1}{2}$ with
\[\frac{(b-\beta)p}{p-1}>1.\] 
From Proposition \ref{prop-bilinear-est} and the inequalities \eqref{eq-gwp-est-1} and \eqref{eq-gwp-a-3}, we have
\begin{align}\label{eq-gwp-11}
\Bigl| \int_0^T \int_{\RTL} (\mathbb{L}_{c^*} v )(\chi_\delta(v,c-c^*)N(v,c,\rho))\,dtdxdy  \Bigr| \lesssim \min \{\norm{v(0)}_{H^1}, \delta\}^3.
\end{align}
From \eqref{eq-gwp-9} and \eqref{eq-gwp-11}, by the energy estimate we obtain \eqref{eq-gwp-est-3}.

Combining above two cases, we have the global well-posedness of the \eqref{LZKeq-1}--\eqref{LZKeq-2} on $H^1(\RTL) \times (0,\infty)\times \R$ and the estimate \eqref{eq-gwp-est-1}--\eqref{eq-gwp-est-3}.

\qed

Next, we define a mobile distance which was introduced in \cite{N S 2}.
Let $C_2$ be a large real constant and $\phi$ be the smooth positive non-deceasing function with 
\[\phi(r)=\begin{cases}
1, & r \leq C_2, \\
r, & r \geq 2 C_2. 
\end{cases}\]
We define $\phi_{\delta}$ by 
\[\phi_{\delta}(u)= \phi\Bigl( \delta^{-1}\norm{P_{\gamma}u}_{E}\Bigr)\]
for $u \in H^1(\RTL)$.
In this paper, to treat the term $\dot{\rho}-c$ in the system \eqref{LZKeq-1}--\eqref{LZKeq-2} which has same order of $v$, we replace a correction term of the mobile distance in \cite{N S 2} by $\delta |q|^2\phi_{\delta}(v_{1-j})^2$.
\begin{definition}
Let $\delta>0$.
We define the mobile distance $\mathfrak{m}_{\delta} : (H^{1}(\RTL)\times (0,\infty))^2 \to [0,\infty)$ by
\begin{align*}
(\mathfrak{m}_{\delta}(\bm{v}_0,\bm{v}_1))^2=&\norm{P_d(v_0-v_1)}_{E}^2+\inf_{q \in \R, j=0,1}\bigl (\norm{P_{\gamma}v_j-\tau_{q}(P_{\gamma}v_{1-j})}_{E}^2+ \delta |q|^2\phi_{\delta}( v_{1-j})^2\bigr)\\
&+|\log c_0-\log c_1|^2
\end{align*}
for $\bm{v}_0=(v_0,c_0),\bm{v}_1=(v_1,c_1) \in H^{1}(\RTL)\times (0,\infty) $.
\end{definition}
In the following lemma, we show $\mathfrak{m}_\delta$ is a complete quasi-distance on $H^1(\RTL)\times (0,\infty) $.
\begin{lemma}\label{lem-md-pro}
Let $0<\delta <1$. 
$\mathfrak{m}_{\delta}$ satisfies the following.
\begin{enumerate}
  \setlength{\parskip}{0.1cm} 
  \setlength{\itemsep}{0.05cm} 
\renewcommand{\labelenumi}{\rm (\roman{enumi})}
\item $\mathfrak{m}_{\delta}(\bm{v}_0,\bm{v}_1)=\mathfrak{m}_{\delta}(\bm{v}_1,\bm{v}_0) \geq 0$, where the equality holds iff $\bm{v}_0=\bm{v}_1$.
\item $\mathfrak{m}_{\delta}(\bm{v}_0,\bm{v}_1) \leq C(\mathfrak{m}_{\delta}(\bm{v}_0,\bm{v}_2)+\mathfrak{m}_{\delta}(\bm{v}_2,\bm{v}_1))$, for some absolute constant $C>0$ which does not depend on $\delta$.
\item If $\mathfrak{m}_{\delta}(\bm{v}_n,\bm{v}_m) \to 0$ $(n,m \to \infty)$, then $\{\bm{v}_n\}_n$ converges in $H^1(\RTL)\times (0,\infty)$.
\item For $\bm{v}_0=(v_0,c_0), \bm{v}_1=(v_1,c_1)\in H^1(\RTL)\times(0,\infty)$
\begin{align*}
|\norm{v_1}_{H^1}-\norm{v_0}_{H^1}|+\norm{v_0-v_1}_{L^2}+|\log c_0-\log c_1|
 \lesssim &\mathfrak{m}_{\delta}(\bm{v}_0,\bm{v}_1) \\
\lesssim & \norm{v_0-v_1}_{H^1}+|\log c_0-\log c_1|,
\end{align*}
where the implicit constants do not depend on $\delta$.
\end{enumerate}
\end{lemma}
\proof
By the definition of $\mathfrak{m}_{\delta}$, we have  (i).
The right inequality of (iv) follows the equivalence of $\norm{\cdot}_{H^1}$ and $\norm{\cdot}_{E}$.
From the inequalities
\begin{align*}
&|\norm{v_0}_{H^1}^2-\norm{v_1}_{H^1}^2|  \\
\lesssim & \inf_{q\in \R }\min_{j=0,1}(| \norm{P_{\gamma}v_j}_{H^1}^2-\norm{\tau_{q}P_{\gamma}v_{j-1}}_{H^1}^2|+\norm{v_0}_{H^1}(\norm{P_{\gamma}v_{j}-\tau_qP_{\gamma}v_{j-1}}_{H^1}+|q|\norm{P_{\gamma}v_{1-j}}_{H^1}))\\\
& +(\norm{v_0}_{H^1}+\norm{v_1}_{H^1})\norm{P_d(v_0-v_1)}_{E}\\
\lesssim & (\norm{v_0}_{H^1}+\norm{v_1}_{H^1})\mathfrak{m}_{\delta} (\bm{v}_0,\bm{v}_1),
\end{align*}
and
\begin{align*}
\norm{v_0-v_1}_{L^2} 
\leq & \norm{P_{\gamma}(v_0-v_1)}_{L^2}+\norm{P_d(v_0-v_1)}_{L^2}\\
\lesssim & \inf_{q\in \R, j=0,1}\{ \norm{P_{\gamma}v_j-\tau_{q} P_{\gamma}v_{1-j}}_{L^2}+|q|\norm{\nabla P_{\gamma}v_{1-j}}_{L^2}\}+\norm{P_d(v_0-v_1)}_{E}\\
\lesssim & \mathfrak{m}_{\delta}(\bm{v}_0,\bm{v}_1),
\end{align*}
we obtain the left inequality of (iv).

Next, we show the quasi-triangle inequality (ii).
Let $\bm{v}_0,\bm{v}_1,\bm{v}_2 \in H^1(\RTL)\times (0,\infty) $.
In the case 
\begin{align*}
\norm{P_{\gamma}v_2}_{H^1}
\ll& \min\{\norm{P_{\gamma}v_0}_{H^1},\norm{P_{\gamma}v_1}_{H^1}\},
\end{align*}
by the inequality
\begin{align*}
\norm{P_{\gamma}v_j}_{H^1} 
\lesssim &  \inf_{q\in \R} \norm{\tau_{q}P_{\gamma}v_2- P_{\gamma}v_j}_{H^1} , \quad (j=0,1)
\end{align*}
we have
\begin{align*}
\mathfrak{m}_{\delta}(\bm{v}_0,\bm{v}_1) 
\lesssim&  \norm{P_d(v_0-v_2)}_{H^1}+\norm{P_{\gamma}u_0}_{H^1}+|\log c_0-\log c_2|\\
&+\norm{P_{d}(v_1-v_2)}_{H^1}+\norm{P_{\gamma}u_1}_{H^1}+|\log c_1-\log c_2|\\
 \lesssim & \mathfrak{m}_{\delta}(\bm{v}_0,\bm{v}_2)+\mathfrak{m}_{\delta}(\bm{v}_2,\bm{v}_1).
\end{align*}
In the case 
\begin{align*}
&\min\{\norm{P_{\gamma}v_0}_{H^1},\norm{P_{\gamma}v_1}_{H^1}\} 
\lesssim \norm{P_{\gamma}v_2}_{H^1}, 
\end{align*}
by the equivalence between $\norm{\cdot}_{E}$ and $\norm{\cdot}_{H^1}$, we obtain 
\begin{align*}
\mathfrak{m}_{\delta}(\bm{v}_0,\bm{v}_1)
\lesssim & \inf_{q_0,q_1\in \R, j=0,1}\bigl( \norm{P_{\gamma}v_j-\tau_{q_j} P_{\gamma} v_2}_{H^1}+\sqrt{\delta}|q_j|\min\{ \phi_{\delta}(v_j),\phi_{\delta}(v_2)\}\\
&+\norm{\tau_{q_j}P_{\gamma}v_2-\tau_{q_j-q_{1-j}}P_{\gamma}v_{1-j}}_{H^1}+\sqrt{\delta}|q_{1-j}|\min\{\phi_{\delta}(v_2),\phi_{\delta}(v_{1-j})\}\bigr)\\
&+\norm{P_d(v_0-v_1)}_{E}+|\log c_0-\log c_1|\\
\lesssim & \mathfrak{m}_{\delta}(\bm{v}_0,\bm{v}_2)+\mathfrak{m}_{\delta}(\bm{v}_2,\bm{v}_1).
\end{align*}
Therefore, (ii) holds true.

Finally, we show the completeness of $\mathfrak{m}_{\delta}$.
Let $\{\bm{v}_n\}_n$ be a sequence in $H^1(\RTL)\times (0,\infty) $ with $\mathfrak{m}_{\delta}(\bm{v}_n,\bm{v}_m) \to 0$ as $n,m \to \infty$.
We show $\{\bm{v}_n\}_n$ has a convergent subsequence.
If $\{ \bm{v}_n\}_n$ has a subsequence $\{\bm{v}_{n_k}\}_k$ with $\norm{P_{\gamma}v_{n_k}}_{H^1(\RTL)} \to 0$ as $k \to \infty$, then $\{\bm{v}_{n_k}\}_k$ is a convergent sequence.
Hence, we assume 
\[\inf_{n} \norm{P_{\gamma}v_n}_{H^1}>c>0.\]
Since $\{\bm{v}_n\}_n$ is a Cauchy sequence in $\mathfrak{m}_{\delta}$, there exists a subsequence $\{\bm{v}_{n_k}\}_k \subset \{\bm{v}_n\}_n$ such that 
\[\mathfrak{m}_{\delta}(\bm{v}_{n_k},\bm{v}_{n_{k+1}}) \leq \frac{1}{2^{k+1}}.\]
Thus, there exist $q_k \in \R$ and $C>0$ such that
\begin{align}\label{lem-md-eq-1}
\norm{\tau_{q_k}P_{\gamma}v_{n_k}-P_{\gamma}v_{n_{k+1}}}_{H^1} + \frac{|q_k|c}{\sqrt{\delta}} \leq & \norm{\tau_{q_k}P_{\gamma}v_{n_k}-P_{\gamma}v_{n_{k+1}}}_{H^1} \notag \\
&+\frac{|q_k|\min\{ \norm{P_{\gamma}v_{n_k}}_{H^1}, \norm{P_{\gamma}v_{n_{k+1}}}_{H^1(\RTL)}\}}{\sqrt{\delta}}
\leq  \frac{C}{2^{k+1}}.
\end{align}
Let
\[p_k =\sum_{j \geq k}q_k.\]
Then, from \eqref{lem-md-eq-1} we have $p_k \to 0$ as $k \to \infty$ and
\[\norm{\tau_{p_k}P_{\gamma}v_{n_k}-\tau_{p_{k+1}}P_{\gamma}v_{n_{k+1}}}_{H^1} \leq \frac{C}{2^{n+1}}.\]
Therefore, the sequence $\{\tau_{p_k}P_{\gamma}v_{n_k}\}_k$ converges to an element $v_* \in H^1(\RTL)$ and $P_{\gamma}v_{n_k}$ also converge to $v_*$ in $H^1$.
Since any Cauchy sequence in $\mathfrak{m}_{\delta}$ has a convergent sequence in $\mathfrak{m}_{\delta}$, $\mathfrak{m}_{\delta}$ is complete.

\qed

In the following lemma, we show the Lipschitz continuity of the flow of the system \eqref{LZKeq-1}--\eqref{LZKeq-2} on the quasi-metric space $(H^1(\RTL),\mathfrak{m}_{\delta})$ and the estimate of the nonlinear term.
To prove the following lemma, we apply the argument in the proof of Lemma 3.2 in \cite{N S 2}.
In our equation, since the order of $\dot{\rho}-c$ is same as $v$, we can not show the statement of Lemma 3.2 in \cite{N S 2} directly.

\begin{lemma}\label{lem-est-md}
There exists $T^*, \delta^*>0$ such that for any $0<\delta<\delta^*$ and solutions $(\bm{v}_j,\rho_j)=(v_j,c_j,\rho_j) $ 
to the system $\eqref{LZKeq-1}$--$\eqref{LZKeq-2}$ given in Theorem \ref{thm-gwp-LZKeq}, we have
\begin{align}
&\sup_{|t|\leq T^*} \mathfrak{m}_{\delta}(\bm{v}_0(t),\bm{v}_1(t)) \lesssim \mathfrak{m}_{\delta}(\bm{v}_0(0),\bm{v}_1(0))\label{lem-est-md-eq-1}\\
&\sup_{|t|\leq T^*} \Bigl( \norm{P_d(v_0(t)-v_1(t)-e^{t\mathcal{A}}(v_0(0)-v_1(0)))}_{E}\notag \\
&+ \Bigl| \inf_{q \in \R, j=0,1}\bigl( \norm{P_{\gamma}v_j(t)-\tau_qP_{\gamma}v_{1-j}(t)}_{E}^2+\delta|q|^2\phi_{\delta}(v_{1-j}(t)) ^2\bigr)^{1/2}\notag \\
& - \inf_{q \in \R, j=0,1}\bigl( \norm{P_{\gamma}v_j(0)-\tau_qP_{\gamma}v_{1-j}(0)}_{E}^2+\delta|q|^2\phi_{\delta}(v_{1-j}(0))^2 \bigr)^{1/2} \Bigr| \Bigr)\notag\\
\lesssim & \delta^{\frac{1}{4}}\mathfrak{m}_{\delta}(\bm{v}_0(0),\bm{v}_1(0)),\label{lem-est-md-eq-2}
\end{align}
where the implicit constants do not depend on $\delta$ and solutions $(\bm{v}_j,\rho_j)$.
\end{lemma}

\proof
To prove Lemma \ref{lem-est-md}, we treat three cases.
Let $C_1,C_0>0$ be large positive numbers with $C_0 \ll C_1 \ll C_2$ and $(\bm{v}_j,\rho_j) $ be solutions to the system \eqref{LZKeq-1}--\eqref{LZKeq-2} given in Theorem \ref{thm-gwp-LZKeq}.


{\bf Case (I)} We consider the case $\norm{v_j(0)}_{H^1}+|c_j(0)-c^*|<C_1\delta$ $(j=0,1)$.
Then, for sufficiently small $\delta>0$ and any $a_1,a_2 \in (c^*-2C_1\delta,c^*+2C_1\delta )$,  we have
\[|\log a_1-\log a_2| \simeq |a_1-a_2|< 4C_1\delta.\]
For $q \in \R$ and $j \in \{0,1\}$, we define 
\[\rho^q_j(t)=q+c^*t+\int_0^t\chi_{\delta}(v_j(s),c_j(s))(\dot{\rho}_j(s)-c^*)ds.\]
Let $\zeta _{\rho,j}=\tau_{\rho}P_{\gamma}v_j$, $Q_c^{\rho}=\tau_{\rho}Q_c$, 
$\chi_{\delta}^{j}=\chi_{\delta}(v_j,c_j-c^*)$ and $N^{\rho,j}=\tau_{\rho}P_{\gamma}N(v_j,c_j,\rho_j)$.
Then, from Theorem \ref{thm-gwp-LZKeq}, $\zeta_{\rho^q_j,j}$ satisfies
\begin{align*}
\partial_t\zeta_{\rho^q_j,j}=&-\partial_x \Delta \zeta_{\rho^q_j,j} +\tau_{\rho^q_j}[P_{\gamma},\partial_x\mathbb{L}_{c^*}]v_j-\partial_x(2Q_{c^*}^{\rho^q_j}\zeta_{\rho^q_j,j})-\chi_{\delta}^j(\dot{\rho}_j-c^*)\partial_x \zeta_{\rho^q_j,j}+\chi_{\delta}^jN^{\rho^q_j,j}
\end{align*}
and
\begin{align*}
\sup_{|t|\leq T^*} \norm{\zeta_{\rho^q_j,j}}_{E} \lesssim \norm{v_j(0)}_{H^1},
\end{align*}
where $[A,B]=AB-BA$.
The difference $\zeta_{\rho^{q_0}_0,0}-\zeta_{\rho^{q_{1}}_{1},1}$ satisfies
\begin{align}\label{eq-est-md-1}
\partial_t(\zeta_{\rho^{q_0}_0,0}-\zeta_{\rho^{q_{1}}_{1},1})=&-\partial_x\Delta(\zeta_{\rho^{q_0}_0,0}-\zeta_{\rho^{q_{1}}_{1},1})+\tau_{\rho^{q_0}_0}[P_{\gamma},\partial_x\mathbb{L}_{c^*}]v_0-\tau_{\rho^{q_{1}}_{1}}[P_{\gamma},\partial_x\mathbb{L}_{c^*}]v_{1}\notag \\
&-2\partial_x(Q_{c^*}^{\rho^{q_0}_0}\zeta_{\rho^{q_{0}}_{0},0}-Q_{c^*}^{\rho^{q_{1}}_{1}}\zeta_{\rho^{q_{1}}_{1},1})-\chi_{\delta}^0(\dot{\rho}_0-c^*)\partial_x \zeta_{\rho^{q_0}_0,0}\notag \\
&+\chi_{\delta}^{1}(\dot{\rho}_{1}-c^*)\partial_x \zeta_{\rho^{q_{1}}_{1},1}+\chi_{\delta}^0N^{\rho^{q_0}_0,0}-\chi_{\delta}^{1}N^{\rho^{q_{1}}_{1},1}.
\end{align}
We estimate each term of \eqref{eq-est-md-1}.
By the simple calculation, we have
\begin{align*}
|\chi_{\delta}^0(t)-\chi_{\delta}^1(t)|\lesssim& \frac{|\norm{v_0(t)}_{H^1}-\norm{v_1(t)}_{H^1}|+|c_0(t)-c_1(t)|}{\delta}\lesssim \frac{\mathfrak{m}_{\delta}(\bm{v}_0(t),\bm{v}_1(t))}{\delta}\\
\norm{\partial_x(v_0(t))^2-\partial_x(v_1(t))^2}_{H^{-2}}
\lesssim & \norm{v_0(t)-v_1(t)}_{L^2}\max_{j=0,1}\norm{v_j(t)}_{H^1} \lesssim \delta \mathfrak{m}_{\delta}(\bm{v}_0(t),\bm{v}_1(t))\\
 |\dot{c}_0(t)-\dot{c}_1(t)| \lesssim& \mathfrak{m}_{\delta}(\bm{v}_0(t),\bm{v}_1(t))\\
 |\dot{\rho}_0(t)-c_0(t)-\dot{\rho}_1(t)+c_1(t)|\lesssim&  \mathfrak{m}_{\delta}(\bm{v}_0(t),\bm{v}_1(t)).
\end{align*}
Therefore, we have
\begin{align}\label{eq-est-md--1}
&\norm{P_d(v_0(t)-v_1(t)-e^{t\mathcal{A}}(v_0(0)-v_1(0)))}_{E}\notag\\
\lesssim & \Bigl| \int_0^t \delta \mathfrak{m}_{\delta}(\bm{v}_0(s),\bm{v}_1(s)) ds \Bigr| \leq |t| \delta \sup_{|s|\leq |t|}\mathfrak{m}_{\delta}(\bm{v}_0(s),\bm{v}_1(s)).
\end{align}
From the boundedness of the operator norm of $e^{t\partial_x\mathbb{L}_{c^*}}$, there exists $k^*>0$ such that for $|t|\leq T^*$
\begin{align*}
\norm{P_de^{t\mathcal{A}}(v_0(0)-v_1(0))}_{E}
 \lesssim & \norm{P_d(v_0(0)-v_1(0))}_{E}+\norm{\zeta_{q_0,0}(0)-\zeta_{q_1,1}(0)}_{L^2}\\
&+|q_0-q_1|\min_{j=0,1}\norm{v_j(0)}_{L^2}
\end{align*}
and
\begin{align}\label{eq-est-md-0}
\norm{P_d(v_0(t)-v_1(t))}_{E} \lesssim & \norm{P_d(v_0(0)-v_1(0))}_{E}+\norm{\zeta_{q_0,0}(0)-\zeta_{q_1,1}(0)}_{L^2}\notag \\
&+|q_0-q_1|\min_{j=0,1}\norm{v_j(0)}_{L^2}+ (1+|t|)\delta \sup_{|s|\leq |t|}\mathfrak{m}_{\delta}(\bm{v}_0(s),\bm{v}_1(s)).
\end{align}
By $P_{\gamma}=I-P_d$, we have
\begin{align}\label{eq-est-md-2}
& \norm{\tau_{\rho^{q_0}_0}[P_{\gamma},\partial_x\mathbb{L}_{c^*}]v_0-\tau_{\rho^{q_{1}}_{1}}[P_{\gamma},\partial_x\mathbb{L}_{c^*}]v_1}_{H^1} \notag \\
\lesssim &\norm{\zeta_{\rho^{q_0}_0,0}-\zeta_{\rho^{q_1}_1,1}}_{L^2}+|q_0-q_1|\min_{j=0,1}\norm{v_j}_{L^2}+\norm{P_d(v_0-v_1)}_{L^2}.
\end{align}
From the definition of the system \eqref{LZKeq-1}--\eqref{LZKeq-2}, the differences $\dot{\rho}_0-c_0-\dot{\rho}_1+c_1$ and $\dot{c}_0-\dot{c}_1$ satisfy
\begin{align}
|\dot{\rho}_0-c_0-\dot{\rho}_1+c_1| \lesssim & \norm{P_{d}(v_0-v_1)}_{L^2}+ \norm{\zeta_{\rho^{q_0}_0,0}-\zeta_{\rho^{q_1}_1,1}}_{L^2} + |q_0-q_1|\min_{j=0,1} \norm{v_j}_{L^2} \notag \\
&+ \delta|c_0(0)-c_1(0)|, \label{eq-est-md-3}
\\
|\dot{c}_0-\dot{c}_1| \lesssim  \delta\Bigl ( |c_0(0)-c_1(0)|& + \norm{P_{d}(v_0-v_1)}_{L^2}+  \norm{\zeta_{\rho^{q_0}_0,0}-\zeta_{\rho^{q_1}_1,1}}_{L^2} + |q_0-q_1|\min_{j=0,1} \norm{v_j}_{L^2}\Bigr).\label{eq-est-md-4}
\end{align}
Applying the similar estimate as \eqref{gwp-eq-1} and \eqref{gwp-eq-5}, for $0<T<T^*/4$ and $w_0,w_1 \in X^{1,b}$ with $\tau_{-\rho_j^{q_j}}w_j(t)=P_{\gamma}v_j(t)$ $(|t|\leq T^*, j=0,1)$ we obtain that
\begin{align}
&\norm{\partial_x(Q_{c^*}^{\rho^{q_0}_0}w_0-Q_{c^*}^{\rho^{q_{1}}_{1}}w_1)}_{X^{1,0}_T}\notag\\
\lesssim & \delta (|q_0-q_1| + |c_0(0)-c_1(0)| + \norm{\dot{\rho}_0-c_0-\dot{\rho}_1+c_1}_{L^2_T}+\norm{\dot{c}_0-\dot{c}_1}_{L^2_T}) + \norm{w_0-w_1}_{X^{1,b}_T} \label{eq-est-md-5}
\end{align}
and
\begin{align}
&\norm{\chi_{\delta}^0\tau_{\rho_0^{q_0}}\partial_x(\tau_{-\rho_0^{q_0}}w_0+P_{d}v_0)^2-\chi_{\delta}^1\tau_{\rho_1^{q_1}}\partial_x(\tau_{-\rho_1^{q_1}}w_1+P_{d}v_1)^2}_{X^{1,2b-\frac{3}{2}}_T} \notag\\
\lesssim_b& \delta \Bigl(\norm{w_0-w_1}_{X^{1,b}_T}+|c_0(0)-c_1(0)|+\norm{\dot{c}_0-\dot{c}_1}_{L^2_T}+\norm{\dot{\rho}_0-c_0-\dot{\rho}_1+c_1}_{L^2_T}\notag\\
&+|q_0-q_1|\sup_{|t|\leq T}\min_{j=0,1}\norm{v_j(t)}_{L^2} + \sup_{|t| \leq T}\norm{P_{d}(v_0(t)-v_1(t))}_{L^2}\Bigr). \label{eq-est-md-6}
\end{align}
By the equations \eqref{eq-est-md-0}--\eqref{eq-est-md-6}, from a priori estimate of \eqref{eq-est-md-1}  we obtain 
\begin{align}\label{eq-est-md-6-1}
\sup_{|t|\leq T^*} \norm{\zeta_{\rho^{q_{0}}_{0},0}(t)-\zeta_{\rho^{q_{1}}_{1},1}(t)}_{E} \lesssim & \norm{P_d(v_0(0)-v_1(0))}_{E}+ \norm{\zeta_{q_{0},0}(0)-\zeta_{q_{1},1}(0)}_{E}\notag \\
&+\delta |q_0-q_1|+\delta |c_0(0)-c_1(0)|+T^* \delta \sup_{|s|\leq |t|} \mathfrak{m}_{\delta}(\bm{v}_0(s),\bm{v}_1(s))
\end{align}
for small $T^*>0$ and $|t| <T^*$.
Therefore, we obtain \eqref{lem-est-md-eq-1} by 
\begin{align}\label{eq-est-md-7}
&\sup_{|t|\leq T^*} \mathfrak{m}_{\delta}(\bm{v}_0(t),\bm{v}_1(t)) ^2\notag\\
 \lesssim & 
\sup_{|t|\leq T^*} \Bigl( \norm{P_d(v_0(t)-v_1(t))}_{E}^2+ \inf_{ q \in \R, j=0,1}\Bigl (  \norm{\zeta_{\rho^{0}_j,j }(t)-\zeta_{\rho^{q}_{1-j},1-j}(t)}_{E}^2+\delta |\rho^{0}_{j}(t)-\rho^{q}_{1-j}(t)|^2 \Bigr)\notag \\
& \quad +|c_0(t)-c_1(t)|^2 \Bigr)\notag \\
\lesssim & \norm{P_d(v_0(0)-v_1(0))}_{E}^2+ \inf_{ q \in \R, j=0,1}\Bigl (  \norm{\zeta_{0,j }(0)-\zeta_{q,1-j}(0)}_{E}^2+\delta |q|^2 \Bigr) +|c_0(0)-c_1(0)|^2 \notag\\
\lesssim& \mathfrak{m}_{\delta}(\bm{v}_0(0),\bm{v}_1(0))^2.
\end{align}
The equations \eqref{eq-est-md--1} and \eqref{eq-est-md-7} yields 
\begin{align} \label{eq-est-md-7-1}
\norm{P_d(v_0(t)-v_1(t)-e^{t\mathcal{A}}(v_0(0)-v_1(0)))}_{E}
\lesssim & \delta \mathfrak{m}_{\delta}(\bm{v}_0(0),\bm{v}_1(0)).
\end{align}

For any $\epsilon >0$, there exist $q=q(\epsilon)>0$ and $j_0 \in \{0,1\}$ such that
\[
\inf_{p \in \R, j=0,1}\bigl( \norm{P_{\gamma}v_j(0)-\tau_pP_{\gamma}v_{1-j}(0)}_{E}^2+\delta|p|^2\phi_{\delta}(v_{1-j}(0))^2 \bigr)\geq \norm{\zeta_{0,j_0}(0)-\zeta_{q,1-j_0}(0)}_{E}^2+\delta|q|^2 -\epsilon.\]
To show the inequality \eqref{lem-est-md-eq-2}, we estimate the right hand side of following inequality.
\begin{align*}
&\inf_{p \in \R, j=0,1}\bigl( \norm{P_{\gamma}v_j(t)-\tau_pP_{\gamma}v_{1-j}(t)}_{E}^2+\delta|p|^2\phi_{\delta}(v_{1-j}(t)) ^2\bigr) \\
& - \inf_{p \in \R, j=0,1}\bigl( \norm{P_{\gamma}v_j(0)-\tau_pP_{\gamma}v_{1-j}(0)}_{E}^2+\delta|p|^2\phi_{\delta}(v_{1-j}(0))^2 \bigr) -\epsilon\\
\lesssim & \norm{\zeta_{0,j_0}(t) -\zeta_{\rho_{1-j_0}^{q}-\rho_{j_0}^{0},1-j_0}(t)}_{E}^2+ \delta|\rho_{j_0}^{0}-\rho_{1-j_0}^{q}|^2-\norm{\zeta_{0,j_0}(0)-\zeta_{q,1-j_0}(0)}_{E}^2-\delta|q|^2\\
\lesssim & \norm{P_d(\zeta_{0,j_0}(t) -\zeta_{\rho_{1-j_0}^{q}-\rho_{j_0}^{0},1-j_0}(t))}_{E}^2-\norm{P_d(\zeta_{0,j_0}(0)-\zeta_{q,1-j_0}(0))}_{E}^2\\
&+ \tbr{\mathbb{L}_{c^*}P_{\gamma}(\zeta_{0,j_0}(t) -\zeta_{\rho_{1-j_0}^{q}-\rho_{j_0}^{0},1-j_0}(t)),P_{\gamma}(\zeta_{0,j_0}(t) -\zeta_{\rho_{1-j_0}^{q}-\rho_{j_0}^{0},1-j_0}(t))}_{H^{-1},H^1}\\
&- \tbr{\mathbb{L}_{c^*}(\zeta_{0,j_0}(t) -\zeta_{\rho_{1-j_0}^{q}-\rho_{j_0}^{0},1-j_0}(t)),\zeta_{0,j_0}(t) -\zeta_{\rho_{1-j_0}^{q}-\rho_{j_0}^{0},1-j_0}(t)}_{H^{-1},H^1}\\
&+ \tbr{\mathbb{L}_{c^*}(\zeta_{0,j_0}(t) -\zeta_{\rho_{1-j_0}^{q}-\rho_{j_0}^{0},1-j_0}(t)),\zeta_{0,j_0}(t) -\zeta_{\rho_{1-j_0}^{q}-\rho_{j_0}^{0},1-j_0}(t)}_{H^{-1},H^1}\\
&-\tbr{\mathbb{L}_{c^*}(\zeta_{0,j_0}(0) -\zeta_{q,1-j_0}(0)),\zeta_{0,j_0}(0) -\zeta_{q,1-j_0}(0)}_{H^{-1},H^1}\\
&+\tbr{\mathbb{L}_{c^*}(\zeta_{0,j_0}(0) -\zeta_{q,1-j_0}(0)),\zeta_{0,j_0}(0) -\zeta_{q,1-j_0}(0)}_{H^{-1},H^1}\\
&-\tbr{\mathbb{L}_{c^*}P_{\gamma}(\zeta_{0,j_0}(0) -\zeta_{q,1-j_0}(0)),P_{\gamma}(\zeta_{0,j_0}(0) -\zeta_{q,1-j_0}(0))}_{H^{-1},H^1}\\
&+\delta( |\rho_{j_0}^{0}-\rho_{1-j_0}^{q}|^2-|q|^2).
\end{align*}
Since $|q| \delta ^{\frac{1}{2}}- \epsilon   \lesssim \mathfrak{m}_{\delta}(\bm{v}_0(0),\bm{v}_1(0))$, we have
\begin{align}\label{eq-est-md-8}
\norm{P_d\bigl(\zeta_{0,j_0}(t) -\zeta_{\rho_{1-j_0}^{q}-\rho_{j_0}^{0},1-j_0}(t)\bigr)}_{E}=&\norm{P_d\bigl(\zeta_{0,1-j_0}(t) -\zeta_{\rho_{1-j_0}^{q}-\rho_{j_0}^{0},1-j_0}(t)\bigr)}_{E}\notag \\
\lesssim &|\rho_{1-j_0}^{q}-\rho_{j_0}^{0}|\norm{\zeta_{0,1-j_0}(t)}_{H^1} \lesssim \delta^{\frac{1}{2}}\bigl(\mathfrak{m}_{\delta}(\bm{v}_0(0),\bm{v}_1(0))+\epsilon \bigr)
\end{align}
and
\begin{align}\label{eq-est-md-10}
&\delta( |\rho_{j_0}^{0}(t)-\rho_{1-j_0}^{q}(t)|^2-|q|^2) \notag \\
\lesssim & \delta\Bigl( |q|\int_0^t|\chi_{\delta}^0(\dot{\rho}_0-c^*)-\chi_{\delta}^1(\dot{\rho}_1-c^*)|ds+\Bigl(\int_0^t|\chi_{\delta}^0(\dot{\rho}_0-c^*)-\chi_{\delta}^1(\dot{\rho}_1-c^*)|ds \Bigr)^2 \Bigr)\notag \\
\lesssim& \delta^{\frac{1}{2}} \bigl(\mathfrak{m}_{\delta}(\bm{v}_0(0),\bm{v}_1(0))+\epsilon \bigr)^2.
\end{align}
The same calculation as \eqref{eq-est-md-8} yields
\begin{align}\label{eq-est-md-9}
\norm{P_d(\zeta_{0,j_0}(0)-\zeta_{q,1-j_0}(0))}_{E}\lesssim \delta^{\frac{1}{2}}\bigl(\mathfrak{m}_{\delta}(\bm{v}_0(0),\bm{v}_1(0))+\epsilon \bigr),
\end{align}
\begin{align}\label{eq-est-md-11}
&\Bigl| \tbr{\mathbb{L}_{c^*}P_{\gamma}(\zeta_{0,j_0}(t) -\zeta_{\rho_{1-j_0}^{q}-\rho_{j_0}^{0},1-j_0}(t)),P_{\gamma}(\zeta_{0,j_0}(t) -\zeta_{\rho_{1-j_0}^{q}-\rho_{j_0}^{0},1-j_0}(t))}_{H^{-1},H^1} \notag \\
&- \tbr{\mathbb{L}_{c^*}(\zeta_{0,j_0}(t) -\zeta_{\rho_{1-j_0}^{q}-\rho_{j_0}^{0},1-j_0}(t)),\zeta_{0,j_0}(t) -\zeta_{\rho_{1-j_0}^{q}-\rho_{j_0}^{0},1-j_0}(t)}_{H^{-1},H^1} \Bigr| \notag \\
\lesssim &\norm{P_d\bigl(\zeta_{0,j_0}(t) -\zeta_{\rho_{1-j_0}^{q}-\rho_{j_0}^{0},1-j_0}(t)\bigr)}_{E}^2\notag \\
\lesssim & \delta \bigl(\mathfrak{m}_{\delta}(\bm{v}_0(0),\bm{v}_1(0))+ \epsilon \bigr)^2
\end{align}
and
\begin{align}\label{eq-est-md-12}
&\Bigl|\tbr{\mathbb{L}_{c^*}(\zeta_{0,j_0}(0) -\zeta_{q,1-j_0}(0)),\zeta_{0,j_0}(0) -\zeta_{q,1-j_0}(0)}_{H^{-1},H^1}\notag \\
&-\tbr{\mathbb{L}_{c^*}P_{\gamma}(\zeta_{0,j_0}(0) -\zeta_{q,1-j_0}(0)),P_{\gamma}(\zeta_{0,j_0}(0) -\zeta_{q,1-j_0}(0))}_{H^{-1},H^1}\Bigr|\notag \\
\lesssim & \delta \bigl(\mathfrak{m}_{\delta}(\bm{v}_0(0),\bm{v}_1(0))+\epsilon \bigr)^2.
\end{align}
Next, we estimate 
\begin{align*}
 &\tbr{\mathbb{L}_{c^*}(\zeta_{0,j_0}(t) -\zeta_{\rho_{1-j_0}^{q}-\rho_{j_0}^{0},1-j_0}(t)),\zeta_{0,j_0}(t) -\zeta_{\rho_{1-j_0}^{q}-\rho_{j_0}^{0},1-j_0}(t)}_{H^{-1},H^1}\\
&-\tbr{\mathbb{L}_{c^*}(\zeta_{0,j_0}(0) -\zeta_{q,1-j_0}(0)),\zeta_{0,j_0}(0) -\zeta_{q,1-j_0}(0)}_{H^{-1},H^1}\\
=& \int_0^t \partial_s \tbr{\mathbb{L}_{c^*}(\zeta_{0,j_0}(s) -\zeta_{\rho_{1-j_0}^{q}-\rho_{j_0}^{0},1-j_0}(s)),\zeta_{0,j_0}(s) -\zeta_{\rho_{1-j_0}^{q}-\rho_{j_0}^{0},1-j_0}(s)}_{H^{-1},H^1}ds.
\end{align*}
Since
\begin{align*}
&\partial_t(\zeta_{0,j_0}-\zeta_{\rho_{1-j_0}^{q}-\rho_{j_0}^{0},1-j_0})\\
=&\partial_x\mathbb{L}_{c^*}(\zeta_{0,j_0}-\zeta_{\rho_{1-j_0}^{q}-\rho_{j_0}^{0},1-j_0})+[P_{\gamma},\partial_x\mathbb{L}_{c^*}]v_{j_0}-\tau_{\rho_{1-j_0}^{q}-\rho_{j_0}^{0}}[P_{\gamma},\partial_x\mathbb{L}_{c^*}]v_{1-j_0}\\
&-2\partial_x((Q_{c^*}-Q_{c^*}^{\rho_{1-j_0}^{q}-\rho_{j_0}^{0},1-j_0})\zeta_{\rho_{1-j_0}^{q}-\rho_{j_0}^{0},1-j_0})\\
&+\Bigl(\chi_{\delta}^{j_0}(\dot{\rho}_{j_0}-c^*)-\chi_{\delta}^{1-j_0}(\dot{\rho}_{1-j_0}-c^*)\Bigr)\partial_x \zeta_{\rho_{1-j_0}^{q}-\rho_{j_0}^{0},1-j_0}\\
&+\chi_{\delta}^{j_0}N^{0,j}-\chi_{\delta}^{1-j_0}N^{\rho_{1-j_0}^{q}-\rho_{j_0}^{0},1-j_0},\\
\end{align*}
by the similar calculation to \eqref{eq-gwp-11} in the proof of Theorem \ref{thm-gwp-LZKeq} we obtain 
\begin{align}\label{eq-est-md-13}
&\Bigl| \tbr{\mathbb{L}_{c^*}(\zeta_{0,j_0}(t) -\zeta_{\rho_{1-j_0}^{q}-\rho_{j_0}^{0},1-j_0}(t)),\zeta_{0,j_0}(t) -\zeta_{\rho_{1-j_0}^{q}-\rho_{j_0}^{0},1-j_0}(t)}_{H^{-1},H^1}\notag \\
&-\tbr{\mathbb{L}_{c^*}(\zeta_{0,j_0}(0) -\zeta_{q,1-j_0}(0)),\zeta_{0,j_0}(0) -\zeta_{q,1-j_0}(0)}_{H^{-1},H^1}\Bigr| \notag \\
\lesssim & \delta ( \mathfrak{m}_{\delta}(\bm{v}_0(0),\bm{v}_1(0))+ \varepsilon  )^2+ \sup_{|s|<|t|}\Bigl|\tbr{\mathbb{L}_{c^*}(\zeta_{0,j_0}(s) -\zeta_{\rho_{1-j_0}^{q}-\rho_{j_0}^{0},1-j_0}(s)), [P_{\gamma},\partial_x\mathbb{L}_{c^*}]v_{j_0}(s)\notag \\
&-\tau_{\rho_{1-j_0}^{q}-\rho_{j_0}^{0}}[P_{\gamma},\partial_x\mathbb{L}_{c^*}]v_{1-j_0}(s)}_{H^{-1},H^1}\Bigl|.
\end{align}
From \eqref{orth-pro-1} and \eqref{orth-pro-2}, we have 
\begin{align}\label{eq-est-md-c1-3}
\mathbb{L}_{c^*}[P_{\gamma},\partial_x \mathbb{L}_{c^*}]= - \mathbb{L}_{c^*}(P_1+P_2)\partial_x\mathbb{L}_{c^*}+\mathbb{L}_{c^*} \partial_x\mathbb{L}_{c^*}(P_1+P_2)=0.
\end{align}
The above equation yields
\begin{align}
&\norm{\mathbb{L}_{c^*}\Bigl([P_{\gamma},\partial_x\mathbb{L}_{c^*}]v_{j_0}-\tau_{\rho_{1-j_0}^{q}-\rho_{j_0}^{0}}[P_{\gamma},\partial_x\mathbb{L}_{c^*}]v_{1-j_0}\Bigr)}_{H^1}\notag \\
\leq &\norm{\mathbb{L}_{c^*}\Bigl( [P_{\gamma},\partial_x\mathbb{L}_{c^*}]v_{1-j_0}-\tau_{\rho_{1-j_0}^{q}-\rho_{j_0}^{0}}[P_{\gamma},\partial_x\mathbb{L}_{c^*}]v_{1-j_0}\Bigr)}_{H^1}+\norm{\mathbb{L}_{c^*}[P_{\gamma},\partial_x\mathbb{L}_{c^*}](v_{j_0}-v_{1-j_0})}_{H^1}\notag \\
\lesssim & \delta^{\frac{1}{2}}\mathfrak{m}_{\delta}(\bm{v}_0(0),\bm{v}_1(0)). \label{eq-est-md-13-1}
\end{align}
Thus, by the inequalities \eqref{eq-est-md-8}--\eqref{eq-est-md-13-1}, we obtain 
\begin{align}\label{eq-est-md-14} 
& \inf_{q \in \R, j=0,1}\bigl( \norm{P_{\gamma}v_j(t)-\tau_qP_{\gamma}v_{1-j}(t)}_{E}^2+\delta|q|^2\phi_{\delta}(v_{1-j}(t)) ^2\bigr) \notag \\
& - \inf_{q \in \R, j=0,1}\bigl( \norm{P_{\gamma}v_j(0)-\tau_qP_{\gamma}v_{1-j}(0)}_{E}^2+\delta|q|^2\phi_{\delta}(v_{1-j}(0))^2 \bigr) 
\lesssim  \delta^{\frac{1}{2}}\mathfrak{m}_{\delta}(\bm{v}_0(0),\bm{v}_1(0))^2
\end{align}
Showing the lower estimate of \eqref{eq-est-md-14} by reversing time, we complete the proof in Case (I).

{\bf Case (II)} We consider the case $\min_{j=0,1} \norm{v_j(0)}_{H^1}+|c_j(0)-c^*|>C_0\delta$.
In this case, since $\chi_\delta (v_j,c_j-c^*)=0$, the solution $(v_j,c_j,\rho_j)$ to system \eqref{LZKeq-1}--\eqref{LZKeq-2} is a solution to the linear system \eqref{Lsys-LZK-1}--\eqref{Lsys-LZK-2}.
Therefore, from the equation \eqref{linear-eq-1} we obtain $c_j(t)=c_j(0)$ and
\begin{align}\label{eq-est-md-c2--1}
&\norm{P_d(v_0(t)-v_1(t))}_{E}
\lesssim  \norm{P_de^{t\mathcal{A}}P_d(v_0(0)-v_1(0))}_{E}
\notag\\
&+\inf_{q \in \R, j=0,1}\Bigl( \norm{P_de^{t\mathcal{A}} (P_{\gamma}v_j(0)-\tau_qP_{\gamma}v_{1-j}(0))}_{E}+\norm{P_de^{t\mathcal{A}}(P_{\gamma}v_{1-j}(0)-\tau_qP_{\gamma}v_{1-j}(0))}_{E}  \Bigr)\notag\\
\lesssim & \mathfrak{m}_{\delta}(\bm{v}_0(0),\bm{v}_1(0))
\end{align}
for sufficiently small $t>0$.
By Lemma \ref{lem-linear-property}, we have
\begin{align}\label{eq-est-md-c2-0}
 \norm{P_{\gamma}v_j(t)}_{E}=\norm{P_{\gamma}v_j(0)}_{E} \mbox{ and } \phi_{\delta}(v_j(t))=\phi_{\delta}(v_j(0))
\end{align}
for $j=0,1$.
Let $\tilde{\zeta}_{q,j}(t)=P_{\gamma}v_j(t)-\tau_qP_{\gamma}v_{1-j}(t)$.
Then, $\tilde{\zeta}_{q,j}$ satisfies
\begin{align}\label{eq-est-md-c2-1}
\partial_t\tilde{\zeta}_{q,j}=&\partial_x\mathbb{L}_{c^*}\tilde{\zeta}_{q,j}+ [P_{\gamma},\partial_x\mathbb{L}_{c^*}]v_j-\tau_{q}[P_{\gamma},\partial_x\mathbb{L}_{c^*}]v_{1-j}+2\partial_x((\tau_qQ_{c^*}-Q_{c^*})\tau_qP_{\gamma}v_{1-j}).
\end{align}
Since $\tau_{q+c^*t} P_{\gamma} v_{1-j}$ is the solution to
\[w_t=-\partial_x\Delta w-2\partial_x((\tau_{q+c^*t}Q_{c^*})w)+\norm{\partial_xQ_{c^*}}_{L^2}^{-2}(w,\tau_{q+c^*t}\mathbb{L}_{c^*}\partial_x^2Q_{c^*})_{L^2}\tau_{q+c^*t}\partial_xQ_{c^*}\]
with the initial data $\tau_q P_{\gamma}v_{1-j}(0)$, we have there exists $T>0$ such that
\begin{align}\label{eq-est-md-c2-1-2-1}
\norm{ \tau_{q+c^*t} P_{\gamma}v_{1-j}}_{X^{1,b}_T} \lesssim \norm{P_{\gamma}v_{1-j}(0)}_{E}.
\end{align}
By the inequality \eqref{eq-linear-est-4} and \eqref{eq-est-md-c2-1-2-1} we obtain
\begin{align}\label{eq-est-md-c2-1-2}
\norm{\partial_x((\tau_qQ_{c^*}-Q_{c^*})\tau_qP_{\gamma}v_{1-j})}_{X^{1,0}_T}=&\norm{\tau_{c^*t}\partial_x((\tau_qQ_{c^*}-Q_{c^*})\tau_qP_{\gamma}v_{1-j})}_{X^{1,0}_T}\notag \\
 \lesssim &
|q| \norm{ \tau_{q+c^*t} P_{\gamma} v_{1-j}}_{X^{1,b}_T} 
\lesssim  |q| \norm{P_{\gamma}v_{1-j}(0)}_{E}
\end{align}
for $b>1/2$.
Combining the inequality \eqref{eq-est-md-c2-1-2} and the similar calculation to \eqref{eq-est-md-2} and \eqref{eq-est-md-6-1}, for small $t>0$ we obtain
\begin{align}\label{eq-est-md-c2-2}
\norm{\tilde{\zeta}_{q,j}(t)}_{E}\lesssim \norm{P_d(v_0(0)-v_1(0))}_{E}+\norm{\tilde{\zeta}_{q,j}(0)}_{E}+|q| \norm{P_{\gamma}v_{1-j}(0)}_{E}  .
\end{align}
The above inequalities and equations \eqref{eq-est-md-c2--1}--\eqref{eq-est-md-c2-2} yield
\begin{align*}
\mathfrak{m}_{\delta}(\bm{v}_0(t),\bm{v}_1(t))^2 =& \norm{P_d(v_0(t)-v_1(t))}_{E}^2+|\log c_0(0)-\log c_1(0)|^2\\
&+\inf_{q \in \R, j=0,1}\Bigl( \norm{\tilde{\zeta}_{q,j}(t)}_{E}^2 + \delta|q|^2\phi_{\delta}(v_{1-j}(t))^2 \Bigr)
\lesssim \mathfrak{m}_{\delta}(\bm{v}_0(0),\bm{v}_1(0))^2.
\end{align*}
The equation \eqref{eq-est-md-c2-0} yields
\begin{align}\label{eq-est-md-c2-2-1}
\norm{P_d\tilde{\zeta}_{j,q}(t)}_{E}=\norm{P_d\tau_qP_{\gamma}v_{1-j}(t)}_{E}\lesssim |q| \norm{P_{\gamma}v_{1-j}(0)}_{E}.
\end{align}
Since $[P_\gamma,\partial_x\mathbb{L}_{c^*}]=-P_1\partial_x\mathbb{L}_{c^*}P_{\gamma}$, we have
\begin{align}
& \norm{(1-\tau_{q})[P_{\gamma},\partial_x\mathbb{L}_{c^*}]v_{1-j}(t) }_{H^1} \notag \\
= & \norm{\partial_xQ_{c^*}}_{L^2}^{-2}\norm{(P_{\gamma}v_{1-j},\mathbb{L}_{c^*}\partial_x^2Q_{c^*})_{L^2}(1-\tau_q)\partial_xQ_{c^*}}_{H^1} 
 \lesssim  |q| \norm{P_{\gamma}v_{1-j}(0)}_{E}. \label{eq-est-md-c2-2-2}
 \end{align}
By the inequalities \eqref{eq-est-md-c2-1}, \eqref{eq-est-md-c2-1-2} and \eqref{eq-est-md-c2-2-2}, the equations \eqref{eq-est-md-c1-3} and the energy estimate, we obtain 
\begin{align}\label{eq-est-md-c2-3}
&\Bigl|\int_{0}^t\partial_t\tbr{\mathbb{L}_{c^*}\tilde{\zeta}_{q,j}(s),\tilde{\zeta}_{q,j}(s)}\, ds\Bigr| 
\lesssim   |q| \norm{P_{\gamma}v_{1-j}(0)}_{E}\norm{\tilde{\zeta}_{q,j}}_{L^{\infty}((-|t|,|t|)E)}.
\end{align}
Therefore, applying the calculation to show \eqref{eq-est-md-14}, by the equation \eqref{eq-est-md-c2-0} and the inequalities \eqref{eq-est-md-c2-1-2}, \eqref{eq-est-md-c2-2-1} and \eqref{eq-est-md-c2-3} we obtain 
\begin{align*}
& \inf_{q \in \R, j=0,1}\bigl( \norm{P_{\gamma}v_j(t)-\tau_qP_{\gamma}v_{1-j}(t)}_{E}^2+\delta|q|^2\phi_{\delta}(v_{1-j}(t)) ^2\bigr) \\
& - \inf_{q \in \R, j=0,1}\bigl( \norm{P_{\gamma}v_j(0)-\tau_qP_{\gamma}v_{1-j}(0)}_{E}^2+\delta|q|^2\phi_{\delta}(v_{1-j}(0))^2 \bigr)  
\lesssim  \delta^{\frac{1}{2}}\mathfrak{m}_{\delta}(\bm{v}_0(0),\bm{v}_1(0))^2
\end{align*}
By reversing time, we obtain \eqref{lem-est-md-eq-2} and complete the proof in Case (II).

{\bf Case (III)} We consider the case $ \norm{v_{1-j_1}(0)}_{H^1}+|c_{1-j_1}(0)-c^*|>C_1\delta \gg C_0 \delta >\norm{v_{j_1}(0)}_{H^1}+|c_{j_1}(0)-c^*| $.
By Theorem \ref{thm-gwp-LZKeq}, we have 
\[ \norm{v_{1-j_1}(t)}_{H^1}+|c_{1-j_1}(t)-c^*| \gtrsim C_1\delta \gg C_0 \delta \gtrsim \norm{v_{j_1}(t)}_{H^1}+|c_{j_1}(t)-c^*| \]
for sufficiently small $t>0$.
Therefore, 
\begin{align*}
\mathfrak{m}_{\delta}(\bm{v}_0(t),\bm{v}_1(t)) \simeq & \norm{v_{1-j_1}(t)}_{H^1}+|\log c_{1-j_1}(t)-\log c^*| \\
\simeq & \norm{v_{1-j_1}(0)}_{H^1}+|\log c_{1-j_1}(0)-\log c^*|
 \simeq  \mathfrak{m}_{\delta}(\bm{v}_0(0),\bm{v}_1(0)) \gtrsim \delta
\end{align*}
for sufficiently small $t>0$.
By the same argument as in Case (I), we obtain \eqref{eq-est-md-7-1}.
Let $\zeta_{q,j}=\tau_qP_{\gamma}v_j$.
Then, we have
\begin{align}
\partial_t(\zeta_{0,j_1}-\zeta_{q-\tilde{\rho}_{j_1}^0,1-j_1})=& \partial_x \mathbb{L}_{c^*}(\zeta_{0,j_1}-\zeta_{q-\tilde{\rho}_{j_1}^0,1-j_1}) +[P_{\gamma},\partial_x\mathbb{L}_{c^*}]v_{j_1}-\tau_{q-\tilde{\rho}_{j_1}^0}[P_{\gamma},\partial_x\mathbb{L}_{c^*}]v_{1-j_1} \notag \\
&+2\partial_x((\tau_{q-\tilde{\rho}_{j_1}^0}Q_{c^*}-Q_{c^*})\zeta_{q-\tilde{\rho}_{j_1}^0,1-j_1}) +\chi_{\delta}^{j_1}(\dot{\rho}_{j_1}-c^*)\partial_x\zeta_{0,j_1}+\chi_{\delta}^{j_1}N^{0,j_1}\label{eq-est-md-c3-1}
\end{align}
and
\begin{align}
\partial_t(\zeta_{0,1-j_1}-\zeta_{\tilde{\rho}_{j_1}^q,j_1})=& \partial_x \mathbb{L}_{c^*}(\zeta_{0,1-j_1}-\zeta_{\tilde{\rho}_{j_1}^q,j_1})+[P_{\gamma},\partial_x\mathbb{L}_{c^*}]v_{1-j_1}-\tau_{\tilde{\rho}_{j_1}^q}[P_{\gamma},\partial_x\mathbb{L}_{c^*}]v_{j_1} \notag \\
&+2\partial_x((\tau_{\tilde{\rho}_{j_1}^q}Q_{c^*}-Q_{c^*})\zeta_{\tilde{\rho}_{j_1}^q,j_1}) -\chi_{\delta}^{j_1}(\dot{\rho}_{j_1}-c^*)\partial_x\zeta_{\tilde{\rho}_{j_1}^q,j_1}-\chi_{\delta}^{j_1}N^{\tilde{\rho}_{j_1}^q ,j_1},\label{eq-est-md-c3-2}
\end{align}
where 
\[\tilde{\rho}^q_j(t)=q+\int_0^t\chi_{\delta}(v_j(s),c_j(s))(\dot{\rho}_j(s)-c^*)ds.\]
For any $\epsilon >0$, there exist $q=q(\epsilon)>0$ and $j_0 \in \{0,1\}$ such that
\begin{align*}
&\inf_{p \in \R, j=0,1}\bigl( \norm{P_{\gamma}v_j(0)-\tau_pP_{\gamma}v_{1-j}(0)}_{E}^2+\delta|p|^2\phi_{\delta}(v_{1-j}(0))^2 \bigr) \\
\geq & \norm{\zeta_{0,j_0}(0)-\zeta_{q,1-j_0}(0)}_{E}^2+\delta|q|^2\phi_{\delta}(v_{1-j_0}(0))^2 -\epsilon.
\end{align*}
By the smallness of $v_{j_1}$, we have
\begin{align*}
\delta|q|^2\phi_{\delta}(v_{1-j_0}(0))^2 - \epsilon \leq & \norm{P_{\gamma}v_{j_0}(0)-P_{\gamma}v_{1-j_0}(0)}_{E}^2-\norm{P_{\gamma}v_{j_0}(0)-\tau_{q}P_{\gamma}v_{1-j_0}(0)}_{E}^2\\
\lesssim& |q|\norm{P_{\gamma}v_{1-j_0}(0)}_{E}\mathfrak{m}_{\delta}(\bm{v}_0(0),\bm{v}_1(0)) 
\end{align*}
Thus, for  $0<\epsilon<\delta^3$ we obtain
\begin{align}\label{eq-lem-qq-1}
|q|\phi_{\delta}(v_{1-j_0}(0)) \lesssim \max \{\delta, \mathfrak{m}_{\delta}(\bm{v}_0(0),\bm{v}_1(0))\}\lesssim \mathfrak{m}_{\delta}(\bm{v}_0(0),\bm{v}_1(0)).
\end{align}
Since
\[\bigl||\tilde{\rho}_{j_0}^q|^2-|q|^2\bigr|+\bigl||q-\tilde{\rho}_{j_0}^0|^2-|q|^2\bigr| \lesssim \delta (|q|+\delta),\]
by the inequality \eqref{eq-lem-qq-1} we have the estimate of the mobile distance part 
\begin{align*}
&\bigl|\delta |\tilde{\rho}_{j_1}^q|^2\phi_{\delta}(v_{1-j_0}(t))^2-\delta |q|^2\phi_{\delta}(v_{1-j_0}(0))^2\bigr|+\bigl|\delta |q-\tilde{\rho}_{j_1}^0|^2\phi_{\delta}(v_{1-j_0}(t))^2-\delta |q|^2\phi_{\delta}(v_{1-j_0}(0))^2\bigr|\\
\lesssim & \delta^2 (|q|+\delta) \Bigl(\frac{\norm{P_{\gamma}v_{1-j_0}(0)}_E}{\delta}+1\Bigr)\phi_{\delta}(v_{1-j_0}(0))\\
\lesssim & \delta\mathfrak{m}_{\delta}(\bm{v}_0(0),\bm{v}_1(0))^2.
\end{align*}
Thus, from the energy estimate for \eqref{eq-est-md-c3-1} and \eqref{eq-est-md-c3-2}, the arguments in Case (I) and Case (II) yield \eqref{lem-est-md-eq-2}.

\qed

\section{Construction of the center stable manifolds}
In this section, we construct the center stable manifolds by applying the Hadamard method in \cite{N S 2}.
Let $\mathcal{H}$ be the complete quasi-metric space $H^1(\RTL) \times (0,\infty)$ with the quasi-distance $\mathfrak{m}_{\delta}$.
We denote $\mathscr{G}_{l,\delta}^{+}$ by 
\begin{align*}
\{G:\mathcal{H} \to P_+ H^1(\RTL);& \, G=G \circ P_{\leq 0}, G(0,c^*)=0, \\
& \norm{G(\bm{v}_0)-G(\bm{v}_1)}_{E}\leq l  \mathfrak{m}_{\delta}(\bm{v}_0,\bm{v}_1)\mbox{ for } \bm{v}_0,\bm{v}_1 \in \mathcal{H}\},
\end{align*}
where $P_{\leq 0} (v,c)= ((I-P_+)v,c)$.
We define the graph $\gbr{G}$ of $G \in \mathscr{G}_{l,\delta}^+$ as
\[\{ (v,c) \in \mathcal{H};P_+ v= G(v,c)\}.\]

In the following lemma, we show the upper estimate of the growth of unstable eigenmode.
\begin{lemma}\label{lem-Lip-l}
There exist $0< T^*<1$ and $ C_L>0$ such that if $l, \delta >0$  satisfy
\begin{align}\label{ass-smallness}
\delta+l \ll 1 \mbox{ and } l^{-1}\delta^{\frac{1}{4}} \ll  1 ,
\end{align}
then for any solutions $(v_j, c_j, \rho_j)$ to the system $\eqref{LZKeq-1}$--$\eqref{LZKeq-2}$ $(j=0,1)$ satisfying
\[ \norm{P_+(v_0(0)-v_1(0))}_{E} \leq l \mathfrak{m}_{\delta}(\bm{v}_0(0),\bm{v}_1(0))\]
one has
\begin{align}\label{lem-eq-2}
\norm{P_+(v_0(t) -v_1(t))}_{E} \leq 
\begin{cases}
C_Ll \mathfrak{m}_{\delta}(\bm{v}_0(t),\bm{v}_1(t)), & |t|\leq T^*, \\
l \mathfrak{m}_{\delta}(\bm{v}_0(t),\bm{v}_1(t)), &  - T^* \leq t \leq -\frac{T^*}{2}.
\end{cases}
\end{align}
\end{lemma}
\proof
By the boundedness of the operator $e^{t\partial_x\mathbb{L}_{c^*}}$ on $H^1(\RTL)$, we have for $t \in \R$
\begin{align}
&\min\{e^{\pm k^*t},e^{\pm k_*t}\}\norm{P_{\pm}(v_0(0)-v_1(0))}_{E} \notag \\
\leq & \norm{P_{\pm}e^{t\mathcal{A}}(v_0(0)-v_1(0))}_{E}
\leq \max\{e^{\pm k^*t},e^{\pm k_*t}\}\norm{P_{\pm}(v_0(0)-v_1(0))}_{E},\label{eq-Lip-l-1}
\end{align}
where $k^*$ and $k_*$ are defined by \eqref{min-max-eigen}.
From Lemma \ref{lem-linear-property}, we have
\begin{align}
(e^{t\mathcal{A}}(v_0(0)-v_1(0)), \partial_xQ_{c^*})_{L^2}
= ( v_0(0)-v_1(0),\partial_xQ_{c^*})_{L^2}
\label{eq-Lip-l-2}
\end{align}
and
\begin{align}\label{eq-Lip-l-3}
(e^{t\mathcal{A}}(v_0(0)-v_1(0)), Q_{c^*})_{L^2}
 =& ( v_0(0)-v_1(0),Q_{c^*})_{L^2}.
\end{align}
By the inequality \eqref{eq-Lip-l-1}, Lemma \ref{lem-est-md} yields that 
\begin{align}\label{eq-Lip-l-5-0}
\norm{P_+(v_0(t)-v_1(t))}_{E} \leq & \norm{P_+e^{t\partial_x\mathbb{L}_{c^*}}(v_0(0)-v_1(0))}_{E} + \delta^{\frac{1}{4}} C\mathfrak{m}_{\delta}(\bm{v}_0(0),\bm{v}_1(0))\notag \\
\leq & (\max\{e^{k^*t},e^{k_*t}\} l+\delta^{\frac{1}{4}}C)\mathfrak{m}_{\delta}(\bm{v}_0(0),\bm{v}_1(0))
\end{align}
and
\begin{align}
&\mathfrak{m}_{\delta}(\bm{v}_0(t),\bm{v}_1(t))^2-|\log c_0(t)- \log c_1(t)|^2 \notag \\
\geq & \norm{P_d(e^{t\mathcal{A}}(v_0(0)-v_1(0)))}_{E}^2\notag \\ 
&+\inf_{q\in \R, j=0,1}(\norm{P_{\gamma}v_j(0)-\tau_qP_{\gamma}v_{1-j}(0)}_{E}^2+\delta|q|^2\phi_{\delta}(v_{1-j}(0))^2)-C\delta^{\frac{1}{2}}\mathfrak{m}_{\delta}(\bm{v}_0(0),\bm{v}_1(0))^2\label{eq-Lip-l-5}
\end{align}
for sufficiently small $|t|$.
Plugging \eqref{eq-Lip-l-1}--\eqref{eq-Lip-l-3} into the estimate \eqref{eq-Lip-l-5}, we have there exist $C , T>0$ such that
\begin{align}\label{eq-Lip-l-6}
\mathfrak{m}_{\delta}(\bm{v}_0(t),\bm{v}_1(t))^2 
 \geq & 
\begin{cases}
(1-l^2+e^{2k^*t}l^2-C\delta^{\frac{1}{2}})\mathfrak{m}_{\delta}(\bm{v}_0(0),\bm{v}_1(0))^2, & -T\leq t \leq 0,\\
(e^{-2k^*|t|}-C\delta^{\frac{1}{2}})\mathfrak{m}_{\delta}(\bm{v}_0(0),\bm{v}_1(0))^2, & |t|\leq T.
\end{cases}
\end{align}
From \eqref{eq-Lip-l-5-0} and \eqref{eq-Lip-l-6} we obtain for sufficiently small $|t|$
\begin{align*}
\norm{P_+(v_0(t)-v_1(t))}_{E} \leq & (e^{k^*|t|}l + C \delta^{\frac{1}{4}}) (e^{k^*|t|} + C \delta^{\frac{1}{4}})\mathfrak{m}_{\delta}(\bm{v}_0(t),\bm{v}_1(t)) \notag \\
\leq & 2le^{2k^*|t|}\mathfrak{m}_{\delta}(\bm{v}_0(t),\bm{v}_1(t)).
\end{align*}
The above inequality yields \eqref{lem-eq-2} in the case with sufficiently small $|t|$.
By the inequalities \eqref{eq-Lip-l-5-0} and \eqref{eq-Lip-l-6}, we have
\begin{align*}
\norm{P_+(v_0(t)-v_1(t))}_{E} 
\leq & l ( 1- k_*T/3+C\delta^{\frac{1}{4}} l^{-1})(1+C(l+\delta^{\frac{1}{4}})) \mathfrak{m}_{\delta}(\bm{v}_0(t),\bm{v}_1(t))\notag \\
\leq & l  \mathfrak{m}_{\delta}(\bm{v}_0(t),\bm{v}_1(t))
\end{align*}
for $-T\leq t\leq -T/2$ and sufficiently small $T$ and $\delta$.
This is the inequality \eqref{lem-eq-2} in the case for $-T^*\leq t\leq -T^*/2$.
\qed

The following lemma shows that the flow map $U_{\delta}(t)$ of the system \eqref{LZKeq-1}--\eqref{LZKeq-2} given by Theorem \ref{thm-gwp-LZKeq} yields the mapping on the set of graphs.
\begin{lemma}\label{lem-flow}
Under the condition $\eqref{ass-smallness}$, the solution map $U_{\delta}(t)$ of the system $\eqref{LZKeq-1}$--$\eqref{LZKeq-2}$ for $|t|\leq T^*$ defines a map $\mathcal{U}_{\delta}(t): \mathscr{G}_{l,\delta}^+ \to \mathscr{G}_{C_Ll,\delta}^+$ uniquely by the relation $U_{\delta}(t)(\gbr{G}\times \R)=\gbr{\mathcal{U}_{\delta}(t)G}\times \R$. 
Moreover, if $-T^*\leq t\leq -T^*/2$, then $\mathcal{U}_{\delta}(t)$ maps $\mathscr{G}_{l,\delta}^+$ into itself.
\end{lemma}
Lemma \ref{lem-flow} follows Lemma \ref{lem-Lip-l} and the similar proof to the proof of Lemma 3.4 in \cite{N S 2}.

Let 
\[ \norm{G}_{\mathscr{G}^+}=\sup_{\bm{v} \in \mathcal{H}\setminus \{(0,c^*)\}} \frac{\norm{G(\bm{v})}_{E}}{\norm{\bm{v}}_{E}},\]
where
\[\norm{(v,c)}_{E}^2=\norm{v}_{E}^2+|\log c-\log c^*|^2.\]
Then, for $G \in \mathscr{G}_{l,\delta}^{+}$, we have 
\[ \norm{G(\bm{v})}_{E} \leq l \mathfrak{m}_{\delta}(\bm{v},(0,c^*)) \leq l \norm{\bm{v}}_{E} \mbox{ for } \bm{v} \in \mathcal{H}.\]
Therefore, $G \in \mathscr{G}_{l,\delta}^+$ satisfies $\norm{G}_{\mathscr{G}^+}\leq l$.
By the definition of $\norm{\cdot}_{\mathscr{G}^+}$, we obtain the ordered pair $(\mathscr{G}_{l,\delta}^+,\norm{\cdot}_{\mathscr{G}^+})$ is the bounded complete metric space.

In the following lemma, we show the mapping $\mathcal{U}_{\delta}(t)$ in Lemma \ref{lem-flow} is a contraction.
\begin{lemma}\label{lem-cotr}
Under the condition $\eqref{ass-smallness}$, the mapping $\mathcal{U}_{\delta}(t)$ is a contraction on $(\mathscr{G}_{l,\delta}^+,\norm{\cdot}_{\mathscr{G}^+})$ for $t < - T^*/2$.
\end{lemma}
\proof
Let $G_0,G_1 \in \mathscr{G}_{l,\delta}^+$ and $  T \in [-T^*,-T^*/2]$.
We define solutions $(v_j,c_j,\rho_j)$ to the system \eqref{LZKeq-1}--\eqref{LZKeq-2} by 
\[(v_j(t), c_j(t),\rho_j(t))=U_{\delta}(t-T)(P_{\leq 0} \psi+(\mathcal{U}_{\delta}(T)G_j)(\psi,\alpha),\alpha ,q)\]
for $j\in \{0,1\}$,  $(\psi,\alpha) \in \mathcal{H}, t \in \R$ and $q \in \R$.
Then, we have
\begin{align}\label{eq-contr-1}
P_{\leq 0}v_0(T)=P_{\leq 0}v_1(T) \mbox{ and } P_+(v_0(T)-v_1(T))=v_0(T)-v_1(T).
\end{align}
The equality
\begin{align*}
 P_+(P_{\leq 0} \psi +\mathcal{U}_{\delta}(T)G_j(\psi,\alpha))=\mathcal{U}_{\delta}(T)G_j(P_{\leq 0} \psi,\alpha) 
 =& \mathcal{U}_{\delta}(T)G_j\bigl(P_{\leq 0} \psi + \mathcal{U}_{\delta}(T) G_j(\psi, \alpha),\alpha\bigr)
\end{align*}
implies 
\[P_{\leq 0} \psi + \mathcal{U}_{\delta}(T)G_j(\psi,\alpha) \in \gbr{\mathcal{U}_{\delta}(T)G_j}.\]
Thus, we have
\[U_{\delta}(-T)(P_{\leq 0} \psi + \mathcal{U}_{\delta}(T) G_j (\psi,\alpha ) ) \in U_{\delta}(-T)\gbr{\mathcal{U}_{\delta}(T)G_j}=\gbr{G_j}.\]
This inclusion yields 
\begin{align*}
P_+v_j(0)=&P_+(U_{\delta}(-T)(P_{\leq 0}\psi + \mathcal{U}_{\delta}(T)G_j(\psi,\alpha)))\notag \\
=&G_j\bigl(U_{\delta}(-T)(P_{\leq 0}\psi+\mathcal{U}_{\delta}(T)G_j(\psi,\alpha)),\alpha\bigr)=G_j(v_j(0),\alpha).
\end{align*}
Therefore, we have
\begin{align}\label{eq-contr-6}
\norm{P_{+}(v_0(0)-v_1(0))}_{E} 
\leq & \norm{G_0-G_1}_{\mathscr{G}^+} \norm{(P_{\leq 0} v_0(0),c_0(0))}_{E}+l \mathfrak{m}_{\delta}(P_{\leq 0} \bm{v}_0(0),P_{\leq 0} \bm{v}_1(0)),
\end{align}
where $\bm{v}_j(t)=(v_j(t),c_j(t))$ for $j=0,1$.
Since $c_0$ satisfies \eqref{LZKeq-2}, if $\max\{|c_0(0)-c^*|,|c_0(T)-c^*|\}>\sqrt{2}\delta$, then $c_0(0)=c_0(T)$. 
Thus, Theorem \ref{thm-gwp-LZKeq} yields that there exists $C>0$ such that
\begin{align}
&\norm{(P_{\leq 0}v_0(0), c_0(0))}_{E}^2 \notag \\
\leq & \norm{(P_-+P_0)(v_0(0)-e^{-T\mathcal{A}}v_0(T))}_{E}^2 + \norm{(P_-+P_0)e^{-T\mathcal{A}}v_0(T)}_{E}^2 \notag \\
&+ 2\norm{(P_-+P_0)(v_0(0)-e^{-T\mathcal{A}}v_0(T))}_{E}   \norm{(P_-+P_0)e^{-T\mathcal{A}}v_0(T)}_{E} \notag \\
&+ \norm{P_{\gamma}v_0(0)}_{E}^2-\norm{P_{\gamma}v_0(T)}_{E}^2+\norm{P_{\gamma}v_0(T)}_{E}^2+|\log c_0(0)-\log c_0(T)|^2 \notag\\
&+ 2|\log c_0(0)-\log c_0(T)||\log c_0(T)-\log c^*| +|\log c_0(T)-\log c^*|^2\notag\\
 \leq& (1+2C\delta)\norm{(P_{\leq 0}v_0(T),c_0(T))}_{E}^2+2C\delta \norm{P_+v_0(T)}_{E}^2. \notag 
\end{align}
Since $\norm{P_+v_0(T)}_{E} = \norm{\mathcal{U}_{\delta}(T)G_0(\psi,\alpha)}_{E}\leq l \norm{(\psi,\alpha)}_{E}$,
by the definition $(v_0(T),c_0(T))$ we have
\begin{align}\label{eq-contr-7}
\norm{(P_{\leq 0}v_0(0), c_0(0))}_{E}^2\leq (1+ 4C\delta) \norm{(\psi,\alpha)}_{E}^2.
\end{align}
Applying Lemma \ref{lem-est-md} from $t=T$, by the equation \eqref{eq-contr-1} we obtain
\begin{align}\label{eq-contr-2}
&\norm{P_d(v_0(t)-v_1(t)-e^{(t-T) \mathcal{A}} (v_0(T)-v_1(T)))}_{E} \notag \\
&+ \inf_{p \in \R, j=0,1}\bigl(\norm{P_{\gamma}v_j(t)-\tau_pP_{\gamma}v_{1-j}(t)}_{E}+\delta^{1/2}|p|\phi_{\delta}(v_{1-j}(t))\bigr) \lesssim \delta^{\frac{1}{4}}\norm{P_+(v_0(T)-v_1(T))}_{E}
\end{align}
for $T\leq t \leq 0$.
By the inequalities \eqref{eq-Lip-l-1} and \eqref{eq-contr-2}, we have that there exists $C>0$ such that 
\begin{align}
\norm{P_+(v_0(0)-v_1(0))}_{E} & \geq \norm{e^{-T\mathcal{A}}P_+(v_0(T)-v_1(T))}_{E} -C\delta^{\frac{1}{4}}\norm{P_+(v_0(T)-v_1(T))}_{E} \notag \\
&\geq (e^{-k_* T}-C\delta^{\frac{1}{4}}) \norm{P_+(v_0(T)-v_1(T))}_{E}\label{eq-contr-3}
\end{align}
and
\begin{align}
&\inf_{p \in \R, j=0,1}\bigl(\norm{P_{\gamma}v_j(0)-\tau_pP_{\gamma}v_{1-j}(0)}_{E}+\delta^{\frac{1}{2}}|p|\phi_{\delta}(v_{1-j}(0))\bigr) + \norm{P_-(v_0(0)-v_1(0))}_{E} \notag \\
&+\norm{P_0(v_0(0)-v_1(0))}_{E} \lesssim \delta^{\frac{1}{4}}\norm{P_+(v_0(T)-v_1(T))}_{E}. \label{eq-contr-4}
\end{align}
Therefore, by \eqref{eq-contr-4} we have 
\begin{align}\label{eq-contr-8}
\mathfrak{m}_{\delta}(P_{\leq 0} \bm{v}_0(0),P_{\leq 0} \bm{v}_1(0)) \lesssim \delta^{\frac{1}{4}}\norm{P_+(v_0(T)-v_1(T))}_{E}.
\end{align}
From the inequalities \eqref{eq-contr-6}, \eqref{eq-contr-7}, \eqref{eq-contr-3} and \eqref{eq-contr-8}, we obtain there exists $0<\Lambda<1$ such that
\begin{align*}
\norm{P_+(v_0(T)-v_1(T))}_{E}\leq & (e^{-k_*T}-C\delta^{\frac{1}{4}})^{-1}(1-Cl\delta^{\frac{1}{4}})^{-1} (1+4C\delta)^{1/2}\norm{G_0-G_1}_{\mathscr{G}^+} \norm{(\psi,\alpha)}_{E}\\
\leq & \Lambda \norm{G_0-G_1}_{\mathscr{G}^+} \norm{(\psi,\alpha)}_{E}.
\end{align*}
Therefore, 
\[\frac{\norm{\mathcal{U}_{\delta}(T)G_0(\psi,\alpha)-\mathcal{U}_{\delta}(T)G_1(\psi,\alpha)}_{E}}{\norm{(\psi,\alpha)}_{E}}=\frac{\norm{P_+(v_0(T)-v_1(T))}_{E}}{\norm{(\psi,\alpha)}_{E}}\leq \Lambda \norm{G_0-G_1}_{\mathscr{G}^+}.\]
Thus, $\mathcal{U}_\delta(T)$ is a contraction.
\qed

Applying Lemma \ref{lem-cotr}, we obtain the existence of the fix point of $\mathcal{U}_{\delta}(t)$.
\begin{proposition}\label{prop-fix-point}
Assume that $l,\delta>0$ satisfy $\eqref{ass-smallness}$.
Then there exists a unique $G_+^{\delta} \in \mathscr{G}_{l,\delta}^+$ such that $\mathcal{U}_{\delta}(t)G_+^{\delta} =G_+^{\delta} $ for all $t<0$.
Moreover, the uniqueness holds for any fixed $t<0$.
\end{proposition}
The proof of Proposition follows the proof of Theorem 3.6 in \cite{N S 2}.
In the following lemma, we show the estimate of the smallness of the modulation parameter $c$.
\begin{lemma}\label{lem-mod-para}
Let $c^*/2<c_0<2c^*$. 
Then, there exists $\delta>0$ such that for $u \in L^2(\RTL)$ satisfying $u=\tau_\rho(v+Q_c) $, $(v,Q_{c^*})_{L^2}=0$, $\norm{v}_{L^2}<\delta$ and $\norm{u}_{L^2}=\norm{Q_{c_0}}_{L^2}$, we have
\[|c_0-c|\lesssim \norm{v}_{L^2}^2+|c_0-c^*|\norm{v}_{L^2}.\]
\end{lemma}
\proof
Since
\begin{align*}
\norm{Q_{c_0}}_{L^2}^2-\norm{Q_c}_{L^2}^2=&\norm{v}_{L^2}^2+2(v,Q_c-Q_{c^*})_{L^2},
\end{align*}
we have
\[|c_0-c| \lesssim \norm{v}_{L^2}^2+(|c_0-c|+|c_0-c^*|)\norm{v}_{L^2}.\]
The smallness of $\norm{v}_{L^2}$ yields the conclusion.
\qed

Let 
\begin{align*}
\mathcal{M}_{cs}^{\delta} (c^*,r)=\{ \tau_{\rho}(w+G_+^{\delta} (w,c)+Q_{c}); &w \in (P_-+P_{\gamma})H^1(\RTL),|c-c^*|<c^*/2,\\
& \inf_{q \in \R}\norm{w+G_+^{\delta} (w,c)+Q_{c}-\tau_qQ_{c^*}}_{H^1}< r, \rho \in \R\}
\end{align*}
and 
\begin{align*}
\tilde{\mathcal{M}}_{cs}^{\delta} (c^*,r)=\{ \tau_{\rho}(w+G_+^{\delta} (w,c)+Q_{c}); &w \in P_{\leq 0}H^1(\RTL), \norm{P_0w}_{H^1(\RTL)}<r^{1/2},\\
|c-c^*|<c^*/2,& \inf_{q \in \R}\norm{w+G_+^{\delta} (w,c)+Q_{c}-\tau_qQ_{c^*}}_{H^1}< r, \rho \in \R\}
\end{align*}
for $r>0$.

We show the stability of $Q_c$ on $\tilde{\mathcal{M}}_{cs}^{\delta}(c,\varepsilon )$.
\begin{theorem}\label{thm-ex-csm}
Let $l,\delta>0$. Assume $\eqref{ass-smallness}$.
For any $\varepsilon >0$, there exists $\tilde{\varepsilon}=\tilde{\varepsilon}(c^*,\varepsilon) >0$ such that for $u_0 \in \tilde{\mathcal{M}}_{cs}^{\delta} (c^*,\tilde{\varepsilon})$ the solution $u$ to the equation $\eqref{ZKeq}$ with the initial data $u_0$ satisfies $u(t) \in  \tilde{\mathcal{M}}_{cs}^{\delta} (c^*,\varepsilon)$ 
for all $t>0$.
\end{theorem}
\proof
Let $l,\delta>0$ satisfying $\eqref{ass-smallness}$.
We prove the stability by contradiction.
We assume there exists $0<\varepsilon _0 \ll \delta^2$ such that for $0<\tilde{\varepsilon}<\varepsilon _0$ there exist $t_0>0$ and the solution $u$ to the equation \eqref{ZKeq} with the initial data $ \tau_{\rho_0(0)}(v_0(0)+G_+^{\delta} (v_0(0),c_0(0))+Q_{c_0(0)}) \in \tilde{\mathcal{M}}_{cs}^{\delta} (c^*,\tilde{\varepsilon})$  satisfying 
\[\sup_{0\leq t \leq t_0} \inf_{q\in \R}\norm{u(t)-\tau_qQ_{c^*}}_{H^1}\leq \varepsilon_0\]
and 
\[\inf_{q\in \R}\norm{u(t_0)-\tau_qQ_{c^*}}_{H^1}= \varepsilon_0.\]
We define the solution $(v_0(t),c_0(t),\rho_0(t))$ to the system \eqref{LZKeq-1}--\eqref{LZKeq-2} with the initial data $(v_0(0)+G_+^{\delta} (v_0(0),c_0(0)), c_0(0),\rho_0(0))$ and the solution $(v_1(t),c_1(t),\rho_1(t))$ to the system \eqref{LZKeq-1}--\eqref{LZKeq-2} with the initial data $(v_1(0), c_1(0),\rho_1(0))\in (I-P_0)H^1(\RTL) \times (0, \infty) \times \R$ satisfying $\tau_{\rho_0(0)}(v_0(0)+G_+^{\delta} (v_0(0),c_0(0))+Q_{c_0(0)})=\tau_{\rho_1(0)}(v_1(0)+Q_{c_1(0)})$.
Then, by the invariance of $U_{\delta}(t)(\gbr{G_+^{\delta} }\times \R )=\gbr{G_+^{\delta} }\times \R $, we have $G_+^{\delta} (v_0(t),c_0(t))=P_+v_0(t)$.
Since 
\begin{align*}
 \norm{\tau_{\rho_0}Q_{c_0}-\tau_qQ_{c^*}}_{H^1} \lesssim & \norm{P_0(Q_{c_0}-\tau_{q-\rho_0}Q_{c^*})}_{H^1}\\
 \lesssim & \norm{P_0v_0}_{H^1}+\norm{\tau_{\rho_0}(v_0+G_+^{\delta}(v_0,c_0)+Q_{c_0})-\tau_qQ_{c^*}}_{H^1},
 \end{align*}
 we have
\begin{align*}
\norm{v_0(t)+G_+^{\delta}(v_0(t),c_0(t))}_{H^1} +|c_0(t)-c^*| \lesssim & \inf_{q\in \R}\norm{u(t)-\tau_qQ_{c^*}}_{H^1}+\norm{P_0v_0(0)}_{H^1} \\
\lesssim & \varepsilon_0 +\tilde{\varepsilon }^{1/2} \ll \delta.
\end{align*}
The smallness of $v_0+G_+^{\delta}(v_0,c_0)$ and $c_0-c^*$  yields 
\[u(t)=\tau_{\rho_0(t)}(v_0(t)+G_+^{\delta}(v_0(t),c_0(t))+Q_{c_0(t)})\]
for $0\leq t \leq t_0$.
Form the similar calculation, we obtain 
\begin{align*}
\norm{v_1(0)}_{H^1} \lesssim & \inf_{q \in \R} (\norm{v_1(0)-Q_{c_1(0)}+\tau_qQ_{c^*}}_{H^1}+\norm{Q_{c_1}(0)-\tau_qQ_{c^*}}_{H^1}) \\
\lesssim & \inf_{q \in \R} \norm{v_1(0)-Q_{c_0(0)}+\tau_qQ_{c^*}}_{H^1}<\tilde{\varepsilon } \lesssim \delta^2.
\end{align*}
Therefore, we have
\[\tau_{\rho_1(t)}(v_1(t)+Q_{c_1(t)})=\tau_{\rho_0(t)}(v_0(t)+G_+^{\delta}(v_0(t),c_0(t))+Q_{c_0(t)}),\]
\[\norm{v_0(t)+G_+^{\delta}(v_0(t),c_0(t))}_{H^1}+|c_0(t)-c^*|+|\rho_0(t)-\rho_1(t)| \lesssim \varepsilon_0+\tilde{\varepsilon }^{1/2}\]
and
\[\norm{v_1(t)}_{H^1}+|c_1(t)-c^*| \lesssim \varepsilon _0\]
for $0 \leq t \leq t_0$.
By the conservation of the action $S_c$, we obtain
\begin{align*}
S_{c^*}(u(t))-S_{c^*}(Q_{c^*})=\sum_{j,k}(\Lambda_{k}^{+,j}(t)\Lambda_k^{-,j}(t))+\frac{1}{2}\tbr{\gamma(t),\mathbb{L}_{c^*}\gamma(t)}_{H^1,H^{-1}}+C_N(u(t)),
\end{align*}
where 
\[\Lambda_k^{\pm,j}(t)=(v_1(t),\mathbb{L}_{c^*}F_k^{\mp,j})_{L^2}, \quad \gamma(t)=v_1(t)-\sum_{j,k}(\Lambda_k^{+,j}(t)F_k^{+,j}+\Lambda_k^{-,j}(t)F_k^{-,j})\]
and
\begin{align}\label{eq-sta--1}
C_N(u(t))=S_{c^*}(u(t))-S_{c^*}(Q_{c^*})-\sum_{j,k}(\Lambda_{k}^{+,j}(t)\Lambda_k^{-,j}(t))-\frac{1}{2}\tbr{\gamma(t),\mathbb{L}_{c^*}\gamma(t)}_{H^1,H^{-1}}.
\end{align}
We define $c_u>0$ as
\[\norm{u(t)}_{L^2}=\norm{Q_{c_u}}_{L^2}.\]
Then, we have
\begin{align}\label{eq-sta-0}
|c_u-c^*|\lesssim |\norm{u(0)}_{L^2}-\norm{Q_{c^*}}_{L^2}| \leq \inf_{q \in \R}\norm{u(0)-\tau_qQ_{c^*}}_{H^1}<\tilde{\varepsilon}.
\end{align}
By the inequality \eqref{eq-sta-0} and Lemma \ref{lem-mod-para}, we obtain
\begin{align}\label{eq-sta-0-1}
|c_1(t)-c^*|\leq |c_1(t)-c_u|+|c_u-c^*| \lesssim \tilde{\varepsilon}+\varepsilon_0^2
\end{align}
and
\begin{align}\label{eq-sta-1}
|C_N(u(t))| \lesssim & |c_1(t)-c^*|^2 + \norm{v_1(t)}_{H^1}^3
\lesssim  \tilde{\varepsilon }^2+\varepsilon _0^3.
\end{align}
Since 
\begin{align*}
&\sum_{j,k}|\Lambda_k^{+,j}(t)|^2 \\
\lesssim & \norm{P_+(v_0(t)+G_+^{\delta}(v_0(t),c_0(t)))}_{E}^2+ |\rho_0(t)-\rho_1(t)|^2\norm{v_0(t)+G_+^{\delta}(v_0(t),c_0(t))}_{E}^2\\
= & \norm{G_+^{\delta}(v_0(t),c_0(t))}_{E}^2+|\rho_0(t)-\rho_1(t)|^2\norm{v_0(t)+G_+^{\delta}(v_0(t),c_0(t))}_{E}^2 \\
\lesssim &(l^2+\varepsilon_0^2+\tilde{\varepsilon})\varepsilon_0^2,
\end{align*}
by \eqref{eq-sta--1} and \eqref{eq-sta-1} we obtain
\begin{align}\label{eq-sta-2}
\tbr{\gamma(t),\mathbb{L}_{c^*}\gamma(t)}_{H^1,H^{-1}}\lesssim \tilde{\varepsilon}^2+(l^2+\varepsilon_0^2+\tilde{\varepsilon})^{1/2}\varepsilon_0^2
\end{align}
for $0\leq t \leq t_0$.
Therefore, by the inequalities \eqref{eq-sta-0-1} and \eqref{eq-sta-2} we have
\begin{align}\label{eq-sta-3}
\sum_{j,k}(\Lambda_{k}^{-,j}(t_0))^2 \simeq \norm{v_1(t_0)}_{H^1}^2 \simeq \varepsilon_0^2
\end{align}
for sufficiently small $\tilde{\varepsilon} \ll \varepsilon_0$.
By the system \eqref{LZKeq-1}--\eqref{LZKeq-2}, we have
\begin{align}\label{eq-sta-4}
\frac{d}{dt}\sum_{j,k}(\Lambda_{k}^{-,j}(t))^2\leq 2\sum_{j,k}-\lambda_k(\Lambda_{k}^{-,j}(t))^2+O(\norm{v_1(t)}_{H^1}^3+|c_1(t)-c^*|^3).
\end{align}
Combining \eqref{eq-sta-3} and \eqref{eq-sta-4}, we obtain
\[\frac{d}{dt}\sum_{j,k}(\Lambda_{k}^{-,j}(t))^2<0\]
for $0\leq t\leq t_0$ such that $\sum_{j,k}(\Lambda_{k}^{-,j}(t))^2 \gg \varepsilon_0^3$.
This contradicts \eqref{eq-sta-3} or $\norm{v_1(0)}_{H^1} \lesssim \tilde{\varepsilon} \ll \varepsilon_0$.
Thus, we obtain the conclusion.

\qed

In the following lemmas, to show the property of solutions to the equation \eqref{ZKeq} off the manifold $\mathcal{M}_{cs}^{\delta}(c^*,\varepsilon )$, we prove the estimate of the growth of unstable modes.
\begin{lemma}\label{lem-est-off}
Let $\delta, l_0>0$. 
Suppose
\begin{align}\label{ass-off}
\delta (1+l_0)^{4} \ll \min\{1,l_0^4\}.
\end{align}
There exists $T^*>0$ such that for any solutions $(v_0,c_0,\rho_0)$ and $(v_1,c_1,\rho_1)$ to the system $\eqref{LZKeq-1}$--$\eqref{LZKeq-2}$ satisfying
\begin{align}\label{eq-est-off-1}
 (\mathfrak{m}_{\delta}(\bm{v}_0(0),\bm{v}_1(0))^2-\norm{P_+(v_0(0)-v_1(0))}_{E}^2)^{\frac{1}{2}}
\leq & l_0 \norm{P_+(v_0(0)-v_1(0))}_{E},
\end{align}
one has
\begin{align*}
 (\mathfrak{m}_{\delta}(\bm{v}_0(t),\bm{v}_1(t))^2-\norm{P_+(v_0(t)-v_1(t))}_{E}^2)^{\frac{1}{2}}
\leq & 
\begin{cases}
2l_0\norm{P_+(v_0(t)-v_1(t))}_{E}, & 0\leq t <T^*/2, \\
l_0\norm{P_+(v_0(t)-v_1(t))}_{E}, & T^*/2 \leq t \leq T^*
\end{cases}
\end{align*}
and
\begin{align*}
\norm{P_+(v_0(t)-v_1(t))}_{E}\geq 
\begin{cases}
\frac{1}{2}e^{k_*t/2}\norm{P_+(v_0(0)-v_1(0))}_{E}, & 0\leq t <T^*/2, \\
e^{k_*t/2}\norm{P_+(v_0(0)-v_1(0))}_{E},  & T^*/2 \leq t \leq T^*,
\end{cases}
\end{align*}
where $k_*$ is defined by the equation $\eqref{min-max-eigen}$.
\end{lemma}
\proof
By the assumption \eqref{eq-est-off-1}, we have there exists $C>0$ such that
\[ \norm{P_+(v_0(0)-v_1(0))}_{E} \leq \mathfrak{m}_{\delta}(\bm{v}_0(0),\bm{v}_1(0)) \leq (1+l_0)\norm{P_+(v_0(0)-v_1(0))}_{E}.\]
Lemma \ref{lem-est-md} yields 
\begin{align*}
& \mathfrak{m}_{\delta}(\bm{v}_0(t),\bm{v}_1(t))^2-\norm{P_+(v_0(t)-v_1(t))}_{E}^2 \notag \\
\leq & e^{-2k_*t}\norm{P_-(v_0(0)-v_1(0))}_{E}^2+ \norm{P_0(v_0(0)-v_1(0))}_{E}^2 \notag \\
&+\inf_{q\in \R,j=0,1}\bigl( \norm{P_{\gamma}v_j(0)-\tau_qP_{\gamma}v_{1-j}(0)}_{E}^2+\delta|q|^2\phi_{\delta}(v_{1-j}(0))^2\bigr) \notag \\
& +|\log c_0(0)-\log c_1(0)|^2+C\delta^{\frac{1}{2}}\mathfrak{m}_{\delta}(\bm{v}_0(0),\bm{v}_1(0))^2\notag \\
\leq & (l_0^2+\delta^{\frac{1}{2}}(1+l_0)^{2})\norm{P_+(v_0(0)-v_1(0))}_{E}^2 \notag \\
\leq & (l_0^2+\delta^{\frac{1}{2}}(1+l_0)^{2})(e^{2k_*t}-\delta^{\frac{1}{2}}(1+l_0)^{2})^{-1}\norm{P_+(v_0(t)-v_1(t))}_{E}^2
\end{align*}
for $0\leq t \leq T^*$.
By the assumption \eqref{ass-off}, we obtain the conclusion.
\qed

\begin{lemma}\label{lem-exit}
Let $\delta, l_0>0$. 
Suppose the assumption $\eqref{ass-off}$.
There exists $\varepsilon_*=\varepsilon _*(c^*,\delta,l_0)>0$ such that for any $0<\varepsilon<\varepsilon_*$ and solutions $u_0(t)$ and $u_1(t)$ to the equation $\eqref{ZKeq}$ satisfying
\begin{align}\label{ass-lem-exit}
\sup_{t \geq 0}\inf_{q \in \R} \norm{u_1(t)-\tau_q Q_{c^*}}_{H^1} < \varepsilon, \quad \inf_{q \in \R} \norm{u_0(0)-\tau_q Q_{c^*}}_{H^1} <\varepsilon
\end{align}
and 
\begin{align}\label{lem-exit-p}
(\mathfrak{m}_{\delta}(\bm{v}_0(0),\bm{v}_1(0))^2-\norm{P_+(v_0(0)-v_1(0))}_{E}^2)^{\frac{1}{2}}< l_0 \norm{P_+(v_0(0)-v_1(0))}_{E},
\end{align}
 one has 
 \begin{align}\label{eq-lem-exit-q}
\inf_{q \in \R}\norm{u_0(t_0)-\tau_qQ_{c^*}}_{H^1}\geq \varepsilon 
\end{align}
 for some $t_0>0$, where $(v_0(0),c_0(0),\rho_0(0))$ and $(v_1(0),c_1(0),\rho_1(0))$ satisfy 
\begin{align}\label{cond-exit-orth}
u_j(0)=\tau_{\rho_j(0)}(v_j(0)+Q_{c_j(0)}), \quad|(v_1(0),\partial_xQ_{c^*})_{L^2}|+|(v_1(0),Q_{c^*})_{L^2}|< \varepsilon ^{1/2}
\end{align}
for $j=0,1$. 
\end{lemma}
\proof
Let $v_j$ be the solution to the system \eqref{LZKeq-1}--\eqref{LZKeq-2} with the initial data $(v_j(0),c_j(0),\rho_j(0))$.
we show the inequality \eqref{eq-lem-exit-q} by the contradiction.
Assume for any $0<\varepsilon_* \ll \delta^2$ there exist $0<\varepsilon<\varepsilon_*$ and solutions $u_0(t)$ and $u_1(t)$ to the equation \eqref{ZKeq} satisfying  \eqref{ass-lem-exit}, \eqref{lem-exit-p}, \eqref{cond-exit-orth}
and
\begin{align}\label{eq-lem-exi-cont}
 \sup_{t\geq 0}\inf_{q \in \R}\norm{u_0(t)-\tau_qQ_{c^*}}_{H^1}<\varepsilon.
\end{align}
Since 
\[ \norm{v_1(t)}_{H^1}+|c_1(t)-c^*| \lesssim \inf_{q \in \R }\norm{u_1(t)-\tau_q Q_{c^*}}_{H^1}+\varepsilon^{1/2} \ll \delta \]
and $u_1(t)=\tau_{\rho_1(t)}(v_1(t)+Q_{c_1(t)})$ as long as $\norm{v_1(t)}_{H^1}^2+|c_1(t)-c^*|^2<\delta^2$, we have $u_1(t)=\tau_{\rho_1(t)}(v_1(t)+Q_{c_1(t)})$ for all $t \geq 0$.
Applying Lemma \ref{lem-est-off} repeatedly, we obtain
\begin{align}\label{eq-repulsive-1}
0<\frac{1}{2}e^{k_*t/2}\norm{P_+(v_0(0)-v_1(0))}_{E}\leq \norm{P_+(v_0(t)-v_1(t))}_{E}
\end{align}
and
\begin{align}\label{eq-repulsive-1-1}
(\mathfrak{m}_{\delta}(\bm{v}_0(t),\bm{v}_1(t))^2-\norm{P_+(v_0(t)-v_1(t))}_{E}^2)^{\frac{1}{2}}< 2l_0 \norm{P_+(v_0(t)-v_1(t))}_{E},
\end{align}
for all $t>0$.
Since 
\[ P_+(\tau_qQ_c)=0\]
for $q \in \R$ and $c>0$, we have 
\begin{align*}
\norm{P_+(v_0(t)-v_1(t))}_{E}\lesssim& \norm{P_+v_0(t)}_{H^1}+\norm{P_+v_1(t)}_{H^1} \notag \\
\lesssim& \inf_{q \in\R}\norm{v_0(t)+Q_{c_0(t)}-\tau_qQ_{c^*}}_{H^1}+\inf_{q \in \R}\norm{v_1(t)+Q_{c_1(t)}-\tau_qQ_{c^*}}_{H^1}.
\end{align*}
Therefore, if $0<\varepsilon \ll  \delta$, then by the assumption \eqref{ass-lem-exit} and the inequality \eqref{eq-repulsive-1}  we have 
\begin{align}\label{eq-repulsive-2}
\frac{1}{2}e^{k_*t/2}\norm{P_+(v_0(0)-v_1(0))}_{E} \leq \norm{P_+(v_0(t)-v_1(t))}_{E} \lesssim \inf_{q \in\R}\norm{u_0(t)-\tau_qQ_{c^*}}_{H^1}+\varepsilon 
\end{align}
for $t>0$ as long as $\norm{v_0(t)}_{H^1}^2+|c_0(t)-c_*|^2<\delta^2$.
By Lemma \ref{lem-est-off}, the inequality \eqref{eq-repulsive-1-1} and the assumption \eqref{ass-lem-exit}, we have
\begin{align}
\norm{v_0(t)}_{H^1}^2+|c_0(t)-c_*|^2 
\lesssim & \mathfrak{m}_{\delta}(\bm{v}_0(t),(0,c^*))^2\notag \\
\lesssim & \mathfrak{m}_{\delta}(\bm{v}_0(t),\bm{v}_1(t))^2+\inf_{q \in \R}\norm{(I-P_0)v_1(t)+Q_{c_1(t)}-\tau_qQ_{c^*}}_{H^1}^2 \notag \\
&+|(v_1(t),\partial_xQ_{c^*})_{L^2}|^2+|(v_1(t),Q_{c^*})_{L^2}|^2 \notag \\
\lesssim & (1+l_0)\norm{P_+(v_0(t)-v_1(t))}_{E}^2+ \varepsilon \label{eq-repulsive-3}
\end{align}
The inequalities \eqref{eq-repulsive-2}, \eqref{eq-repulsive-3} and \eqref{eq-repulsive-2} contradict the assumption \eqref{eq-lem-exi-cont} for sufficiently small $\varepsilon_*>0$.
Thus, the proof was completed.
\qed

In the following corollary, we show that solutions to the equation \eqref{ZKeq} off the center--stable manifold exit neighborhoods of a  line solitary wave.
\begin{corollary}\label{cor-repulsive}
Let $\delta,l>0$. Suppose $\eqref{ass-smallness}$ and $\eqref{ass-off}$.
There exists $\varepsilon ^*=\varepsilon ^*(c^*,\delta,l_0)>0$ such that for $u(0) \in N_{c^*}(\varepsilon ^*) \setminus \mathcal{M}_{cs}^{\delta}(c^*,\varepsilon ^*)$, the solution $u$ of the equation $\eqref{ZKeq}$ corresponding to the initial data $u(0)$ satisfies 
\[\inf_{q \in \R}\norm{u(t_0)-\tau_qQ_{c^*}}_{H^1}\geq \varepsilon ^*,\]
for some $t_0 \geq 0$.
\end{corollary}
\proof
Let $u(0) \in N_{c^*}(\varepsilon ^*) \setminus \mathcal{M}_{cs}^{\delta}(c^*,\varepsilon ^*)$ and $u$ be the solution to the equation \eqref{ZKeq} corresponding to the initial data $u(0)$.
By applying Lemma \ref{lem-orth}, we define $v(0)=\tau_{-\rho(u(0))}u(0)-Q_{c(u(0))}$ and the solution $(v_1,c_1,\rho_1)$ to the system \eqref{LZKeq-1}--\eqref{LZKeq-2} corresponding to the initial data $(( P_-+P_{\gamma})v(0)+G_+(( P_-+P_{\gamma})v(0), c(u(0))), c(u(0)),\rho(u(0)))$.
Then, the solution $(v_0,c_0,\rho_0)$ to the system \eqref{LZKeq-1}--\eqref{LZKeq-2} corresponding to the initial data $(v(0),c(u(0)),\rho(u(0)))$ satisfies $u(t)=\tau_{\rho_0(t)}(v_0(t)+Q_{c_0(t)})$ as long as $\norm{v_0(t)}_{H^1}^2+|c_0(t)-c^*|^2<\delta^2$.
Since
\[(P_-+P_{\gamma})(v_0(0)-v_1(0))=0,\]
we have
\[\mathfrak{m}_{\delta}(\bm{v}_0(0),\bm{v}_1(0))^2-\norm{P_+(v_0(0)-v_1(0))}_{E}^2 =0\leq l_0^2 \norm{P_+(v_0(0)-v_1(0))}_{E}^2.\]
Thus,  the conclusion follows Lemma \ref{lem-exit} and the inequality \eqref{eq-repulsive-2}.
\qed

In the following corollary, we show the correspondence between $\mathcal{M}_{cs}^{\delta}(c^*,\varepsilon) $ and $\tilde{\mathcal{M}}_{cs}^{\delta}(c^*,\varepsilon)$.
\begin{corollary}\label{cor-corr}
Let $\delta,l>0$. Suppose $\eqref{ass-smallness}$.
There exist $\varepsilon_1=\varepsilon_1(c^*,\delta)>0$ such that for $0<\varepsilon<\varepsilon_1$ 
\begin{align*}
\mathcal{M}_{cs}^{\delta}(c^*,\varepsilon)=\tilde{\mathcal{M}}_{cs}^{\delta}(c^*,\varepsilon)
\end{align*}
and
\begin{align*}
&\{w+G_+^{\delta}(w,c)+Q_{c} ;w \in P_{\leq 0}H^1(\RTL), |c-c^*|\leq c^*/2, \norm{P_0(w+Q_c-Q_{c^*})}_{H^1}<\varepsilon^{1/2}, \\
& \qquad \inf_{q \in \R} \norm{w+G_+^{\delta}(w,c)+Q_c-\tau_qQ_{c^*}}_{H^1} <\varepsilon\}\\
=&\{ \tau_{\rho}(w+G_+^{\delta}(w,c)+Q_c); w \in (P_-+P_{\gamma})H^1(\RTL), |c-c^*|\leq c^*/2, \rho \in \R, \\
&  \norm{P_0(\tau_\rho(w+G_+^{\delta}(w,c)+Q_c)-Q_{c^*})}_{H^1}<\varepsilon^{1/2} , \inf_{q \in \R} \norm{w+G_+^{\delta}(w,c)+Q_c-\tau_qQ_{c^*}}_{H^1} <\varepsilon\}.
\end{align*}
Moreover, 
\begin{align}\label{eq-corr}
\mathcal{M}_{cs}^{\delta}&(c^*,\varepsilon)=\bigcup_{q \in \R} \{w+G_+^{\delta}(w,c)+Q_{c} ;w \in P_{\leq 0}H^1(\RTL), |c-c^*|\leq c^*/2, \notag \\
&\norm{P_0(w+Q_c-Q_{c^*})}_{H^1}<\varepsilon^{1/2}, 
\inf_{q \in \R} \norm{w+G_+^{\delta}(w,c)+Q_c-\tau_qQ_{c^*}}_{H^1} <\varepsilon\}.
\end{align}
\end{corollary}
\proof
By the definitions of $\mathcal{M}_{cs}^{\delta}(c^*,\varepsilon)$ and $\tilde{\mathcal{M}}_{cs}^{\delta}(c^*,\varepsilon)$, we have 
\[\mathcal{M}_{cs}^{\delta}(c^*,\varepsilon) \subset \tilde{\mathcal{M}}_{cs}^{\delta}(c^*,\varepsilon).\]
From Theorem \ref{thm-ex-csm}, solutions $u(t)$ to the equation \eqref{ZKeq} with an initial data $u(0) \in \tilde{\mathcal{M}}_{cs}^{\delta}(c^*,\varepsilon)$ satisfy
\[\sup_{t\geq 0}\inf_{q \in \R} \norm{u(t)-\tau_qQ_{c^*}}_{H^1}<\varepsilon_*\]
 for sufficiently small $\varepsilon>0$, where $ \varepsilon_*$ is defined in Corollary \ref{cor-repulsive}.
 Therefore, Corollary \ref{cor-repulsive} yields $\tilde{\mathcal{M}}_{cs}^{\delta}(c^*,\varepsilon) \cap (H^1(\RTL) \setminus \mathcal{M}_{cs}^{\delta}(c^*,\varepsilon))=\emptyset$.
Thus, we obtain $\mathcal{M}_{cs}^{\delta}(c^*,\varepsilon)=\tilde{\mathcal{M}}_{cs}^{\delta}(c^*,\varepsilon)$.
Since $\mathcal{M}_{cs}^{\delta}(c^*,\varepsilon)=\tilde{\mathcal{M}}_{cs}^{\delta}(c^*,\varepsilon)$, we have
\begin{align*}
&\{w+G_+^{\delta}(w,c)+Q_{c} ;w \in P_{\leq 0} H^1(\RTL), |c-c^*|\leq c^*/2, \norm{P_0(w+Q_c-Q_{c^*})}_{H^1}<\varepsilon^{1/2}, \\
& \qquad \inf_{q \in \R} \norm{w+G_+^{\delta}(w,c)+Q_c-\tau_qQ_{c^*}}_{H^1} <\varepsilon\}\\
\subset &\{ \tau_{\rho}(w+G_+^{\delta}(w,c)+Q_c); w \in (P_-+P_{\gamma})H^1(\RTL), |c-c^*|\leq c^*/2, \rho \in \R, \\
&  \norm{P_0(\tau_\rho(w+G_+^{\delta}(w,c)+Q_c)-Q_{c^*})}_{H^1}<\varepsilon^{1/2} , \inf_{q \in \R} \norm{w+G_+^{\delta}(w,c)+Q_c-\tau_qQ_{c^*}}_{H^1} <\varepsilon\}.
\end{align*}
 Let $\tau_{\rho}(w+G_+^{\delta}(w,c)+Q_c)$ satisfy that $w \in (P_-+P_{\gamma})H^1(\RTL)$, 
\[\norm{P_0(\tau_\rho(w+G_+^{\delta}(w,c)+Q_c)-Q_{c^*})}_{H^1}<\varepsilon^{1/2} \]
 and 
\[ \inf_{q \in \R} \norm{w+G_+^{\delta}(w,c)+Q_c-\tau_qQ_{c^*}}_{H^1} <\varepsilon.\]
We define the solution $u_0(t)$ to the equation \eqref{ZKeq} with the initial data $\tau_{\rho}(w+G_+^{\delta}(w,c)+Q_c)$ and the solution $u_1(t)$ to the equation \eqref{ZKeq} with the initial data 
\[P_{\leq 0} w_0+G_+(P_{\leq 0} w_0,c^*)+Q_{c^*},\]
where $w_0=\tau_{\rho}(w+G_+^{\delta}(w,c)+Q_c)-Q_{c^*}$.
By Theorem \ref{thm-ex-csm}, there exists $l_0$ such that $l_0$ and $\delta$ satisfies $\eqref{ass-off}$ and
\[\sup_{t\geq 0}\inf_{q\in \R}\norm{u_j(t)-\tau_qQ_{c^*}}_{H^1}<\varepsilon_*(c^*,\delta ,l_0)\]
for $j=0,1$ and sufficiently small $\varepsilon>0$, where $\varepsilon_*$ is defined in Lemma \ref{lem-exit}.
Since $u_1(0)$ satisfy the assumption \eqref{cond-exit-orth} for $\varepsilon<\varepsilon_*$, by Lemma \ref{lem-exit} we have 
\begin{align*}
l_0\norm{P_+(u_0(0)-u_1(0))}_{E}\lesssim & \norm{P_{\leq 0} (u_0(0)-u_1(0))}_{H^1}=0
\end{align*}
and $u_0(0)=u_1(0)$.
Therefore, we obtain 
 \begin{align*}
&\{w+G_+^{\delta}(w,c)+Q_{c} ;w \in P_{\leq 0}H^1(\RTL), |c-c^*|\leq c^*/2, \norm{P_0(w+Q_c-Q_{c^*})}_{H^1}<\varepsilon^{1/2}, \\
& \qquad \inf_{q \in \R} \norm{w+G_+^{\delta}(w,c)+Q_c-\tau_qQ_{c^*}}_{H^1} <\varepsilon\}\\
\supset &\{ \tau_{\rho}(w+G_+^{\delta}(w,c)+Q_c); w \in (P_-+P_{\gamma})H^1(\RTL), |c-c^*|\leq c^*/2, q\in \R, \\
&  \norm{P_0(\tau_\rho(w+G_+^{\delta}(w,c)+Q_c)-Q_{c^*})}_{H^1}<\varepsilon^{1/2} , \inf_{q \in \R} \norm{w+G_+^{\delta}(w,c)+Q_c-\tau_qQ_{c^*}}_{H^1} <\varepsilon\}.
\end{align*}
Since
\[\norm{P_0(w+Q_c-Q_{c^*})}_{H^1} \lesssim \inf_{q \in \R} \norm{w+G_+^{\delta}(w,c)+Q_c-\tau_qQ_{c^*}}_{H^1}\]
for $w \in (P_-+P_{\gamma})H^1(\RTL)$ and $|c-c^*|<c^*/2$, we have the equation \eqref{eq-corr} for sufficiently small $\varepsilon >0$.

\qed

\section{Smoothness of the center stable manifolds}
In this section, we show the center stable manifolds has the $C^1$ regularity by applying the argument in \cite{K N S}.

The following lemma shows the local uniqueness of $\mathcal{M}_{cs}^{\delta}$.
\begin{lemma}\label{lem-uniqueness-M}
Let $l,l_0, \delta_0,\delta_1>0$. Assume $(\delta_0,l)$ and $(\delta_1,l)$ satisfy $\eqref{ass-smallness}$ and assume $(\delta_0,l_0)$ and $(\delta_1,l_0)$ satisfy $\eqref{ass-off}$.
Then, there exists $r_0=r_0(\delta_0,\delta_1)>0$ such that 
\[G_+^{\delta_0}(w,c)=G_+^{\delta_1}(w,c),\]
where $(w,c) \in P_{\leq 0}H^1(\RTL)\times (0,\infty)$ satisfying 
\[\inf_{q\in \R, j=0,1}\norm{w+G_+^{\delta_j}(w,c)+Q_c-\tau_qQ_{c^*}}_{H^1}<r_0.\]
Moreover, $\mathcal{M}_{cs}^{\delta_0}(c^*,r_0)=\mathcal{M}_{cs}^{\delta_1}(c^*,r_0)$. 
\end{lemma}
\proof
By Theorem \ref{thm-ex-csm}, there exists $r_1>0$ such that for $t \geq 0$, $j \in \{0,1\}$ and $\phi_j \in \mathcal{M}_{cs}^{\delta_j}(c^*,r_1)$,
\[v_j(t) \in \mathcal{M}_{cs}^{\delta_j}(c^*,\varepsilon_0),\]
where $\varepsilon _0=\min_{j =0,1} \varepsilon _*(c^*,\delta_j,l_0)/2$, $v_j$ is the solution to the equation \eqref{ZKeq} with $v_j(0)=\phi_j$ and the constant $\varepsilon _*$ is defined in Lemma \ref{lem-exit}.
We prove the conclusion by contradiction.
Assume for any $r>0$ there exist $j_0 \in \{0,1\}$ and 
\[\tau_{\rho}(w+G_+^{\delta_{j_0}}(w,c)+Q_c) \in \mathcal{M}_{cs}^{\delta_{j_0}}(c^*,r) \setminus \mathcal{M}_{cs}^{\delta_{1-j_0}}( c^*,r)\]
such that $w \in (P_-+P_{\gamma}) H^1(\RTL)$ and 
\[\inf_{q \in \R} \norm{w+G_+^{\delta_{j_0}}(w,c)+Q_c-\tau_qQ_{c^*}}_{H^1}<r.\]
Then, $G_+^{\delta_{j_0}}(w,c)\neq G_+^{\delta_{1-j_0}}(w,c)$ and
\begin{align*}\label{eq-uni-M-1}
\inf_{q \in \R} \norm{w+G_+^{\delta_{1-j_0}}(w,c)+Q_c-Q_{c^*}}_{H^1} \leq r+ \norm{G_+^{\delta_{j_0}}(w,c)-G_+^{\delta_{1-j_0}}(w,c)}_{H^1} \lesssim r.
\end{align*}
Without loss of generality, we can choose $r$ satisfying 
\[\max_{j=0,1} \inf_{q \in \R} \norm{w+G_+^{\delta_j}(w,c)+Q_c-\tau_qQ_{c^*}}_{H^1}<\min\{r_1,\varepsilon _0\}\]
and
\[\norm{w}_{H^1}^2+|c-c^*|^2<\min_{j=0,1}\delta_j.\]
Let $u_j$ be the solution to \eqref{ZKeq} with $u_j(0)=\tau_{\rho}(w+G_+^{\delta_j}(w,c)+Q_c)$.
Then, we have for $t \geq 0$ and $j \in \{0,1\}$
\begin{align}\label{eq-uni-M-3}
u_j (t) \in \mathcal{M}_{cs}^{\delta_j}(c^*, \varepsilon _0).
\end{align}
By the inequality 
\begin{align*}
&\bigl(\mathfrak{m}_{\delta_j}((w+G_+^{\delta_0}(w,c),c),(w+G_+^{\delta_1}(w,c),c))^2 - \norm{G_+^{\delta_0}(w,c)-G_+^{\delta_1}(w,c)}_{E}^2\bigr)^{1/2} \\
<& l_0 \norm{G_+^{\delta_0}(w,c)-G_+^{\delta_1}(w,c)}_{E}
\end{align*}
and Lemma \ref{lem-exit}, there exists $t_0>0$ such that 
\[\inf_{q \in \R}\norm{u_{1-j}(t_0)-\tau_qQ_{c^*}}_{H^1}=\varepsilon_*(c^*,\delta_j,l_0)\]
which contradict \eqref{eq-uni-M-3}.
Thus, there exists $r>0$ such that
\[ \mathcal{M}_{cs}^{\delta_0}(c^*,r)=\mathcal{M}_{cs}^{\delta_1}(c^*,r).\]
\qed

The following corollary shows the tangent plain of the center stable manifolds $\mathcal{M}_{cs}$ at $\tau_q Q_{c^*}$ is $\tau_q (P_{\leq 0}H^1(\RTL)+Q_{c^*})$. 
\begin{corollary}\label{cor-Lip}
Let $\delta, l>0$.
Assume $\eqref{ass-smallness} $.
Then, for any $l_*>0$ there exists $\delta_{l_*}>0$ such that 
\[\norm{G_+^{\delta}(w_0,c_0)-G_+^{\delta}(w_1,c_1)}_{E}\leq l_* (\norm{w_0-w_1}_{E}+|\log c_0 -\log c_1|)\]
for $w_0,w_1 \in P_{\leq 0}H^1(\RTL)$ and $c_0, c_1>0$ satisfying
\[\max_{j=0,1}(\norm{w_j}_{E}+|\log c_j -\log c^{*}|)\leq \delta_{l_*}.\]
\end{corollary}
\proof
For any $l_*>0$, there exists $\delta_0>0$ such that $G_+^{\delta_0} \in \mathscr{G}_{l_*,\delta_0}^{+}$.
By Lemma \ref{lem-uniqueness-M}, there exists $r>0$ such that $G_+^{\delta}(w,c)=G_+^{\delta_0}(w,c)$ for $(w,c) \in P_{\leq 0}H^1(\RTL)\times (0,\infty)$ satisfying 
\[\inf_{q\in \R}\norm{w+G_+^{\delta}(w,c)+Q_c-\tau_qQ_{c^*}}_{H^1}<r.\]
Since there exists $C>0$ such that
\[\inf_{q\in \R}\norm{w+G_+^{\delta}(w,c)+Q_c-\tau_qQ_{c^*}}_{H^1} \leq C (\norm{w}_{H^1}+|\log c -\log c^*|), \]
we have
\begin{align*}
\norm{G_+^{\delta}(w_0,c_0)-G_+^{\delta}(w_1,c_1)}_{E}=&\norm{G_+^{\delta_0}(w_0,c_0)-G_+^{\delta_0}(w_1,c_1)}_{E}\\
\leq & l_* (\norm{w_0-w_1}_{E}+|\log c_0 -\log c_1|)
\end{align*}
for $(w_0,c_0),(w_1,c_1) \in P_{\leq 0}H^1(\RTL)\times (0,\infty)$ satisfying 
\[\max_{j=0,1}(\norm{w_j}_{H^1}+|\log c_j -\log c^*|)< C^{-1}r.\]
\qed

In the rest of this section, we prove that $G_+^{\delta}$ is at least $C^1$ in $P_{\leq 0} H^1(\RTL) \times (0,\infty)$ by applying the argument in the section 2.3 in \cite{K N S}.
Let $\varepsilon , a_0,a_1>0$ and $\psi_0, \psi_1 \in (P_-+P_{\gamma}+P_1)H^1(\RTL)$ with $\norm{\psi_0}_{H^1(\RTL)}<\varepsilon $ and $|a_0-c^*|<\varepsilon $. 
We consider solution $(v_0,c_0,\rho_0)$ to the system \eqref{vZKeq} and \eqref{eq-orth-1} such that 
\[v_0(0)=\psi_0+G_+^{\delta}(\psi_0,a_0), \quad c_0(0)=a_0.\]
Let $v_h$ be a solution to the equation
\begin{align}\label{eq-v0}
v_t=\partial_x \mathbb{L}_{c^*}v+(\dot{\rho}_0-c^*)\partial_xv+2\partial_x((Q_{c^*}-Q_{c_h})v)+(\dot{\rho}_0-c_h)\partial_xQ_{c_h}-\dot{c}_0\partial_cQ_{c_h}-\partial_x(v^2)
\end{align} 
with the initial data $v_h(0)=\psi_0+h\psi_1+G_+^{\delta}(\psi_0+h\psi_1,a_0+ha_1)$, where 
\[c_h(t)=c_0(t)+ha_1.\]
Then, $\tau_{\rho_0}(v_0+Q_{c_0})$ and $\tau_{\rho_0}(v_h+Q_{c_h})$ are solutions to the equation \eqref{ZKeq}.
By the Lipschitz continuity of $G_+^{\delta}$, for any sequence $\{h_n\}_n$ with $h_n \to 0$  as $n \to \infty$ there exist a subsequence $\{h_n'\}_n \subset \{h_n\}_n$ and $\psi_+ \in P_+H^1(\RTL)$ such that 
\[\frac{G_+^{\delta}(\psi_0+h'_n\psi_1,a_0+h'_na_1)-G_+^{\delta}(\psi_0,a_0)}{h'_n} \to \psi_+  \mbox{ as } n \to \infty.\] 
Let $w_h=\tau_{\rho_0}v_h$ for $h \geq 0$.
Then, for $h\geq 0$, $w_h$ is the solution to the equation
\[w_t=-\partial_x\Delta w - 2\partial_x(\tau_{\rho_0}Q_{c_h}w)+(\dot{\rho}_0-c_h)\tau_{\rho_0}\partial_xQ_{c_h}-\dot{c}_0\tau_{\rho_0}\partial_cQ_{c_h}-\partial_x(w^2)\]
 with $w_h(0)=v_h(0)$.
From the well-posedness result of the equation \eqref{ZKeq} in \cite{M P}, we have there exists $b_0 > \frac{1}{2}$ such that  for $T>0$ and $\frac{1}{2}<b<b_0$ there exists $C=C(T,b)>0$ satisfying
\begin{align}\label{eq-w0}
\norm{w_0}_{X^{1,b}_T} \leq C \norm{v_0(0)}_{H^1}.
\end{align}
We define $\xi$ as the solution to the equation
\begin{align}\label{eq-xi-eta}
\xi_t=&-\partial_x\Delta\xi-2\partial_x(\tau_{\rho_0}Q_{c_0}\xi)-2a_1\partial_x(\tau_{\rho_0}\partial_cQ_{c_0}w_0)-a_1\tau_{\rho_0}\partial_xQ_{c_0}\notag \\
&+a_1(\dot{\rho}_0-c_0)\tau_{\rho_0}\partial_x\partial_cQ_{c_0}-\dot{c}_0a_1\tau_{\rho_0}\partial_c^2Q_{c_0}-2\partial_x(w_0\xi)
\end{align}
with the initial data $\xi(0)=\psi_1+\psi_+$.
By the smoothness of the flow map of the equation \eqref{ZKeq} given by \cite{M P}, we have that for $T>0$
\begin{align}\label{eq-c1-0}
&\norm{\frac{w_{h'_n}-w_0}{h'_n}-\xi}_{L^{\infty}((-T,T),H^1)} \to 0 \mbox{ as } n \to \infty.
\end{align}
Let $\eta=\tau_{-\rho_{0}}\xi$.
Then, $\eta$ satisfies the equation
\begin{align}\label{eq-eta}
\eta_t=&\partial_x\mathbb{L}_{c^*}\eta-2\partial_x((Q_{c_0}-Q_{c^*})\eta)+(\dot{\rho}_0-c^*)\partial_x\eta -2a_1\partial_x(\partial_cQ_{c_0}v_0)\notag \\
&-a_1\partial_xQ_{c_0}+a_1(\dot{\rho}_0-c_0)\partial_x\partial_cQ_{c_0}-\dot{c}_0a_1\partial_c^2Q_{c_0}-2\partial_x(v_0\eta).
\end{align}
We define the norm $\norm{\cdot}_{E_{\kappa}}$ by 
\[ \norm{u}_{E_{\kappa}}=\norm{(I-P_1)u+\kappa P_1u}_{E_{\kappa}}.\]
In the following lemma, we show the behavior of solutions of the equation \eqref{eq-eta} (see Lemma 2.4 in \cite{K N S}.)
\begin{lemma}\label{lem-prop-eta}
Let $C,\kappa,K_0>0$. 
There exists $K_1(C,K_0),\kappa_0(C,K_0)>0$ satisfying the following property.
Let $(v_0,c_0,\rho_0)$ be a solution to the system $\eqref{vZKeq}$ and $\eqref{eq-orth-1}$ satisfying $\sup_{t\geq 0}(\norm{v_0(t)}_{H^1}+|c_0(t)-c^*|)\leq C\kappa$ and $\tau_{\rho_0}(v_0+Q_{c_0})$ is a solution to the equation $\eqref{ZKeq}$.
Then, for $a_1 \in \R$ the Cauchy problem of the equation $\eqref{eq-eta}$ is global well-posed in $H^1(\RTL)$.
Precisely, there exists $b> \frac{1}{2}$ such that for any  $\eta_0 \in H^1(\RTL)$ there exists a unique solution $\xi$ to the equation $\eqref{eq-xi-eta}$ satisfying that $\xi(0)=\tau_{-\rho_0(0)}\eta_0$, 
\[w \in X^{1,b}_{T} \mbox{ for } T>0\]
and $\tau_{-\rho_0}\xi$ is a solution to the equation $\eqref{eq-eta}$ with initial data $\eta_0$.
Moreover, if a solution $\eta$ to the equation $\eqref{eq-eta}$ with initial data $\eta(0) \in H^1(\RTL)$ satisfies $0<\kappa<\kappa_0$ and  
\begin{align}\label{cond-eta}
K_0\kappa^{\frac{1}{3}}(\norm{P_{\leq 0}\eta(t_0)}_{E_{\kappa^{1/3}}}+|a_1|)< \norm{P_+\eta(t_0)}_{E}
\end{align}
at some $t_0\geq 0$, then for $t \geq t_{0}+1/2$ 
\begin{align}\label{ineq-eta-1}
3 \norm{P_+\eta(t)}_{E} > e^{\frac{k_*}{2}(t-t_0)}(\norm{P_+\eta(t_0)}_{E}+K_0\kappa^{1/3}(\norm{P_{\leq 0}\eta(t)}_{E_{\kappa^{1/3}}}+|a_1|)).
\end{align}
On the other hand, if $\eqref{cond-eta}$ fails for $t_0\geq 0$, then for $t\geq 0$
\begin{align}\label{ineq-eta-2}
\norm{P_+\eta(t)}_{E}\lesssim \kappa^{\frac{1}{3}}(\norm{P_{\leq 0}\eta(t)}_{E_{\kappa^{1/3}}}+|a_1|)\lesssim e^{K_1\kappa^{1/6}t}\kappa^{\frac{1}{3}}(\norm{P_{\leq 0}\eta(0)}_{E_{\kappa^{1/3}}}+|a_1|).
\end{align}
\end{lemma}
\proof
From the same manner of the proof of Theorem \ref{thm-gwp-LZKeq}, we obtain the global well-posedness of the equation \eqref{eq-eta} in $H^1(\RTL)$ and for  $a_1 \in \R$.
Moreover, for any solutions $\eta$ to the equation \eqref{eq-eta} and $s\geq 0$ there exists $\xi_s \in X^{1,b}$ such that for $t \in (s-1,s+1)$ we have $\eta (t) = \tau_{\rho_0(s)-\rho_0(t)}\xi_s(t)$ and
\begin{align}\label{eq-lem-eta-1}
\norm{\eta}_{L^{\infty}((s-1,s+1),H^1)}\lesssim \norm{\xi_s}_{X^{1,b}}\lesssim (\norm{\eta(s)}_{H^1}+|a_1|).
\end{align}
By the inequalities \eqref{eq-gwp-a-3}, \eqref{eq-w0} and \eqref{eq-lem-eta-1}, for $t_1,t_2 \geq 0$ with $|t_1-t_2|<1$ we have
\begin{align}\label{eq-lem-eta-2}
|\norm{P_{\gamma}\eta(t_2)}_{E}^2-\norm{P_{\gamma}\eta (t_1)}_{E}^2| \lesssim& \kappa (\norm{\eta(t_1)}_{E}+|a_1|)^2\notag \\
 \lesssim& \kappa^{\frac{1}{3}}\norm{P_{\leq 0}\eta(t_1)}_{E_{\kappa^{1/3}}}^2+\kappa |a_1|^2+\kappa \norm{P_+\eta(t_1)}_{E}^2
\end{align}
Since
\[ \norm{P_-\partial_t\eta(t)-\partial_x\mathbb{L}_{c^*}P_-\eta (t)}_{E}+|(Q_{c^*},\partial_t\eta(t))_{L^2}|\lesssim \kappa (\norm{\eta(t)}_{E}+|a_1|),\]
we have
\begin{align}
|\norm{P_2\eta(t_2)}_{E_{\kappa^{1/3}}}-\norm{P_2\eta(t_1)}_{E_{\kappa^{1/3}}}| \lesssim \kappa^{2/3}\norm{\eta(t_1)}_{E_{\kappa^{1/3}}}+\kappa |a_1| \label{eq-lem-eta-3}
\end{align}
and
\begin{align}
\norm{P_-\eta(t_2)}_{E} -e^{-k_*(t_2-t_1)}\norm{P_-\eta(t_1)}_{E} \lesssim \kappa^{2/3}\norm{\eta(t_1)}_{E_{\kappa^{1/3}}}+\kappa |a_1| \label{eq-lem-eta-3-1}
\end{align}
 for $t_1,t_2 \in \R$ with $|t_1-t_2|<1$.
By the inequality 
\[|(\partial_xQ_{c^*},\partial_t\eta(t))_{L^2}|\lesssim \norm{P_{\gamma}\eta(t)}_{E}+\norm{P_2\eta(t)}_{E}+|a_1|+\kappa \norm{P_{d}\eta(t)}_{E}\]
and  \eqref{eq-lem-eta-1}--\eqref{eq-lem-eta-3} we obtain for $t_1,t_2 \in \R$ with $|t_1-t_2|<1$
\begin{align}
&|\norm{P_1\eta(t_2)}_{E_{\kappa^{1/3}}}-\norm{P_1\eta(t_1)}_{E_{\kappa^{1/3}}}|\notag \\
\lesssim & \kappa^{1/3} (\norm{(P_{\gamma}+P_0)\eta(t_1)}_{E_{\kappa^{1/3}}} +|a_1|)+\kappa^{1/2}\norm{(P_d-P_1)\eta(t_1)}_{E_{\kappa^{1/3}}}, \label{eq-lem-eta-4}
\end{align}
for small $\kappa$.
From the inequality \eqref{eq-lem-eta-2}--\eqref{eq-lem-eta-4}, there exists $C>0$ such that
\begin{align}\label{eq-lem-eta-5}
\norm{P_{\leq 0}\eta(t_2)}_{E_{\kappa^{1/3}}}\leq (1+C\kappa^{1/6})\norm{P_{\leq 0}\eta(t_1)}_{E_{\kappa^{1/3}}}+C\kappa^{1/3}|a_1|+C\kappa^{1/2} \norm{P_+\eta(t_1)}_{E}
\end{align}
for $t_1,t_2 \geq 0$ with $|t_1-t_2|<1$.
The inequality 
\begin{align}\label{eq-lem-eta-6}
\norm{P_+\partial_t\eta(t)-\partial_x\mathbb{L}_{c^*}P_+\eta(t)}_{E} \lesssim \kappa (\norm{\eta(t)}_{E}+|a_1|)
\end{align}
implies that there exists $C>0$ such that for $t_1,t_2 \geq 0$ with $|t_1-t_2|<1$ 
\begin{align}\label{eq-lem-eta-7}
\partial_t\norm{P_+\eta(t_2)}_{E} \geq k_*\norm{P_+\eta(t_2)}_{E} -C\kappa^{2/3}(\norm{\eta(t_1)}_{E_{\kappa^{1/3}}}+|a_1|).
\end{align}

Suppose \eqref{cond-eta} for some $t_0$.
By the assumption \eqref{cond-eta} and the inequality \eqref{eq-lem-eta-7}, we have
\begin{align}\label{eq-lem-eta-8}
\norm{P_+\eta(t)}_{E} \geq& e^{k_*(t-t_0)}\norm{P_+\eta(t_0)}_{E} -(e^{k_*(t-t_0)}-1) C \kappa^{1/3}\norm{P_+\eta(t_0)}_{E}\notag \\
 \geq& (1-C\kappa^{1/3})e^{k_*(t-t_0)}\norm{P_+\eta(t_0)}_{E}
\end{align}
for $t_0\leq t <t_0+1$.
From the assumption \eqref{cond-eta} and the inequalities \eqref{eq-lem-eta-5} and \eqref{eq-lem-eta-8}, we obtain
\begin{align}
\norm{P_+\eta(t)}_E
> & (1-2C\kappa^{1/3})(1+C\kappa^{1/6})^{-1}e^{k_*(t-t_0)}K_0\kappa^{1/3}(\norm{P_{\leq 0}\eta(t)}_{E_{\kappa^{1/3}}}+|a_1|)\label{eq-lem-eta-8-1}
\end{align}
for $t_0\leq t<t_0+1$ and small $\kappa>0$.
Thus, we have
\begin{align*}
\norm{P_+\eta(t)}_{E}> K_0\kappa^{1/3}(\norm{P_{\leq 0}\eta(t)}_{E_{\kappa^{1/3}}}+|a_1|)
\end{align*}
for $t_0+1/2\leq t<t_0+1$ and small $\kappa>0$.
Applying this manner repeatedly, by the inequality \eqref{eq-lem-eta-8} and \eqref{eq-lem-eta-8-1} we obtain the inequality \eqref{ineq-eta-1} for $t>t_0+1/2$.

Suppose \eqref{cond-eta} fails for $t\geq 0$.
Then, the inequality \eqref{eq-lem-eta-5} yields the inequality \eqref{ineq-eta-2} for all $t \geq 0$ and some $K_1>0$.
\qed

Next we prove the uniqueness of solutions to the equation \eqref{eq-eta} not satisfying \eqref{ineq-eta-2} for some $t\geq 0$.
\begin{lemma}\label{lem-eta-exis}
Let $ K_0>0$. Then, there exists $\kappa_1>0 $ such that for $0<\kappa<\kappa_1$, $a_1 \in \R$ and  a solution $(v_0,c_0,\rho_0)$ to the system $\eqref{vZKeq}$ and $\eqref{eq-orth-1}$ with $\sup_{t\geq 0}(\norm{v_0(t)}_{H^1}+|c_0(t)-c^*| ) \leq \kappa$ and for the solutions $\eta_1$ and $\eta_2$ to the equation $\eqref{eq-eta}$ with $P_{\leq 0}\eta_1(0)=P_{\leq 0} \eta_2(0)$ not satisfying that $\eqref{ineq-eta-2}$ for some $t\geq 0$, we have $P_+\eta_1(0)=P_+\eta_2(0)$.
\end{lemma}
\proof
Assume there exist $0<\kappa\ll \kappa_0(1,K_0)$, $a_1 \in \R$, a solution $(v_0,c_0,\rho_0)$ to the system $\eqref{vZKeq}$ and $\eqref{eq-orth-1}$ with $\sup_{t\geq 0}(\norm{v_0(t)}_{H^1}+|c_0(t)-c^*| ) \leq \kappa$ and solutions $\eta_1,\eta_2$ such that $P_{\leq 0}\eta_1(0)=P_{\leq 0} \eta_2(0)$, $P_+\eta_1(0)\neq P_+\eta_2(0)$ and $\eta_1$ and $\eta_2$ do not satisfy that \eqref{ineq-eta-2} for some $t\geq 0$.
Then, $\eta=\eta_1-\eta_2$ is the solution to the equation \eqref{eq-eta} with $a_1=0$ 
\begin{align*}
\eta_t=&\partial_x\mathbb{L}_{c^*}\eta-2\partial_x((Q_{c_0}-Q_{c^*})\eta)+(\dot{\rho}_0-c^*)\partial_x\eta -2\partial_x(v_0\eta).
\end{align*}
Since
\[\kappa^{1/3}K_0\norm{P_{\leq 0}\eta(0)}_{E_{\kappa^{1/3}}} < \norm{P_+\eta(0)}_{E},\]
by Lemma \ref{lem-prop-eta} we have for $t \geq 1/2$ 
\begin{align}\label{eq-lem-exis-1}
3 \norm{P_+\eta(t)}_{E} \geq e^{\frac{k_*}{2}t}(\norm{P_+\eta(0)}_{E} +\kappa^{1/3}\norm{P_{\leq 0}\eta(t)}_{E_{\kappa^{1/3}}}).
\end{align}
On the other hand, by \eqref{ineq-eta-2} we have for $t\geq 0$ 
\[\norm{P_+\eta(t)}_{E} \lesssim e^{K_1\kappa^{1/6}t}(\norm{P_{\leq 0}\eta_1(0)}_{E_{\kappa^{1/3}}}+\norm{P_{\leq 0}\eta_2(0)}_{E_{\kappa^{1/3}}}+|a_1|)\]
and $K_1\kappa^{1/6}\ll k_*$, where $K_1=K_1(1,K_0)$ is defined in Lemma \ref{lem-prop-eta}.
This contradicts the inequality \eqref{eq-lem-exis-1}.
Thus, $P_+\eta_1(0)= P_+\eta_2(0)$.
\qed 
 
First, we prove the  G\^ateaux differentiability of $G_+^{\delta}$.
Let $\delta, l>0$ satisfying \eqref{ass-smallness}.
By Theorem \ref{thm-ex-csm}, for small $\varepsilon >0$ we have
\begin{align}\label{eq-c1-1}
\sup_{t \geq 0}\inf_{q\in \R} \norm{\tau_{\rho_0(t)}(v_h(t)+Q_{c_h(t)})-\tau_q Q_{c^*}}_{H^1} < \min\{ \varepsilon _*,\varepsilon _*^{1/2}, \delta\},
\end{align}
where $\varepsilon _*=\varepsilon _*(c^*,\delta,\delta^{-1/6}$ is defined in Lemma \ref{lem-exit}.
Since $|(v_0,\partial_xQ_{c^*})_{L^2}|+|(v_0,Q_{c^*})_{L^2}|\lesssim \varepsilon _*^{1/2}$, by the inequality \eqref{eq-c1-1} and Lemma \ref{lem-exit} we obtain 
\[ \mathfrak{m}_{\delta}((P_{\leq 0} v_0(t),c_0(t)),(P_{\leq 0}v_h(t),c_h(t))) > \delta^{-1/6}\norm{P_+(v_0(t)-v_h(t))}_{E}\]
for $t \geq 0$.
Therefore, we have
\begin{align}\label{eq-c1-2}
\frac{\norm{P_+(v_h(t)-v_0(t))}_{E}}{\norm{P_{\leq 0}(v_h(t)-v_0(t))}_{E_{\delta^{1/3}}}+h|a_1|} \leq \frac{\delta^{1/6}\mathfrak{m}_{\delta}((P_{\leq 0} v_0(t),c_0(t)),(P_{\leq 0}v_h(t),c_h(t)))}{\norm{P_{\leq 0}(v_h(t)-v_0(t))}_{E_{\delta^{1/3}}}+h|a_1|} \lesssim \delta^{-1/6}
\end{align}
for $t \geq 0$.
On the other hand, the convergence \eqref{eq-c1-0} yields 
\begin{align}\label{eq-c1-3}
\frac{\norm{(h_n')^{-1}P_+(v_{h_n'}(t)-v_0(t))}_E}{\norm{(h_n')^{-1}P_{\leq 0}(v_{h_n'}(t)-v_0(t))}_{E_{\delta^{1/3}}}+|(h_n')^{-1}(c_{h_n'}(t)-c_0(t))|} \to \frac{\norm{P_+\eta(t)}_E}{\norm{P_{\leq 0}\eta(t)}_{E_{\delta^{1/3}}}+|a_1|} 
\end{align}
as $n \to \infty$ for $t \geq 0$.
Since 
\[\frac{\norm{P_+\eta(t)}_E}{\norm{P_{\leq 0}\eta(t)}_{E_{\delta^{1/3}}}+|a_1|} \lesssim \delta^{-1/6}\]
by the inequality \eqref{eq-c1-2} and the convergence \eqref{eq-c1-3},  Lemma \ref{lem-prop-eta} yields that $\eta$ does not satisfy  $\eqref{ineq-eta-2}$ for some $t\geq 0$.
Thus, by the uniqueness in Lemma \ref{lem-eta-exis}, for any sequence $\{h_n\}_n$ with $h_n \to 0$ as $n \to \infty$ there exists subsequence $\{\tilde{h}_n\}_n \subset \{h_n\}_n$ such that
\[\frac{G_+^{\delta}(\psi_0+\tilde{h}_n\psi_1,a_0+\tilde{h}_na_1)-G_+^{\delta}(\psi_0,a_0)}{\tilde{h}_n} \to \psi_+ \mbox{ as } n \to \infty.\]
Therefore, we obtain the convergence
\begin{align}\label{eq-c1-4}
\frac{G_+^{\delta}(\psi_0+h\psi_1,a_0+ha_1)-G_+^{\delta}(\psi_0,a_0)}{h} \to \psi_+ \mbox{ as } h \to 0
\end{align}
which shows $G_+^{\delta}$ is G\^ateaux differentiable at $(\psi_0,a_0)$.
The linearity of the G\^ateaux derivative of $G_+^{\delta}$ follows the linearity of the solution to the equation \eqref{eq-eta} with respect to $(\psi,a_1)$.
The boundedness of the G\^ateaux derivative of $G_+^{\delta}$ follows the Lipschitz property of $G_+^{\delta}$.

Next we prove the continuity of the G\^ateaux derivative of $G_+^{\delta}$.
Let 
\[0<\varepsilon \ll \min\{\delta,\kappa_0(1,1),K_1(1,1)^{-6}\}\]
 and $\{(\psi_n,a_n)\}_{n=0}^{\infty} \subset P_{\leq 0}H^1(\RTL) \times (0,\infty)$ with $(\psi_n,a_n) \to (\psi_0,a_0)$ in $H^1(\RTL) \times (0,\infty)$ as $n \to \infty$ and $\sup_{n \in \N \cup \{0\}}(\norm{\psi_n}_{H^1}+|a_n-c^*|) < \tilde{\varepsilon}(c^*,\varepsilon )$, where $\tilde{\varepsilon}$ is defined in Theorem \ref{thm-ex-csm} and $\kappa_0(C,K_0)$ and $K_1(C,K_0)$ are defined in Lemma \ref{lem-prop-eta}.
Then, by Theorem \ref{thm-ex-csm} we have 
\[\sup_{t \geq 0, n \in \N \cup \{0\}}(\norm{v_n(t)}_{H^1}+|c_n(t)-c^*|) < \varepsilon,\]
where $(v_n,c_n,\rho_n)$ is the solution to the system \eqref{vZKeq} and \eqref{eq-orth-1} with $(v_n(0),c_n(0),\rho_n(0))=(\psi_n,a_n,0)$.
We define $\eta_n^{\psi,a}$ as the solution to the equation
\begin{align}\label{eq-eta-n}
\partial_t \eta=&\partial_x \mathbb{L}_{c^*}\eta -2\partial_x((Q_{c_n}-Q_{c^*})\eta)+(\dot{\rho}_n-c^*)\partial_x \eta-2a \partial_x(\partial_cQ_{c_n}v_n) \notag \\
&-a \partial_xQ_{c_n}+a(\dot{\rho}_n-c_n)\partial_x\partial_cQ_{c_n}-\dot{c}_na\partial_c^2Q_{c_n}-2\partial_x(v_n\eta)
\end{align}
with the initial data $\psi$.
By the convergence of $\{(\tau_{-\rho_n}v_n,c_n,\rho_n)\}_n$ local in time, for $T,C>0$ we obtain the convergence 
\begin{align}\label{eta-conv-1}
\norm{\eta_n^{\psi,a}-\eta_0^{\psi,a}}_{L^{\infty}((0,T),H^1)} \to 0
\end{align}
as $n \to \infty$ uniformly on $\{(\psi,a) \in H^1(\RTL)\times \R:\norm{\psi}_{H^1}+|a|<C\}$.
For $T>0$, by the boundedness of $\{ \norm{\tau_{-\rho_n}v_n}_{X^{1,b}_T}+\norm{c_n-c^*}_{L^{\infty}(0,T)}\}_n$, we have  the convergence
\begin{align}\label{eta-conv-2}
\sup_{n \geq 0} \norm{\eta_n^{\varphi,b}-\eta_n^{\psi,a}}_{L^{\infty}((0,T),H^1)} \to 0
\end{align}
 as $(\varphi,b) \to (\psi,a)$ in $H^1\times \R$.
Let the G\^ateaux derivative of $G_+^{\delta}$ at $(\psi_n,a_n)$ be $\partial G_+^{\delta,n}$.
Applying Lemma \ref{lem-prop-eta} to $\eta_0^{\psi+\partial G_+^{\delta,0}(\psi,a),a}$, we obtain there exists $C>0$ such that 
\begin{align}\label{eq-c1-5}
3\norm{P_+\eta_0^{\psi+\partial G_+^{\delta,0}(\psi,a),a}(t)}_{E} \leq & 3 \varepsilon^{1/3}\Bigl(\norm{P_{\leq 0} \eta_0^{\psi+\partial G_+^{\delta,0}(\psi,a),a}(t)}_{E_{\varepsilon^{1/3}}}+|a|\Bigr) \notag \\
\leq  & Ce^{K_1\varepsilon^{1/6} t}\varepsilon^{1/3}(\norm{\psi}_{E_{\varepsilon^{1/3}}}+|a|)
\end{align}
for $t>0$, $\psi \in P_{\leq 0}H^1(\RTL)$ and $a \in \R$ with $\norm{\psi}_{E}+|a|\leq 1$.
Applying Lemma \ref{lem-prop-eta} to $\eta_0^{\psi_+,0}$, we have 
\begin{align}\label{eq-c1-6}
3\norm{P_+\eta_0^{\psi_+,0}(t)}_{E}> e^{\frac{k_*t}{2}}\norm{\psi_+}_{E}+e^{\frac{k_*}{2}t}\varepsilon^{1/3}\norm{P_{\leq 0}\eta_0^{\psi_+,0}(t)}_{E_{\varepsilon ^{1/3}}}
\end{align}
for $t > 1/2$ and $\psi_+ \in P_+H^1(\RTL)\setminus \{0\}$.
By the inequalities \eqref{eq-c1-5} and \eqref{eq-c1-6}, we obtain 
\begin{align}
&3\norm{P_+\eta_0^{\psi+\partial G_+^{\delta,0}(\psi,a)+\psi_+,a}(t)}_{E}\notag \\
\geq& 3\norm{P_+\eta_0^{\psi_+,0}(t)}_{E}-3\norm{P_+\eta_0^{\psi+\partial G_+^{\delta,0}(\psi,a),a}(t)}_{E} \notag \\
> & e^{\frac{k_*t}{2}}\norm{\psi_+}_{E}-2e^{K_1\varepsilon^{1/6}t}C\varepsilon^{1/3}(\norm{\psi}_{E_{\varepsilon^{1/3}}}+|a|)+\varepsilon^{1/3}\Bigl(\norm{P_{\leq 0}\eta_0^{\psi+\partial G_+^{\delta,0}(\psi,a),a}(t)}_{E_{\varepsilon^{1/3}}}+|a|\Bigr)\notag \\
&+e^{\frac{k_*t}{2}}\varepsilon^{1/3}\norm{P_{\leq 0}\eta_0^{\psi_+,0}(t)}_{E_{\varepsilon^{1/3}}}
\end{align}
for $t > 1/2$, $\psi_+ \in P_+H^1(\RTL)\setminus \{0\}$, $\psi \in P_{\leq 0}H^1(\RTL)$ and $a \in \R$ with $\norm{\psi}_{E}+|a|\leq 1$.
By the convergence $\eta_n^{\varphi,b}(t) \to \eta_0^{\varphi,b}(t)$ as $n \to \infty$ in $L^{\infty}((0,T),H^1)$ uniformly on $(\varphi,b) \in H^1\times \R$ with $\norm{\varphi}_{E_{\varepsilon^{1/3}}}+|b|\leq 1$ for each $T>0$, we obtain there exists $n_T>0$ such that for $n \geq n_T$ and $1/2 < t\leq T$
\begin{align*}
6\norm{P_+\eta_n^{\psi+\partial G_+^{\delta,0}(\psi,a)+\psi_+,a}(t)}_{E} \geq& e^{\frac{k_*t}{2}}(\norm{\psi_+}_{E_{\varepsilon^{1/3}}}-2e^{-\frac{k_*t}{2}+K_1\varepsilon^{1/6}t}C\varepsilon^{1/3}(\norm{\psi}_{E_{\varepsilon^{1/3}}}+|a|))\notag \\
&+\varepsilon^{1/3}\Bigl(\norm{P_{\leq 0}\eta_n^{\psi+\partial G_+^{\delta,0}(\psi,a)+\psi_+,a}(t)}_{E_{\varepsilon^{1/3}}}+|a|\Bigr).
\end{align*}
Thus, for $\sigma>0$ there exists $T_{\sigma}>0$ such that $\eta_n^{\psi+\partial G_+^{\delta,0}(\psi,a)+\psi_+,a}$ satisfies  \eqref{cond-eta} at $T_{\sigma}$ for $n\geq n_{T_{\sigma}}$, $a \in \R$, $\psi \in P_{\leq 0} H^1(\RTL)$ and $\psi_+ \in P_+H^1(\RTL)$ with $\norm{\psi}_E+|a|\leq 1$ and $\norm{\psi_+}_E\geq \sigma$.
Since $\eta_n^{\psi+\partial G_+^{\delta,n}(\psi,a),a}$ do not satisfies \eqref{cond-eta} for some $t\geq 0$, we obtain for $n \geq n_{T_{\sigma}}$ 
\[\norm{\partial G_+^{\delta,n}(\psi,a)-\partial G_+^{\delta,0}(\psi,a)}_{E}<\sigma\]
which implies the continuity of the G\^ateaux derivative of $G_+^{\delta}$ at $(\psi_0,a_0)$ in the sense of the operator norm from $P_{\leq 0}H^1(\RTL) \times \R$ to $P_+H^1(\RTL)$.
Therefore, $G_+^{\delta}$ is $C^1$ class on $P_{\leq 0}H^1(\RTL)\times (0,\infty)$ in the sense of the Fr\'echet differential.
By the equation \eqref{eq-corr} in Corollary \ref{cor-corr}, we obtain the $C^1$ regularity of the manifold $\mathcal{M}_{cs}(c^*,\tilde{\varepsilon}(c^*,\varepsilon ))$.

Finally, we construct a forward flow invariant manifold containing $\mathcal{M}_{cs}(c^*,\varepsilon)$.
By Theorem \ref{thm-ex-csm}, Corollary \ref{cor-repulsive}, Corollary \ref{cor-Lip} and the $C^1$ regularity of the manifold $\mathcal{M}_{cs}(c^*,\varepsilon)$, there exist $l_0, \delta, \varepsilon_m>0$ such that the manifold $\mathcal{M}_{cs}(c^*,\varepsilon_m)$ satisfies (i), (iii), (iv) and (v) in Theorem \ref{thm-main} and 
\[U(t) \mathcal{M}_{cs}(c^*,\varepsilon_m) \subset \mathcal{M}_{cs}(c^*,\varepsilon^*(c^*,\delta,l_0))\]
for $t\geq 0$, where $\varepsilon_*$ is defined in Corollary \ref{cor-repulsive}.
We define
\[ \mathcal{M}_{cs}(c^*)=\{ u(t) \in H^1(\RTL); t\geq 0, u \mbox{ is the solution to } \eqref{ZKeq} \mbox{ with } u(0) \in \mathcal{M}_{cs}(c^*,\varepsilon_m)\}.\]
Then, $\mathcal{M}_{cs}(c^*)$ is a forward flow invariant manifold satisfying (i), (iii), (iv) and (v) in Theorem \ref{thm-main} and $\tau_q \mathcal{M}_{cs}(c^*) = \mathcal{M}_{cs}(c^*)$ for $q \in \R$.
Corollary \ref{cor-main} follows (iv) of Theorem \ref{thm-main} and (ii) of Theorem 1.5 in \cite{YY4}.




\section*{Acknowledgments}
The author would like to express his great appreciation to Professor Yoshio Tsutsumi for a lot to helpful advices and encouragements.
The author would like to thank Professor Nikolay Tzvetkov for his helpful encouragements.
The author is supported by JSPS Research Fellowships for Young Scientists 18J00947.






\end{document}